\documentclass[aap,seceqn,citesort,MSNbibl,dvips]{arximspdf}

\usepackage{graphicx}

\doi{10.1214/10-AAP681}
\volume{21}
\issue{2}
\pubyear{2011}
\firstpage{484}
\lastpage{545}

\makeatletter

\newcommand{\ve}{\varepsilon}
\newcommand{\ra}{\rightarrow}
\newcommand{\E}{\mathbb{E}}
\newcommand{\PP}{\mathbb{P}}
\newcommand{\R}{\mathbb{R}}
\newcommand{\N}{\mathbb{N}}
\newcommand{\ind}{\mathbb{I}}
\newcommand{\osc}{\operatorname{Osc}}
\newcommand{\sMM}{\mathcal{M}}
\newcommand{\sB}{\mathcal{B}}
\newcommand{\sW}{\mathcal{W}}
\newcommand{\sQ}{\mathcal{Q}}
\newcommand{\sU}{\mathcal{U}}

\newcommand{\Shat}{\widehat{S}}
\newcommand{\Ahat}{\widehat{A}}
\newcommand{\Vhat}{\widehat{V}}
\newcommand{\Nhat}{\widehat{N}}
\newcommand{\Ihat}{\widehat{I}}
\newcommand{\What}{\widehat{W}}
\newcommand{\Fhat}{\widehat{F}}
\newcommand{\Chat}{\widehat{C}}
\newcommand{\Qhat}{\widehat{Q}}
\newcommand{\sWhat}{\widehat{\mathcal{W}}}
\newcommand{\sUhat}{\widehat{\mathcal{U}}}
\newcommand{\Uhat}{\widehat{U}}

\newcommand{\hsigman}{\widehat{\sigma}{}^{(n)}}
\newcommand{\htaun}{\widehat{\tau}{}^{(n)}}
\newcommand{\hyn}{\widehat{Y}^{(n)}}
\newcommand{\ven}{\ve^{(n)}}
\newcommand{\busy}{[\hsigman_k,\htaun_k)}
\newcommand{\idle}{[\htaun_{k-1},\hsigman_k)}
\newcommand{\sma}{\Delta}
\newcommand{\newsmall}{\delta}

\newcommand{\kpn}{K_+^{(n)}}
\newcommand{\knn}{K_-^{(n)}}
\newcommand{\hkpn}{\widehat{K}_+^{(n)}}

\newcommand{\renu}{R^{(n)}_U}
\newcommand{\renw}{R^{(n)}_W}
\newcommand{\hrenu}{\widehat{R}^{(n)}_U}
\newcommand{\hrenw}{\widehat{R}^{(n)}_W}
\newcommand{\hrenq}{\widehat{R}^{(n)}_Q}
\newcommand{\hrew}{R_W^*}
\newcommand{\sS}{\mathcal{S}}
\newcommand{\Khat}{\widehat{K}}
\newcommand{\Dbar}{\overline{D}}

\newcommand{\define}{\stackrel{\Delta}{=}}
\newcommand{\inprob}{\stackrel{P}{\longrightarrow}}

\newtheorem{theorem}{Theorem}[section]
\newtheorem{lemma}[theorem]{Lemma}
\newtheorem{proposition}[theorem]{Proposition}
\newproclaim{example}[theorem]{Example}
\newtheorem{corollary}[theorem]{Corollary}
\newproclaim{remark}[theorem]{Remark}

\makeatother

\begin{document}
\begin{frontmatter}

\title{Heavy traffic analysis for EDF queues with reneging}
\runtitle{EDF queues with reneging}

\begin{aug}
\author[A]{\fnms{\L ukasz} \snm{Kruk}\thanksref{t1}\ead
[label=e1]{lkruk@hektor.umcs.lublin.pl}},
\author[B]{\fnms{John} \snm{Lehoczky}\thanksref{t2}\ead
[label=e2]{jpl@stat.cmu.edu}},
\author[C]{\fnms{Kavita} \snm{Ramanan}\thanksref{t3}\ead
[label=e3]{Kavita\_Ramanan@brown.edu}}\\ and
\author[D]{\fnms{Steven} \snm{Shreve}\corref{}\thanksref{t4}\ead
[label=e4]{shreve@andrew.cmu.edu}}
\runauthor{Kruk, Lehoczky, Ramanan and Shreve}
\affiliation{Maria Curie-Sklodowska University and Polish Academy of Sciences,
Carnegie Mellon University, Brown University and Carnegie Mellon University}
\address[A]{\L. Kruk\\
Department of Mathematics\\
Maria Curie-Sklodowska University \\
Lublin\\
Poland \\
and\\
Institute of Mathematics\\
Polish Academy of Sciences\\
Warsaw\\
Poland\\
\printead{e1}} 
\address[B]{J. Lehoczky\\
Department of Statistics\\
Carnegie Mellon University\\
Pittsburgh, Pennsylvania 15213\hspace*{27.1pt}\\
USA\\
\printead{e2}}
\address[C]{K. Ramanan\\
Division of Applied Mathematics\\
Brown University\\
Providence, Rhode Island 02912\\
USA\\
\printead{e3}}
\address[D]{S. Shreve\\
Department of Mathematical Sciences\\
Carnegie Mellon University\\
Pittsburgh, Pennsylvania 15213\\
USA\\
\printead{e4}}
\end{aug}

\thankstext{t1}{Supported in part by the State
Committee for Scientific Research of Poland, Grant 2 P03A 012 23
and the EC FP6 Marie Curie ToK programme SPADE 2 at IMPAN,
Poland.}

\thankstext{t2}{Supported in part by ONR and DARPA under MURI
Contract N00014-01-1-0576.}

\thankstext{t3}{Supported in part by the
NSF Grants DMS-04-06191 and DMS-04-05343 and CMMI-1059967 (formerly CMMI-0728064).}

\thankstext{t4}{Supported in part by the NSF Grants DMS-04-04682
and DMS-09-03475.}

\received{\smonth{12} \syear{2007}}
\revised{\smonth{8} \syear{2009}}

%
\begin{abstract}
This paper presents a heavy-traffic analysis of the behavior of a
single-server queue under an Earliest-Deadline-First (EDF) scheduling
policy in which customers have deadlines and are served only until
their deadlines elapse. The performance of the system is measured by
the fraction of reneged work (the residual work lost due to elapsed
deadlines) which is shown to be minimized by the EDF policy. The
evolution of the lead time distribution of customers in queue is
described by a measure-valued process. The heavy traffic limit of this
(properly scaled) process is shown to be a deterministic function of
the limit of the scaled workload process which, in turn, is identified
to be a doubly reflected Brownian motion. This paper complements
previous work by Doytchinov, Lehoczky and Shreve on the EDF discipline
in which customers are served to completion even after their deadlines
elapse. The fraction of reneged work in a heavily loaded system and
the fraction of late work in the corresponding system without reneging
are compared using explicit formulas based on the heavy traffic
approximations. The formulas
are validated by simulation results.
\end{abstract}

%
\begin{keyword}[class=AMS]
\kwd[Primary ]{60K25}
\kwd[; secondary ]{60G57}
\kwd{60J65}
\kwd{68M20}.
\end{keyword}
\begin{keyword}
\kwd{Due dates}
\kwd{heavy traffic}
\kwd{queueing}
\kwd{reneging}
\kwd{diffusion limits}
\kwd{random measures}
\kwd{real-time queues}.
\end{keyword}

\end{frontmatter}

\section{Introduction} \label{S.Intro}

\subsection{Background and the reneging EDF model}\label{S1.1}
In the last decade, attention has been paid to queueing
systems in which customers have deadlines. Examples
include telecommunication systems carrying digitized
voice or video traffic, tracking systems and real-time control
systems. In the case of voice or video, packetized
information must be received, processed and displayed within
stringent timing bounds so that the integrity of the
transmission is maintained. Similarly, there are processing
requirements for tracking systems that guarantee that a track can be
successfully followed. Real-time control systems
(e.g., those
associated with modern avionics systems, manufacturing plants or
automobiles) also gather data that must be
processed within stringent
timing requirements in order for the system to
maintain stability or react to
changes in the operating environment. We refer to queueing systems
that process tasks with deadlines as ``real-time queueing systems.''

The performance of a real-time queueing system is measured by its
ability to meet the deadlines of the customers. This is in contrast
to ordinary queueing systems in which the measure of performance is
often customer delay, queue length or
utilization of a service facility. We use the fraction
of ``reneged work,''
defined as the residual work not
serviced due to elapsed deadlines,
as our performance measure.
To minimize this quantity, it is
necessary to use a scheduling policy that takes
deadlines into account.
We use the Earliest-Deadline-First (EDF)
policy, which reduces to the more
familiar First-In-First-Out (FIFO) policy
when all customers have the same deadline.
Under general assumptions,
we prove that EDF is optimal with respect to this
performance measure. A related result for
$G/M/c$ queues, in which the number
of reneging customers is used as a performance measure,
was obtained by Panwar and Towsley \cite{pat}.

Heavy traffic analysis of a single real-time queue
was initiated by Lehocz\-ky~\cite{jpl96}.
This was put on a firm mathematical foundation
by Doytchinov, Lehoczky and Shreve (DLS) \cite{bogdan}.
The accuracy of heavy traffic approximations
was developed in
\cite{kls,kls1}. The results of DLS
were generalized to the case of
acyclic networks in
\cite{acyclic}.
In these papers it was assumed that all
customers are served to completion. The case
in which late customers
leave the system
and their residual work is lost is addressed here.
The main result of this paper is a heavy traffic
convergence theorem, from which is derived a simple and
practically useful approximation for the fraction of lost
work when the system is heavily loaded.

The mathematical formulation used by DLS and related papers is
based on random measures.
In addition to the usual queue length and
workload processes associated with the queueing system, to model the
evolution of a real-time queueing system, one
must keep track of the lead time of each customer, that is, the time
until the customer's deadline elapses. This is done by
measure-valued queue length and workload processes. The
measure-valued queue length process puts unit mass on the real line
at the lead time of each customer in the system, while the
measure-valued workload
process puts mass equal to the
remaining service time of each~customer at the lead time of that
customer. These measures evolve dynamically as customers
arrive, age
and depart.
Under the usual heavy traffic assumptions,
since customers are served to completion in the DLS framework,
it is easy to see that the ordinary scaled workload
process
converges weakly to a reflected Brownian motion with drift.
DLS showed that the suitably scaled workload and queue length
measure-valued processes converge to an explicit deterministic
function of the workload process.

In this paper customers
leave the system when their deadlines elapse,
which we refer to as reneging.
Due to the
preemptive nature of the EDF policy,
it is not possible to determine at the time
of admission whether a customer will be fully serviced
before its deadline elapses. It is thus natural to have the
controller make the decision only at the time
the deadline elapses.
The system with reneging shows
marked improvement in performance
over the DLS system, in the sense that the
fraction of reneged work
in this system is much less than the fraction
of work that becomes late in the DLS system.
This improvement is because
once a customer misses its deadline, the processor
devotes no further effort to it, but rather
turns its attention to customers that are not late.

The system with reneging is
considerably more difficult to analyze than the DLS system.
In the reneging system, the evolution
of the scalar total workload process
depends on the entire lead time distribution of customers in queue
and the nature of the EDF discipline.
This is in stark contrast to the DLS system, where the
total workload process is independent of the scheduling
discipline, and is identical to that of any
$GI/G/1$ queue with a work-conserving
scheduling discipline.
A key ingredient of our analysis is
a mapping on the space of measure-valued functions
which, when applied to the DLS system, yields another system
(that we call the reference system) whose difference from
the reneging system vanishes in heavy traffic.
This mapping can be viewed as a generalization of the
scalar double reflection map to measure-valued processes, and, using
its continuity properties, we identify the
heavy traffic limit of the reference and hence the reneging
systems.
Specifically, we show that the limit of the scaled workload process
is a doubly reflected Brownian motion with lower barrier zero
and upper barrier at the mean of the lead time distribution.
We also show that, conditional on the limiting workload, the
resulting limiting measure-valued
workload process is the same limiting process as when
customers are served to completion, that is, in the DLS system.
However, the workload processes in these two systems differ,
and so the unconditional limiting lead-time profiles of these
two systems differ accordingly. In particular, unlike in the
DLS system, the measure-valued workload process in the reneging
system is always concentrated on the positive
real line due to the absence of late work in the
reneging system.

\subsection{Prediction formulas}\label{S1.2}
The results of this paper suggest a simple formula for the fraction
of lost work in the EDF system with reneging.
In particular, consider
a single-server queue with traffic intensity
$\rho=\lambda/\mu$ that is near one,
where $1/\lambda$ is the mean
interarrival time and $1/\mu$ is the mean service time.
Let $\alpha$ and $\beta$ be the standard deviations of
the interarrival times and service times, respectively,
and set $\sigma^2=\lambda(\alpha^2+\beta^2)$,
which we assume is nonzero.
Let $\Dbar$ denote the mean lead time for
arriving customers. Finally, set
$\theta=2(1-\rho)/\sigma^2$.
Under these circumstances,
%
\begin{equation}\label{frw}
\mbox{Fraction of lost work in reneging system}\approx
e^{-\theta\Dbar}\biggl(\frac{1-\rho}{\rho(1-e^{-\theta\Dbar})}\biggr).
\end{equation}
This formula is derived in Section \ref{Subsection7.3}
and compared with simulations in Section~\ref{Subsection7.2}.
If $\rho=1$, in place of
(\ref{frw}) we have
%
\begin{equation}\label{frwo}
\mbox{Fraction of lost work in reneging system}\approx
\frac{\sigma^2}{2\Dbar}.
\end{equation}
Analysis of the
limit of the standard (nonreneging) system suggests
that when $\rho<1$
[see (\ref{lwcs2}) and (\ref{Ratio})],
%
\begin{equation}\label{flw}
\mbox{Fraction of late work in standard system}
\approx
e^{-\theta\Dbar},
\end{equation}
which, together with (\ref{frw}), yields the approximation
%
\begin{equation}\label{7.3}
\frac{\mbox{Lost work in reneging system}}
{\mbox{Late work in standard system}}
\approx
\frac{1-\rho}{\rho(1-e^{-\theta\Dbar})}.
\end{equation}
If $\rho\geq1$ then all work is late
in the limiting standard (nonreneging) system, which leads to the
approximation
%
\begin{equation}\label{Ratioo}
\frac{\mbox{Lost work in reneging system}}
{\mbox{Late work in standard system}}
\approx\frac{\sigma^2}{2\Dbar}.
\end{equation}
When plotted\vspace*{1pt} on a log scale,
the fraction of lost work
in the reneging system
and the fraction of late work in the standard system
will be linear in $\Dbar$, provided
that $e^{\theta\Dbar}\gg 1$,
and these two plots will be separated by
$\log((1-\rho)/\rho)$.
When performance is measured in terms of the work whose service
requirement is not met by the time its deadline elapses, then
the reneging system is far
superior to the nonreneging system.
We refer the reader to the simulations
in Section \ref{Subsection7.2}.

The situation with reneged \textit{customers} as opposed to
reneged work is more complicated.
DLS shows that the number of customers
in the limiting standard system at any time is just $\lambda$ times
the amount of work, the number of late customers is
$\lambda$ times the amount of late work
and hence
%
\begin{eqnarray}\label{lwcs}
&&\mbox{Fraction of late customers in the standard system}
\nonumber\\[-8pt]\\[-8pt]
&&\qquad\approx
\mbox{Fraction of late work in the standard system}\nonumber
\end{eqnarray}
[see also (\ref{lwcs2}) and its derivation for the case $\rho<1$].
In the limiting reneging system, the number of customers who arrive
by a certain time
and the number of customers in system at that time
is $\lambda$ times
the amount of arrived work and $\lambda$ times the amount
of work in the system (Corollary \ref{C3.5a}), respectively, but the
number of customers who renege by a certain time is not
necessarily $\lambda$ times the amount of reneged work
by that time (see Remark \ref{RQ}). In particular,
we do not have a formula like (\ref{lwcs}) for the reneging
system. If
the arrival process is Poisson, the
fraction of lost customers in the reneging system
can be estimated by a heuristic argument
[see (\ref{7.6})] which gives instead
%
\begin{eqnarray} \label{flcr}
&&\mbox{Fraction of lost customers in reneging system}
\nonumber\\[-8pt]\\[-8pt]
&&\qquad\approx
\frac{2}{\mu^2\beta^2+1}\times
(\mbox{Fraction of lost work in reneging system}).\nonumber
\end{eqnarray}

\subsection{Related work and outline of paper}
\label{S1.3}
Measure-valued processes have recently gained prominence
in queueing theory.
Decreusefond and Moyal \cite{DM} use such processes
to obtain the fluid limit of an EDF $M/M/1$ queue with
reneging. Unlike our scaling (\ref{2.4}) of lead
times by $\sqrt{n}$, they scale lead times by $n$
and obtain a characterization of the limiting
lead-time measure-valued process via a transport
equation.
In a different setting, Ward and Glynn \cite{WaGl,WaGl2} find
limits of FIFO queues with reneging. Measure-valued processes
have also proved
useful in the heavy traffic analysis of queues
with scheduling disciplines
other than EDF such as last-in-first-out \cite{lim00}, processor
sharing \mbox{\cite{gro04,grokru07}}, and shortest remaining processing time
\mbox{\cite{DGP09,GKP}}. As dynamical systems, queueing systems present
a mathematical challenge due to discontinuities in their evolution at
boundaries (which denote empty queues). The heavy traffic analysis of
queueing systems described by $\R^n$-valued processes has been
facilitated by the use of representations in terms of continuous
mappings on $\R^n$ \cite{cheman,dupram2,harrei,ramrei1,Whitt}. This
work demonstrates that this perspective can also be useful when the
queueing system is represented by a more complicated, measure-valued
process (see also \cite{kasram1} for recent work that takes a similar
perspective).

Section \ref{S2} introduces our model. Section \ref{S3} summarizes the main
results, and proofs of these results are given in
Section \ref{S5}.
Section \ref{S4} introduces the reference workload process and its
decomposition, and describes its evolution.
This reference workload process is easier to analyze
than the workload process with reneging but the two
are shown to have
the same asymptotic behavior. Comparisons between
the reference workload process and the reneging workload process
are presented in Section \ref{S5a}.
Section \ref{S.sim}
presents simulation results. A proof of optimality of EDF, that may
be of independent interest, is in the \hyperref[app]{Appendix}.

\section{The model, assumptions and notation}\label{S2}

\subsection{Notation}\label{S2.0}

Let $\R$ be the set of real numbers.
For $a,b\in\R$, $a\vee b$ is the maximum of $a$ and $b$,
$a\wedge b$ is the minimum
and $a^+$ is the maximum
of $a$ and $0$.
Also, $\inf\{\varnothing\}$ should be understood as $+\infty$, while
$\sup\{\varnothing\}$ and $\max\{\varnothing\}$
should be understood as $-\infty$. Moreover, if $a < b$, then
the interval $[b,a]$ is understood to be~$\varnothing$.

Denote by $\sMM$ the set of all finite, nonnegative measures
on $\sB(\R)$, the Borel subsets of $\R$.
Under the weak topology, $\sMM$ is a Polish space.
We denote the measure in $\sMM$ that puts one unit of mass
at the point $x\in\R$, that is, the Dirac measure at
$x$, by $\delta_x$.
When $\nu\in\mathcal{M}$ and
$B$ is an interval
$(a,b]$ or a singleton $\{a\}$,
we will simply write $\nu(a,b]$ and $\nu\{a\}$
instead of $\nu((a,b])$ and $\nu(\{a\})$.

Let $T>0$ be given.
Given a Polish space $X$, we use $D_X[0,\infty)$
(resp., $D_X[0,T]$) to
denote the space of
right-continuous functions with left-hand limits
(RCLL functions) from $[0,\infty)$ (resp., $[0,T]$)
to $X$, equipped with the Skorokhod $J_1$ topology.
See \cite{ek} for details.
When dealing with $D_X[0,\infty)$ or $D_X[0,T]$, we
typically consider $X=\R$ or $\R^d$,
with appropriate dimension $d$ for vector-valued functions, or $X =
\sMM$,
unless explicitly\vspace*{1pt} stated otherwise. When $X=\R$ or $\sMM$,
for $t>0$ and $x\in D_X[0,\infty)$, we write $x(t-)$
for the\vspace*{-1pt} left-hand limit $\lim_{s\uparrow t} x(s)$, and
we define $\bigtriangleup x(t)$ to be the jump in $x$ at time $t$, that is,
$\bigtriangleup x(t)\define x(t)-x(t-)$.
Finally, given $D_X[0,\infty)$-valued random variables
$Z_n, n \in\N$,
defined, respectively, on the probability spaces
$(\Omega_n, \mathcal{F}_n, \PP_n)$, $n \in\N$,
and a $D_X[0,\infty)$-valued random variable $Z$
defined on a probability space
$(\Omega, \mathcal{F}, \PP)$, we say $Z^{(n)}$ converges
in distribution to~$Z$ and write
$Z_n \Rightarrow Z$, if
for every bounded continuous function $f$ on $D_X[0,\infty)$,
$\lim_{n \ra\infty} \E_n[ f(Z_n)] = \E[ f(Z)]$.
Here $\E_n$ and $\E$ are expectations taken with
respect to $\PP_n$ and $\PP$, respectively.

\subsection{The model with reneging}\label{S2.1}

We have a sequence of single-station queueing
systems, each serving one class of customers. The queueing
systems are indexed by superscript $(n)$.
The inter-arrival times for the customers are
$\{u_j^{(n)}\}_{j=1}^{\infty}$,
a sequence of strictly positive,
independent, identically distributed
random variables with common mean
$\frac{1}{\lambda^{(n)}}$ and standard deviation
$\alpha^{(n)}$. The service times
are $\{v_j^{(n)}\}_{j=1}^{\infty}$,
another sequence of positive, independent,
identically distributed random variables with
common mean
$\frac{1}{\mu^{(n)}}$
and standard deviation $\beta^{(n)}$.

If the initial condition of the $n$th queue were not zero, then
we would need to specify an initial workload measure-valued
process and frontier
[these terms are defined in (\ref{mvwp}) and (\ref{frontier}) below]
in such a way that these have limits under the heavy traffic
scaling. However, if the limit of the initial scaled
workload process were not of the form appearing in
Theorem \ref{T2.2S} below,
then the workload process would be expected to have a jump
at time zero. To avoid these complications, we assume that
each queue is empty at time zero.

We define the \textit{customer arrival times}
%
\begin{equation}\label{2.1}
S^{(n)}_{0} \define0,\qquad S^{(n)}_k \define\sum_{i=1}^k
u_{i}^{(n)},\qquad
k \geq1,
\end{equation}
the \textit{customer arrival process}
%
\begin{equation}\label{2.2}
A^{(n)}(t)\define\max\bigl\{k;S^{(n)}_k\leq t\bigr\},\qquad
t\geq0,
\end{equation}
and the \textit{work arrival process}
%
\begin{equation}\label{2.3}
V^{(n)}(t)\define\sum_{j=1}^{\lfloor t\rfloor}
v_j^{(n)},\qquad t\geq0.
\end{equation}
The work that has arrived to the queue by time
$t$ is then
$V^{(n)}(A^{(n)}(t))$.

Each customer\vspace*{1pt} arrives with an initial lead time
$L_j^{(n)}$, the time between the arrival time and
the deadline for completion of service for
that customer.
These initial lead times are independent
and identically distributed with
%
\begin{equation}\label{2.4}
\PP\bigl\{L_j^{(n)}\leq\sqrt{n}y\bigr\}=G(y),
\end{equation}
where $G$ is a right-continuous cumulative
distribution function.
We define
%
\begin{equation}\label{2.5}
y_*\define\inf\{y\in\R| G(y)>0\},\qquad
y^*\define\min\{y\in\R| G(y)=1\}
\end{equation}
and assume that $0<y_*\leq y^*<+\infty$.
We assume that for every $n$, the sequences
$\{u_j^{(n)}\}_{j=1}^\infty$, $\{v_j^{(n)}\}_{j=1}^\infty$
and $\{L_j^{(n)}\}_{j=1}^\infty$ are mutually independent.
See\vspace*{2pt} Remark \ref{rem-ystar} for a discussion of these assumptions.

We assume that customers
are served using
the Earliest-Deadline-First (EDF) queue discipline,
that is, the customer with
the shortest lead time receives service.
Preemption occurs when a customer more urgent than the customer in
service arrives
(we assume preempt-resume). There is no set up,
switch-over, or other type of overhead. If the $j$th customer
is still present in the system (either waiting for service or receiving
it) when his deadline passes, that is, at the time
$S^{(n)}_j+L^{(n)}_j$, he leaves the queue immediately. This
may be interpreted as
either reneging or the result of an action of an external controller.

We define $W^{(n)}(t)$, the \textit{workload process} at time $t$, as the remaining
processing time of all the customers
in the system at this time.
We define $\renw(t)$ to be the amount of work that reneges
in the time interval $[0,t]$. The
\textit{queue length process} $Q^{(n)}(t)$ is
the number of customers in the queue at time $t$.
The queueing system described above will be referred to as the \textit{EDF
system with reneging}.

\subsection{The standard EDF model}

We also have a sequence, indexed by superscript
$(n)$, of \textit{standard EDF systems},
with the same stochastic
primitives as the EDF systems with reneging.
In each of these standard systems, the server
serves the customer with
the shortest lead time, preemption occurs
as in the reneging system,
but late customers (customers with negative
lead times) stay in the system until served to
completion.
The
performance processes associated with the
standard system will be denoted by the same symbols as their
counterparts from the system with reneging, but with additional
subscript $S$. For example, $W^{(n)}_S(t)$ denotes the workload in the
standard system at time $t$. The arrival
processes $A^{(n)}(t)$ and $V^{(n)}(t)$ are the same for the both
systems, so we will not attach the subscript $S$ to them.

The standard EDF system is easier to analyze than the EDF system with
reneging in several ways.
For instance, the workload $W^{(n)}_S$ in the standard
system coincides with the workload of a corresponding G/G/1 queue (with
the same primitives) under
any nonidling scheduling policy.
More precisely, in the standard system
the \textit{netput process}
%
\begin{equation}\label{2.15}
N^{(n)}(t)\define V^{(n)}\bigl(A^{(n)}(t)\bigr)-t
\end{equation}
measures the amount of work in queue at time
$t$ provided that the server is never idle up to time $t$, and
the \textit{cumulative idleness process}
%
\begin{equation}\label{2.16}
I^{(n)}_S(t)\define-\inf_{0\leq s\leq t}N^{(n)}(s)
\end{equation}
gives the amount of time the server
is idle. Adding these two processes together, we obtain the
workload process for the standard system
%
\begin{equation}\label{2.18}
W^{(n)}_S(t)=N^{(n)}(t)+I^{(n)}_S(t).
\end{equation}
(All the above processes are RCLL.)
In contrast, the evolution of the
workload $W^{(n)}$ in the reneging system is
more complex and depends
not only on the residual service times but
also on the lead times
of all customers in the queue.
Our analysis of the reneging system will
be facilitated by
results from \cite{bogdan} on the heavy
traffic analysis of the standard EDF system.

\subsection{Heavy traffic assumptions}\label{S2.2}
We assume that the following
limits exist:
%
\begin{eqnarray}\label{3.1}
\lim_{n\rightarrow\infty}\lambda^{(n)}&=&\lambda,\qquad
\lim_{n\rightarrow\infty}\mu^{(n)}=\lambda,\nonumber\\[-8pt]\\[-8pt]
\lim_{n\rightarrow\infty}\alpha^{(n)}&=&\alpha,\qquad
\lim_{n\rightarrow\infty}\beta^{(n)}=\beta,\nonumber
\end{eqnarray}
and, moreover, $\lambda>0$ and $\alpha^2+\beta^2>0$.
Define the \textit{traffic intensity} $\rho^{(n)}\define\frac{\lambda
^{(n)}}{\mu^{(n)}}$.
We make the \textit{heavy traffic assumption}
%
\begin{equation}\label{3.3}
\lim_{n\rightarrow\infty}\sqrt{n}
\bigl(1-\rho^{(n)}\bigr)=\gamma
\end{equation}
for some $\gamma\in\R$.
We also impose the \textit{Lindeberg
condition} on the inter-arrival and service times: for every $c>0$,
%
\begin{eqnarray}\label{3.5}
&&\lim_{n\rightarrow\infty}
\E\bigl[\bigl(u_j^{(n)}-\bigl(\lambda^{(n)}\bigr)^{-1}\bigr)^2
\ind_{\{|u_j^{(n)}-(\lambda^{(n)})^{-1}|
>c\sqrt{n}\}} \bigr]\nonumber\\[-8pt]\\[-8pt]
&&\qquad=
\lim_{n\rightarrow\infty}
\E\bigl[\bigl(v_j^{(n)}-\bigl(\mu^{(n)}\bigr)^{-1}\bigr)^2
\ind_{\{|v_j^{(n)}-(\mu^{(n)})^{-1}|
>c\sqrt{n}\}}\bigr]= 0. \nonumber
\end{eqnarray}

We introduce the \textit{heavy traffic scaling}
for the idleness process in the standard system
and the workload and queue length processes for both EDF systems
%
\begin{eqnarray*}
\Ihat_S^{(n)}(t) &=&
\frac{1}{\sqrt{n}}I_S^{(n)}(nt),\qquad
\What_S^{(n)}(t)=
\frac{1}{\sqrt{n}}W^{(n)}_S(nt),\\
\Qhat_S^{(n)}(t) &=& \frac{1}{\sqrt{n}}Q_S^{(n)}(nt), \qquad
\What^{(n)}(t) =
\frac{1}{\sqrt{n}}W^{(n)}(nt),\\
\Qhat^{(n)}(t) &=& \frac{1}{\sqrt{n}}Q^{(n)}(nt)
\end{eqnarray*}
and the \textit{centered heavy traffic scaling}
for the arrival processes
\begin{eqnarray*}
\Shat{}^{(n)}(t)
&=&
\frac{1}{\sqrt{n}}\sum_{j=1}^{\lfloor nt\rfloor}
\biggl(u_j^{(n)}
-\frac{1}{\lambda^{(n)}}\biggr),\qquad
\Vhat^{(n)}(t)=
\frac{1}{\sqrt{n}}\sum_{j=1}^{\lfloor nt\rfloor}
\biggl(v_j^{(n)}
-\frac{1}{\mu^{(n)}}\biggr),
\\
\Ahat^{(n)}(t) &=&
\frac{1}{\sqrt{n}}
\bigl[A^{(n)}(nt)-\lambda^{(n)} nt\bigr].
\end{eqnarray*}
The scaled netput process (which is the same for both systems) is given by
%
\begin{equation}\label{3.10}
\Nhat^{(n)}(t)
=\frac{1}{\sqrt{n}}
\bigl[V^{(n)}\bigl(A^{(n)}(nt)\bigr)- nt\bigr].
\end{equation}
Note that, by (\ref{2.18}), $\What_S^{(n)}(t)=\Nhat^{(n)}(t)+\Ihat
_S^{(n)}(t)$.

It follows from Theorem 3.1 in \cite{pr} and Theorem 7.3.2 in
\cite{Whitt} that
%
\begin{equation}\label{An}
\bigl(\Shat^{(n)},\widehat{A}^{(n)}\bigr)
\Rightarrow(S^*,A^*),
\end{equation}
where $A^*$
is a zero-drift Brownian motion with
variance $\alpha^2 \lambda^3$ per unit time and
%
\begin{equation}\label{SArel}
S^*(\lambda t)=-\frac{1}{\lambda}A^*(t),\qquad t\geq0.
\end{equation}
It is a standard result \cite{IW} that
%
\begin{equation}\label{w*}
\bigl(\Nhat^{(n)},\Ihat_S^{(n)},\What_S^{(n)}\bigr) \Rightarrow
(N^*,I^*_S,W^*_S),
\end{equation}
where $N^*$ is a Brownian motion with variance
$(\alpha^2 +\beta^2)\lambda$ per unit time
and drift $-\gamma$,
%
\begin{equation}\label{2.16c}
I^*_S(t)\define-\min_{0\leq s\leq t}N^*(s),\qquad
W_S^*(t)=N^*(t)+I_S^*(t).
\end{equation}
In other words,
$W_S^*$ is a Brownian motion reflected at $0$ with
variance $(\alpha^2+\beta^2)\lambda$ per unit time
and drift $-\gamma$ and $I_S^*$ causes the reflection.

\subsection{Measure-valued processes and frontiers}\label{S2.3}
To study whether tasks or customers meet their
timing requirements, one must keep track of customer
lead times.
The action of the EDF discipline requires knowledge
of the current lead times
of all customers in system.
We represent this information
via a collection of measure-valued stochastic processes.

\textit{Customer arrival measure-valued process}:
%
\[
\mathcal{A}^{(n)}(t)(B) \define\left\{
\begin{array}{l}
\mbox{Number of
arrivals by time $t$, whether or not still in
the}\\
\mbox{system at time $t$,
having lead times at time }
t\mbox{ in } B \in\sB(\R)
\end{array}
\right\}.
\]

\textit{Workload arrival measure-valued process}:
%
\[
\mathcal{V}^{(n)}(t)(B)\define\left\{
\begin{array}{l}
\mbox{Work arrived by time $t$,
whether or not still in the sys-}\\
\mbox{tem at time $t$, having lead times at time }t
\mbox{ in } B \in\sB(\R)
\end{array}
\right\}.
\]

\textit{Queue length measure-valued process}:
\[
\mathcal{Q}^{(n)}(t)(B) \define
\left\{
\begin{array}{l}
\mbox{Number of customers in
the queue at time}\\
t \mbox{ having lead times at time }t\mbox{ in } B \in\sB(\R)
\end{array}
\right\}.
\]

\textit{Workload measure-valued process}:
%
\begin{equation}\label{mvwp}\qquad
\mathcal{W}^{(n)}(t)(B) \define
\left\{
\begin{array}{l}
\mbox{Work in the queue at time }t
\mbox{ associated with cus-}\\
\mbox{tomers having lead times at time }
t\mbox{ in } B \in\sB(\R)
\end{array}
\right\}.
\end{equation}
The latter two processes describe the behavior of the EDF system with
reneging. Their counterparts for the standard EDF
system will be denoted by $\mathcal{Q}^{(n)}_S(t)$ and $\mathcal
{W}^{(n)}_S(t)$, respectively.
The following relationships easily follow:
%
\begin{eqnarray*}
A^{(n)}(t) & = & \mathcal{A}^{(n)}(t)(\R),\qquad
V^{(n)}\bigl(A^{(n)}(t)\bigr)  =  \mathcal{V}^{(n)}(t)(\R),
\\
W^{(n)}(t) & = & \mathcal{W}^{(n)}(t)(0,\infty),\qquad
Q^{(n)}(t)  =  \mathcal{Q}^{(n)}(t)(0,\infty), \\
W_S^{(n)}(t) & = & \mathcal{W}_S^{(n)}(t)(\R),\qquad
Q_S^{(n)}(t)  =  \mathcal{Q}_S^{(n)}(t)(\R).
\end{eqnarray*}
In addition, we can represent the reneged work
as follows:
%
\begin{equation}
\label{def-renw}
\renw(t) = \sum_{0 < s \leq t} \mathcal{W}^{(n)} (s-) \{0\}.
\end{equation}

In order to study the behavior of the EDF queue discipline, it is
useful to keep track of
the largest lead time of all customers,
whether present or departed, who have
ever been in service. We define
the \textit{frontier}
%
\begin{equation}\label{frontier}
F^{(n)}(t) \define
\left\{
\begin{array}{l}
\mbox{The maximum of the largest lead time of}\\
\mbox{all customers who have ever been in service,}\\
\mbox{whether still present or not, and } \sqrt{n} y^*-t
\end{array}
\right\}
\end{equation}
for the EDF system with reneging, and its counterpart $F_S^{(n)}(t)$
for the standard EDF system.
Prior to arrival of the first customer, $F^{(n)}(t)$ and $F^{(n)}_S(t)$
equal $\sqrt{n} y^*-t$.
For the EDF system with reneging, we define
the \textit{current lead time}
%
\[
C^{(n)}(t) \define
\left\{
\begin{array}{l}
\mbox{Lead time of the customer in service}\\
\mbox{or }
F^{(n)}(t)\mbox{ if the queue is empty}
\end{array}
\right\}.
\]
In the reneging system, there is no customer 
with lead time smaller than
$C^{(n)}(t)$, and there has never been a customer
in service whose lead time,
if the customer were still present, would
exceed $F^{(n)}(t)$. Furthermore,
$C^{(n)}(t)\leq F^{(n)}(t)$ for all $t\geq0$.
The processes $C^{(n)}$, $F^{(n)}$ and $F_S^{(n)}$
are RCLL.

We introduce heavy traffic scalings. For the real-valued
processes $Z^{(n)} = C^{(n)}, F^{(n)}, F_S^{(n)}, W^{(n)},
Q^{(n)}, R_W^{(n)}$, we define
$\widehat{Z}^{(n)} (t) \define\frac{1}{\sqrt{n}}
Z^{(n)} (nt)$
and for the measure-valued processes
$\mathcal{Z}^{(n)} = \mathcal{Q}^{(n)}, \mathcal{W}^{(n)},
\mathcal{Q}_S^{(n)}, \mathcal{W}_S^{(n)},
\mathcal{A}^{(n)}, \mathcal{V}^{(n)}$, we define
$\widehat{\mathcal{Z}}^{(n)} (t) (B)
\define\frac{1}{\sqrt{n}} \mathcal{Z}^{(n)} (nt) (\sqrt{n} B)$
for every Borel set $B\subset\R$.

\section{Main results}\label{S3}

Before stating our main results, we summarize
the results for the standard EDF system that
were obtained
in \cite{bogdan}---in particular, we recall
Proposition 3.10 and
Theorem 3.1 of \cite{bogdan} which characterize the limiting
distributions of
the workload measure and the queue length measure in the standard system.
Let
%
\begin{equation}\label{5.1}
H(y)\define\int_{y}^{\infty}\bigl(1-G(\eta)\bigr) \,d\eta
=\cases{
\displaystyle \int_y^{y^*}\bigl(1-G(\eta)\bigr) \,d\eta, &\quad
if $y\leq y^*$,\vspace*{2pt}\cr
0, &\quad if $y>y^*$.}
\end{equation}
The function $H$ maps $(-\infty,y^*]$ onto
$[0,\infty)$ and is strictly decreasing and Lipschitz
continuous with Lipschitz constant $1$ on
$(-\infty,y^*]$. Therefore, there exists a
continuous inverse function $H^{-1}$
that maps $[0,\infty)$ onto $(-\infty,y^*]$.
\begin{proposition}[(Proposition 3.10 \cite{bogdan})]\label{P2.1S}
We have $ \widehat{F}_S^{(n)} \Rightarrow F_S^*$ as $n\rightarrow
\infty$,
where the
limiting scaled frontier process $F_S^*$ for the standard EDF system
is explicitly given by
%
\begin{equation}\label{5.2S}
F_S^*(t)\define H^{-1}(W_S^*(t)),\qquad
t\geq0,
\end{equation}
with $W^*_S$ equal to
Brownian motion
with variance
$(\alpha^2 +\beta^2)\lambda$ per unit time
and drift $-\gamma$, reflected at $0$.
\end{proposition}
\begin{theorem}[(Theorem 3.1 \cite{bogdan})]\label{T2.2S}
Let $\mathcal W_S^*$ and $\mathcal Q_S^*$
be the measure-valued processes defined, respectively, by
%
\begin{equation}\label{5.3S}\qquad
\mathcal{W}_S^*(t)(B)
\define\int_{B\cap[F_S^*(t),\infty)}
\bigl(1-G(y)\bigr) \,dy,\qquad
\mathcal{Q}_S^*(t)(B)
\define\lambda\mathcal{W}_S^*(t)(B),
\end{equation}
for all Borel sets $B\subseteq\R$.
Then $\widehat{\mathcal{W}}_S^{(n)} \Rightarrow\mathcal{W_S^*}$ and
$\widehat{\mathcal{Q}}_S^{(n)} \Rightarrow\mathcal{Q_S^*}$, as
$n \ra\infty$.
\end{theorem}
\begin{remark}\label{R3.3}
The proofs in \cite{bogdan} can be
modified to show that the
convergences in (\ref{5.3S}) are in fact joint,
that is,
$(\widehat{\mathcal{W}}_S^{(n)},
\widehat{\mathcal{Q}}_S^{(n)})\Rightarrow
(\mathcal{W}^*_S,\mathcal{Q}_S^*)$.
\end{remark}

There is lateness in the standard EDF system if
and only if the
measure-valued workload process has positive
mass on the negative half line.
Theorem \ref{T2.2S} shows that, in the heavy
traffic limit, this occurs
exactly when the limiting scaled frontier
process $F_S^*$ lies to the
left of $0$ or, equivalently (by Proposition
\ref{P2.1S}), when
$W_S^*$ is greater than
$H(0)=\E[L^{(n)}_j/\sqrt{n}]$, the mean
of the scaled lead-time distribution.
In the reneging system, there is no
lateness, and the amount of work that reneges is precisely the amount
required to prevent lateness. Thus it is natural to expect that
the limiting workload in the reneging system will be constrained to
remain below~$H(0)$.
Let $W^*$ be a Brownian motion with variance
$(\alpha^2 +\beta^2)\lambda$ per unit time
and drift $-\gamma$, reflected at $0$ and $H(0)$.
The first main result of this paper is
that $W^*$ is the limiting workload in the
reneging system.
%
%
\begin{theorem} \label{C2.1}
As $n \ra\infty$, $\widehat{W}^{(n)} \Rightarrow W^*$.
\end{theorem}

The next two results of this paper are the following counterparts of
Proposition~\ref{P2.1S} and
Theorem \ref{T2.2S} for the EDF system
with reneging.
\begin{proposition}\label{P2.1}
We have $\widehat{F}^{(n)} \Rightarrow F^*$ as $n\rightarrow\infty$, where
%
\begin{equation}\label{5.2}
F^*(t)\define H^{-1}(W^*(t)),\qquad
t\geq0.
\end{equation}
\end{proposition}

In other words, the process $F^*$
defined by (\ref{5.2}) is the limiting scaled
frontier process for the EDF system with
reneging.
\begin{theorem}\label{T2.2}
Let $\mathcal W^*$ and $\mathcal Q^*$ be the measure-valued processes
defined by
%
\begin{equation}\label{5.3}\hspace*{28pt}
\mathcal{W}^*(t)(B)
\define\int_{B\cap[F^*(t),\infty)}
\bigl(1-G(y)\bigr) \,dy,\qquad
\mathcal{Q}^*(t)(B)
\define\lambda\mathcal{W}^*(t)(B),
\end{equation}
for all Borel sets $B\subseteq\R$.
Then $(\widehat{\mathcal{W}}^{(n)}, \widehat{\mathcal{Q}}^{(n)})
\Rightarrow
(\mathcal{W^*},\mathcal{Q}^*)$ as $n \ra\infty$.
\end{theorem}

By Theorem \ref{T2.2}, the total masses of
$\sW^{(n)}$ and $\sQ^{(n)}$ must converge jointly
to the total masses of $\sW^*$ and $\sQ^{(n)}$, respectively.
Substituting
$B =\R$ in (\ref{5.3}) and using
(\ref{5.1}) and (\ref{5.2}), we see that
$\sW^* (t) (\R)= H (F^*(t)) = W^*(t)$
and we recover Theorem \ref{C2.1}.
In fact, we have a stronger result.
\begin{corollary}\label{C3.5a}
As $n\rightarrow\infty$, $(\widehat{W}^{(n)},\widehat
{Q}^{(n)})\Rightarrow
(W^*,\lambda W^*)$.
\end{corollary}

Theorem \ref{T2.2} also shows that
the limiting instantaneous lead-time profiles of
customers in the EDF system with reneging
\textit{conditioned on the value of
the} (\textit{limiting}) \textit{workload in the system} are the
same as in the case of the
standard EDF system. However, the limiting real-valued
workload process for the EDF
system with reneging is $W^*$, the
\textit{doubly reflected} Brownian motion and
the \textit{unconditional} limiting lead-time profiles
for these two systems differ accordingly.

We also have a characterization of the
limiting amount of reneged work.
\begin{theorem}\label{T3.1}
As $n \ra\infty$, $\hrenw\Rightarrow\hrew$, where $\hrew$ is
the local time at $H(0)$ of the doubly reflected
Brownian motion $W^*$.
\end{theorem}

Although these results are intuitive in light of the behavior of the
standard EDF system, the proofs are challenging. Moreover, counter to
what one might expect, the result for queue lengths analogous to
Theorem \ref{T3.1} is false. Specifically, although Corollary
\ref{C3.5a} shows that $\widehat{Q}^{(n)}$ converges to the doubly
reflected Brownian motion $Q^* \define\lambda W^*$ on $[0,\lambda
H(0)]$, the scaled sequence $\hrenq, n \in\N$, of reneged customers
\textit{does not} converge to the local time $\lambda R_W^*$ of $Q^*$
at $\lambda H(0)$. This observation, which is elaborated upon in
Section \ref{S.sim}, emphasizes the need for a rigorous justification
of intuitive statements.

The proof of Theorem \ref{C2.1} is in Section \ref{subsub-pre}, the
proofs of Proposition \ref{P2.1} and Theorem \ref{T2.2} are in Section
\ref{S.mainp}, and Section \ref{subs-htren} contains the proof of
Theorem~\ref{T3.1}. We also establish an optimality property for EDF,
Theorem \ref{t.EDFopt}.
\begin{remark}\label{rem-ystar}
The assumption made in (\ref{2.5}) that the support of the lead time
distribution is bounded above by $y^* < \infty$ is mainly technical. It
is expected that the analysis in \cite{lk} for the standard EDF system
under a weaker second moment condition can be applied to the reneging
system as well. On the other hand, the lower bound $y_* > 0$ on the
lead time distribution or some restriction on the behavior of the
density of the lead time distribution at $0$ appears to be necessary.
Indeed, the work of Ward and Glynn \cite{WaGl,WaGl2} on FIFO queues
with reneging suggests that in the absence of such an assumption, the
limiting workload process may no longer be a reflected Brownian motion,
and its properties may exhibit strong sensitivity to the density of the
lead-time distribution near $0$. From a modeling point of view, it is
reasonable to impose a strictly positive lower bound $y_* > 0$ so as to
avoid nonnegligible ``intrinsic lateness,'' in which an arriving
customer has such a small initial lead time that he would be late even
if there were no other customers in the system.

In \cite{lk} the assumption of independence between
the sequence of interarrival times and lead times
is also removed and a more complex version
of Theorem~\ref{T2.2S} is obtained. Starting
from that more complex result, the limit
of the reneging system can be obtained along
the lines of this paper.
\end{remark}

\section{The reference system}\label{S4}

In this section we introduce an auxiliary reference
workload measure-valued process $\mathcal{U}^{(n)}$ and
the corresponding real-valued reference
workload process $U^{(n)}$.
In the special case of constant initial lead times (i.e., $y_*=y^*$),
in which EDF reduces to the well-known FIFO service
discipline, $\mathcal{U}^{(n)}$~and $U^{(n)}$ coincide
with $\mathcal{W}^{(n)}$ and $W^{(n)}$, respectively.
In general, these processes do not coincide (see Example \ref{Ex4.3}),
but, as we will show in Section \ref{subs-htwork}, the difference
between the
diffusion-scaled versions of $U^{(n)}$ and $W^{(n)}$
is negligible under heavy-traffic conditions.
The advantage of working with the reference system, rather than
the reneging system, is that $\mathcal{U}^{(n)}$ can be represented
explicitly as a certain mapping $\Phi$
of the measure-valued workload process $\mathcal{W}_S^{(n)}$ in the
standard system. As shown in
Section \ref{subs-htwork}, continuity properties of the mapping $\Phi$
enable an easy characterization of the limiting distributions of
$\mathcal{U}^{(n)}$ and
$U^{(n)}$ in heavy traffic.

We begin with Section \ref{subs-refdef}, where we define the reference
system and
provide a useful decomposition of the process $U^{(n)}$.
In Section \ref{subs-refdyn} we provide a detailed description of the
evolution of $\mathcal{U}^{(n)}$.

\subsection{Definition and properties of the reference workload}
\label{subs-refdef}

In Section \ref{subsub-mapmvp}, we introduce a deterministic mapping
on the space of measure-valued
functions that is used to define the reference workload.
Then, in Section \ref{subsub-Kdcomp}, we provide a decomposition
of the reference workload process.

\subsubsection{A mapping $\Phi$ of measure-valued processes}
\label{subsub-mapmvp}

We define a sequence of \textit{reference workload measure-valued processes}
for the EDF system with reneging by the formula
%
\begin{equation}\label{refworm}
\mathcal{U}^{(n)} \define\Phi\bigl(\mathcal{W}^{(n)}_S\bigr),
\end{equation}
where the mapping
$\Phi\dvtx D_\sMM[0,\infty)
\mapsto D_\sMM[0,\infty)$
is defined by
%
\begin{eqnarray}\label{Phi}
&&\Phi(\mu)(t) (-\infty,y]\nonumber\\[-8pt]\\[-8pt]
&&\qquad\define
\Bigl[ \mu(t)(-\infty,y]
- \sup_{s \in[0,t]}
\Bigl(\mu(s)(-\infty,0] \wedge\inf_{u \in[s,t]}
\mu(u)(\R)\Bigr)\Bigr]^+
\nonumber
\end{eqnarray}
for every $\mu\in D_\sMM[0,\infty)$, $t\geq0$ and $y\in\R$.
(The claim that $\Phi$ does indeed map $D_\sMM[0,\infty)$ into
$D_\sMM[0,\infty)$
is justified in Lemma \ref{lem-Phi} below.)
We also define the (real-valued) \textit{reference workload process}
$U^{(n)}$ as
the total mass of $\mathcal{U}^{(n)}$, that is,
%
\begin{equation}\label{refwor}
U^{(n)}(t) \define\mathcal{U}^{(n)}(t)(\R)\qquad \forall t \in[0,\infty).
\end{equation}
The frontier $F_S^{(n)}$ defined in Section \ref{2.3}
played a crucial role in the
description and analysis of the evolution of the standard system
in \cite{bogdan}.
In a similar fashion, it will be useful to define the \textit{reference frontier}
%
\begin{equation}\label{E}
E^{(n)}(t) \define
\cases{
\inf\bigl\{y\in\R|\mathcal{U}^{(n)}(t)(-\infty,y]>0\bigr\}, &\quad
if $U^{(n)}(t)>0$,\cr
+\infty, &\quad if $U^{(n)}(t)=0$.}
\end{equation}
By definition, $E^{(n)}(t)$ is the leftmost point of
support of the random measure $\mathcal{U}^{(n)}(t)$
[understood as $\infty$
if
$\mathcal{U}^{(n)}(t)\equiv0$].
The process $E^{(n)}$
has RCLL paths.

From (\ref{refworm})--(\ref{refwor}) we have
%
\begin{eqnarray}
\label{calU}
\mathcal{U}^{(n)}(t)(-\infty,y]&=&
\bigl[ \mathcal{W}_S^{(n)}(t)(-\infty,y]- K^{(n)}(t)\bigr]^+,
\\
\label{U}
U^{(n)}(t) &=& W_S^{(n)}(t) - K^{(n)}(t),
\end{eqnarray}
where
%
\begin{equation}\label{C}
K^{(n)}(t) \define\max_{s \in[0,t]}
\Bigl\{\mathcal{W}^{(n)}_S(s)(-\infty,0] \wedge
\inf_{u \in[s,t]} W^{(n)}_S(u) \Bigr\}.
\end{equation}
In (\ref{C}) we may write maximum rather than supremum
because the process $\sW^{(n)}_S(\cdot)(-\infty,0]$ never
jumps down. Note from (\ref{U}) and (\ref{C}) that $0 \leq K^{(n)}
(t) \leq W_S^{(n)} (t)$ and
so for all $t \geq0$,
%
\begin{equation}
\label{eq-uws}
0 \leq U^{(n)} (t) \leq W_S (t).
\end{equation}
According to (\ref{U}), the reference workload process $U^{(n)}$ is
the standard
workload process $W^{(n)}_S$ with mass $K^{(n)}$ removed.
Equation (\ref{calU}) shows that this mass is removed
from the left-hand side of the support of $\sW_S^{(n)}$.
Moreover, since $\mathcal{U}^{(n)}(t) (-\infty, y] > 0$ for all $y$
to the
right of the frontier $E^{(n)}(t)$,
it is clear from (\ref{refworm}) and (\ref{Phi}) that for
$t \in[0,\infty)$,
$y_2 \geq y_1 > E^{(n)}(t)$,
%
\begin{equation}
\label{e.19}
\mathcal{U}^{(n)}(t) (y_1, y_2] = \mathcal{U}^{(n)}(t)(-\infty,y_2] -
\mathcal{U}^{(n)}(t) (-\infty, y_1] = \mathcal{W}_S^{(n)}(t) (y_1,
y_2],\hspace*{-32pt}
\end{equation}
which shows that $\mathcal{U}^{(n)}$ coincides with
$\mathcal{W}_S^{(n)}$
strictly to the right of $E^{(n)}$.

In the following lemma, we establish some basic properties of
$\Phi$ that show, in particular, that
$\mathcal{U}^{(n)}(t)$, $t\geq0$, and $U^{(n)} (t)$, $t \geq0$,
are stochastic processes with sample paths in
$D_\sMM[0,\infty)$ and $D_{\R_+} [0,\infty)$,
respectively. Although $\Phi$ is not continuous
on $D_{\sMM}[0,\infty)$,
the lemma shows that it satisfies a certain
continuity property that will be sufficient for our purposes.
\begin{lemma}
\label{lem-Phi}
For every $t \in[0,\infty)$, $\Phi(\mu)(t) (-\infty, 0] = 0$.
Moreover, $\Phi$ maps $D_\sMM[0,\infty)$ to
$D_\sMM[0,\infty)$.
Furthermore, if a sequence $\mu_n, n \in\N$, in
$D_{\sMM}[0,\infty)$ converges to $\mu\in D_{\sMM}[0,\infty)$,
where $\mu$ is continuous and for every $t \in[0,\infty)$,
$\mu(t)\{0\} = 0$,
then $\Phi(\mu_n)$ converges to $\Phi(\mu)$ in $D_{\sMM}[0,\infty)$.
\end{lemma}
\begin{pf}
The first statement follows from the simple observation that, due to
the nonnegativity of $\mu$ and (\ref{Phi}),
\[
0 \leq\Phi(\mu) (t) (-\infty, 0] \leq\bigl[ \mu(t) (-\infty, 0] - \mu
(t) (-\infty, 0] \wedge\mu(t) (\R) \bigr]^+ = 0.
\]
Also, since the right-hand side of (\ref{Phi})
is nondecreasing and right-continuous in~$y$, we know that
$\Phi(\mu) (t) \in\sMM$ for every $t \geq0$.
Now, observe that $\Phi(\mu) (t) = \Psi(\mu(t), \Gamma(\mu) (t))$,
where $\Psi\dvtx \mathcal{M} \times\R\mapsto\mathcal{M}$ is the mapping
$\Psi(\nu, x) (-\infty, y] \define
( \nu(-\infty, y] - x )^+$
for all $y \in\R$
and $\Gamma\dvtx D_\sMM[0,\infty) \mapsto\R$ is defined by
\[
\Gamma(\mu) (t) \define\sup_{s \in[0,t]} \Bigl( \mu(s) (-\infty, 0]
\wedge\inf_{u \in[s,t]}
\mu(u) (\R) \Bigr)\qquad \forall t \in[0,\infty).
\]
Using the fact that weak convergence of measures on $\R$ is equivalent to
convergence of the cumulative distribution functions at continuity
points of the limit,
one can verify that $\Psi$ is continuous on
$\mathcal{M} \times\R$.
To show that $\Phi(\mu) \in D_\sMM[0,\infty)$, it suffices to show that
$\Gamma(\mu) \in D[0,\infty)$. For this, we fix $t \in[0,\infty)$
and write
\begin{eqnarray*}
&&
\Gamma(\mu) (t + \ve) - \Gamma(\mu) (t) \\
&&\qquad=
\sup_{s \in[0,t]} \Bigl[ \mu(s) (-\infty, 0] \wedge\inf_{u \in[s,t]}
\mu(u) (\R) \wedge
\inf_{u \in[t,t+\ve] } \mu(u) (\R) \Bigr] \vee Z(\mu,\ve) (t) \\
&&\qquad\quad{} - \sup_{s \in[0,t]}
\Bigl[ \mu(s) (-\infty, 0] \wedge\inf_{u \in[s,t]} \mu(u) (\R) \Bigr],
\end{eqnarray*}
where we define
\[
Z(\mu, \ve) (t) \define\sup_{s \in[t,t+\ve]} \Bigl[ \mu(s) (-\infty,
0] \wedge
\inf_{u \in[s, t+\ve]} \mu(u) (\R) \Bigr].
\]
Since $\mu\in D_{\sMM}[0,\infty)$ implies $\mu(u)$ converges weakly to
$\mu(t)$ as $u \downarrow t$,
we have
$\lim_{u \downarrow t} \mu(u) (\R)
= \mu(t) (\R)$ and
$\mu(t) (-\infty, 0] \geq
\limsup_{s \downarrow t} \mu(s) (-\infty, 0]$
by Portmanteau's theorem.
This, in turn, implies that $\lim_{\ve\ra0} Z(\mu,\ve)(t)=\mu
(t)(-\infty,0]$ for all $t \geq0$.
Combining the above properties, it is easy to deduce that $\Gamma(\mu
) (t+ \ve) - \Gamma(\mu) (t)
\ra0$ as $\ve\downarrow0$, and the right-continuity of
$\Phi(\mu)$ follows. The existence of left limits
for $\Gamma(\mu)$, and hence for $\Phi(u)$, can be established
by an analogous but simpler argument.

Now, suppose $\mu_n$ converges to $\mu$ in $D_{\sMM}[0,\infty)$
and $\mu$ is continuous with $\mu(t)\{0\} = 0$
for every $t \geq0$. Then
$\mu_n(t)$ converges weakly to $\mu(t)$
uniformly for $t$ in compact sets (u.o.c.)
(see \cite{bill}). Since
$0$ is a continuity point for $\mu(t)$, this implies
$\mu_n (t)(-\infty, 0]$ and $\mu_n (t)(\R)$ converge u.o.c. to
$\mu(t) (-\infty, 0]$ and $\mu(t) (\R)$, respectively.
This shows that $\Gamma(\mu_n)(t)$ converges u.o.c. to $\Gamma(\mu)(t)$,
which, when combined with the continuity of $\Psi$, shows that $\Phi
(\mu_n) (t)$ converges
weakly u.o.c. to $\Phi(\mu)(t)$. In particular, this shows
$\Phi(\mu_n)$ converges to $\Phi(\mu)$ in $D_{\sMM}[0,\infty)$.
\end{pf}

As an immediate consequence of the lemma,
the definitions of $\mathcal{U}^{(n)}$ and
$E^{(n)}$, and the fact that $\sU^{(n)}(t)$
is a purely atomic measure, we have, for all $t \geq0$,
%
\begin{equation}
\label{E>0}
\mathcal{U}^{(n)}(t) (-\infty,0] =0 \quad\mbox{and}\quad
E^{(n)} (t) > 0.
\end{equation}

\subsubsection{A decomposition of the reference workload}
\label{subsub-Kdcomp}

We establish a decomposition of $K^{(n)}$
into its increasing and decreasing parts. Define $\sigma
^{(n)}_0\define0$ and $W_S^{(n)}(0-)\define0$. For
$k=0,1,2,\ldots,$
define recursively
%
\begin{eqnarray}
\label{2}\quad
\tau^{(n)}_k
&\define&
\min\Bigl\{t\geq\sigma_k^{(n)}|
W_S^{(n)}\bigl(\sigma^{(n)}_k-\bigr)
\vee\max_{s \in[\sigma_k^{(n)}, t]}
\mathcal{W}_S^{(n)}(s)(-\infty,0]\nonumber\\[-8pt]\\[-8pt]
&&\hspace*{20.9pt}\hspace*{172.9pt}{} \geq W^{(n)}_S(t)
\Bigr\},\nonumber\\
\label{3}
\sigma^{(n)}_{k+1}
&\define&
\min\bigl\{t\geq\tau_k^{(n)}
|W^{(n)}_S(t)>W^{(n)}_S(t-)\bigr\}.
\end{eqnarray}
In addition, for $t \in[0,\infty)$, define
%
\begin{eqnarray}\qquad
\label{def-kpn}
\kpn(t)
&\define&
\sum_{k \in\N} \Bigl[
W_S^{(n)} \bigl(\sigma_k^{(n)} -\bigr)
\vee\max_{s \in[\sigma_k^{(n)},t \wedge
\tau_{k}^{(n)}]} \sW_S^{(n)}(s) (-\infty,0]\nonumber\\[-8pt]\\[-8pt]
&&\hspace*{150.3pt}{}
- W_S^{(n)} \bigl(\sigma_k^{(n)} -\bigr) \Bigr],
\nonumber\\
\label{def-knn}
\knn(t)
&\define&
-\sum_{k \in\N}\bigl[\bigl( W_S^{(n)} \bigl(\tau_{k-1}^{(n)}\bigr)
- \bigl(\sigma_{k}^{(n)} \wedge
t - \tau_{k-1}^{(n)}\bigr)\bigr)^+\nonumber\\[-8pt]\\[-8pt]
&&\hspace*{120.9pt}{}
- W_S^{(n)} \bigl(\tau_{k-1}^{(n)}\bigr) \bigr].
\nonumber
\end{eqnarray}
\begin{theorem}
\label{th-Kdcomp}
We have
%
\begin{equation}
\label{eq-Kdcomp}
K^{(n)} = \kpn- \knn,
\end{equation}
where $\kpn$ and $\knn$ are the positive and negative variations of $K^{(n)}$.
Moreover,
%
\begin{equation}
\label{complementarity}
\int_{[0,\infty)} \ind_{\{U^{(n)}(s) > 0\}} \,d \knn(s) =0.
\end{equation}
\end{theorem}

The theorem is easily deduced from
Propositions \ref{P.3} and \ref{p4} and
Remark \ref{r5} below.
The rest of the section is devoted to
establishing these results.

Observe that the late work $\mathcal{W}_S^{(n)}(s)(-\infty,0]$
is right-continuous in $s$, remaining constant or moving
down at rate one and jumping up. Therefore,
the maximum on the right-hand side of (\ref{2}) is obtained.
Additionally, because of the right-continuity of
$\mathcal{W}_S^{(n)}$ and $W_S^{(n)}$, the minimum in
this equation is also obtained. Finally,
$\mathcal{W}_S^{(n)}(s)(-\infty,0]$ can never exceed
$W_S^{(n)}(s)=\mathcal{W}_S^{(n)}(s)(\R)$,
and $W_S^{(n)}$ never jumps down, so we must in fact have
%
\begin{equation}\label{4}
W_S^{(n)}\bigl(\sigma^{(n)}_k-\bigr)\vee
\max_{s \in[\sigma^{(n)}_k, \tau^{(n)}_k]}
\mathcal{W}_S^{(n)}(s)(-\infty,0]
= W^{(n)}_S\bigl(\tau^{(n)}_k\bigr).
\end{equation}
For $k\geq1$, $\sigma^{(n)}_k$ is the first arrival
time after $\tau^{(n)}_{k-1}$. We thus have
%
\begin{equation}\label{5}
W^{(n)}_S(t)
= \bigl(W^{(n)}_S\bigl(\tau^{(n)}_{k-1}\bigr)-
\bigl(t-\tau^{(n)}_{k-1}\bigr)\bigr)^+,\qquad
\tau^{(n)}_{k-1}\leq t <\sigma^{(n)}_k.
\end{equation}
We further have
%
\begin{equation}\label{6}
0=\sigma^{(n)}_0=\tau^{(n)}_0< \sigma^{(n)}_1 <\tau^{(n)}_1
<\sigma^{(n)}_2<\cdots.
\end{equation}
\begin{proposition}\label{P.3}
For each $k\geq1$, we have
%
\begin{equation}\label{7}
K^{(n)}(t) = W_S^{(n)}\bigl(\sigma^{(n)}_k-\bigr)
\vee\max_{s \in[\sigma^{(n)}_k, t]}
\mathcal{W}_S^{(n)}(s)(-\infty,0]
\end{equation}
for $t\in[\sigma^{(n)}_{k},\tau^{(n)}_k]$.
In particular, $K^{(n)}$ is nondecreasing on the interval
$[\sigma^{(n)}_{k},\tau^{(n)}_k]$.
\end{proposition}
\begin{pf}
We proceed by induction on $k$. For the base case $k=1$,
note that the standard EDF system is empty before the time $\sigma^{(n)}_1$.
Therefore, $W^{(n)}_S(\sigma^{(n)}_1-)=0$, and to prove (\ref{7}),
we must show that
%
\begin{equation}\label{8}
K^{(n)}(t)
=\max_{s \in[0, t]} \mathcal{W}_S^{(n)}(s)(-\infty,0],\qquad
\sigma^{(n)}_{1}\leq t \leq\tau^{(n)}_1.
\end{equation}
For $t\in[\sigma^{(n)}_1,\tau^{(n)}_1]$,
let $s^{(n)}(t)$ be the largest number in
$[\sigma^{(n)}_1,t]$ satisfying
%
\begin{equation}\label{9}
\mathcal{W}_S^{(n)}\bigl(s^{(n)}(t)\bigr)(-\infty,0]
=\max_{s \in[0,t]} \mathcal{W}_S^{(n)}(s)(-\infty,0].
\end{equation}
For $u\in[s^{(n)}(t),t]$, we have
\[
\mathcal{W}_S^{(n)}\bigl(s^{(n)}(t)\bigr)(-\infty,0]
=\max_{s \in[\sigma^{(n)}_1, u]} \mathcal{W}_S^{(n)}(s)(-\infty,0],
\]
which is less than or equal to $W_S^{(n)}(u)$ by the definition
of $\tau^{(n)}_1$ and equation (\ref{4}). Therefore,
\[
\max_{s \in[0,t]} \mathcal{W}_S^{(n)}(s)(-\infty,0]
=\mathcal{W}_S^{(n)}\bigl(s^{(n)}(t)\bigr)(-\infty,0]
\leq\inf_{u \in[s^{(n)}(t),t]}W_S^{(n)}(u).
\]
Equation (\ref{8}) follows from (\ref{C}).

We assume (\ref{7}) holds for some $k$
and prove it for $k+1$. For $t\in[\sigma^{(n)}_{k+1},\tau^{(n)}_{k+1}]$,
%
\begin{eqnarray} \label{10}
K^{(n)}(t) &=& \max_{s \in[0, \sigma^{(n)}_{k+1}) }
\Bigl\{\mathcal{W}^{(n)}_S(s)(-\infty,0] \wedge
\inf_{u \in[s, t]} W^{(n)}_S(u) \Bigr\}\nonumber\\[-8pt]\\[-8pt]
&&{} \vee\max_{s \in[\sigma^{(n)}_{k+1},t]}
\Bigl\{\mathcal{W}^{(n)}_S(s)(-\infty,0] \wedge
\inf_{u \in[s,t]} W^{(n)}_S(u) \Bigr\}.\nonumber
\end{eqnarray}
Equation (\ref{7}) with $k$ replaced by $k+1$ will
follow once we
show that
%
\begin{equation}\label{11}
\max_{s \in[0,\sigma^{(n)}_{k+1}) }
\Bigl\{ \mathcal{W}^{(n)}_S(s)(-\infty,0] \wedge
\inf_{u \in[s,t]} W^{(n)}_S(u) \Bigr\}
= W_S^{(n)}\bigl(\sigma^{(n)}_{k+1}-\bigr)
\end{equation}
and
%
\begin{eqnarray}\label{12}
&&\max_{s \in[\sigma^{(n)}_{k+1}, t]}
\Bigl\{\mathcal{W}^{(n)}_S(s)(-\infty,0] \wedge
\inf_{u \in[s,t]} W^{(n)}_S(u) \Bigr\}\nonumber\\[-8pt]\\[-8pt]
&&\qquad= \max_{s \in[\sigma^{(n)}_{k+1}, t]}
\mathcal{W}_S^{(n)}(s)(-\infty,0].\nonumber
\end{eqnarray}

For (\ref{11}), we observe that because
$\mathcal{W}^{(n)}_S(s)(-\infty,0]$ and $\inf_{s\leq u\leq t} W^{(n)}_S(u)$,
regarded as functions of $s$, cannot increase except by a jump,
the maximum on the left-hand side of (\ref{11}) is attained.
Let $s^{(n)}_k$ be the largest number in $[0,\sigma^{(n)}_{k+1})$
attaining this maximum. We have
\begin{eqnarray*}
&&\max_{s \in[0, \sigma^{(n)}_{k+1}) }
\Bigl\{ \mathcal{W}^{(n)}_S(s)(-\infty,0] \wedge
\inf_{u \in[s, t]} W^{(n)}_S(u) \Bigr\}
\\
&&\qquad= \mathcal{W}^{(n)}_S\bigl(s^{(n)}_k\bigr)(-\infty,0] \wedge
\inf_{u \in[s^{(n)}_k, t]} W^{(n)}_S(u)
\leq W^{(n)}_S(u) \qquad\forall u\in\bigl[s^{(n)}_k,\sigma^{(n)}_{k+1}\bigr),
\end{eqnarray*}
and so
%
\begin{equation}\label{13}
\max_{s \in[0, \sigma^{(n)}_{k+1}) }
\Bigl\{ \mathcal{W}^{(n)}_S(s)(-\infty,0] \wedge\inf_{u \in[s, t]}
W^{(n)}_S(u) \Bigr\} \leq W_S^{(n)}\bigl(\sigma^{(n)}_{k+1}-\bigr).
\end{equation}
On the other hand, by the inequalities $\tau^{(n)}_k<\sigma
^{(n)}_{k+1}\leq t \leq\tau^{(n)}_{k+1}$, definition (\ref{C}), the
induction hypothesis, and equation (\ref{4}),
we have
\begin{eqnarray*}
&&\max_{s \in[0,\sigma^{(n)}_{k+1}) }
\Bigl\{ \mathcal{W}^{(n)}_S(s)(-\infty,0]
\wedge\inf_{u \in[s, t]} W^{(n)}_S(u) \Bigr\}\\
&&\qquad\geq
\max_{s \in[0, \tau^{(n)}_{k}] }
\Bigl\{ \mathcal{W}^{(n)}_S(s)(-\infty,0]
\wedge\inf_{u \in[s, \tau^{(n)}_k]} W^{(n)}_S(u)
\wedge\inf_{u \in[\tau^{(n)}_k, t]}
W^{(n)}_S(u)\Bigr\}\\
&&\qquad=
K^{(n)}\bigl(\tau^{(n)}_k\bigr) \wedge
\inf_{u \in[\tau^{(n)}_k, t]} W^{(n)}_S(u)\\
&&\qquad=
\Bigl(W_S^{(n)}\bigl(\sigma^{(n)}_k-\bigr) \vee
\max_{s \in[\sigma^{(n)}_k, \tau^{(n)}_k]}
\mathcal{W}_S^{(n)}(s)(-\infty,0]\Bigr)
\wedge\inf_{u \in[\tau^{(n)}_k, t]}
W^{(n)}_S(u) \\
&&\qquad= W_S^{(n)}\bigl(\tau^{(n)}_k\bigr)
\wedge\inf_{u \in[\tau^{(n)}_k, t]}
W^{(n)}_S(u)=
\inf_{u \in[\tau^{(n)}_k, t]} W^{(n)}_S(u).
\end{eqnarray*}
Equation (\ref{5}) implies
$W^{(n)}_S(u)\geq W^{(n)}_S(\sigma_{k+1}^{(n)}-)$
for $\tau^{(n)}_k\leq u<\sigma^{(n)}_{k+1}$. For
$\sigma_{k+1}^{(n)}\leq u\leq t<\tau^{(n)}_{k+1}$,
(\ref{2}) implies that
\[
W_S^{(n)}\bigl(\sigma^{(n)}_{k+1}-\bigr)\vee
\max_{s \in[\sigma_{k+1}^{(n)}, u]}
\mathcal{W}_S^{(n)}(s)(-\infty,0]
\leq W^{(n)}_S(u),
\]
and so again we have
$W^{(n)}_S(u)\geq W^{(n)}_S(\sigma_{k+1}^{(n)}-)$.
Finally, if $u=t=\tau^{(n)}_{k+1}$,
then (\ref{4}) implies that
$W^{(n)}_S(u)\geq W^{(n)}_S(\sigma_{k+1}^{(n)}-)$.
It follows that
\[
\inf_{u \in[\tau^{(n)}_k, t]}
W^{(n)}_S(u)\geq W^{(n)}_S\bigl(\sigma_{k+1}^{(n)}-\bigr).
\]
This gives the reverse of the
inequality (\ref{13}), and thus (\ref{11})
is proved.

For (\ref{12}), we let $t^{(n)}_k$ attain the maximum
in $\max_{s \in[\sigma^{(n)}_{k+1}, t]}\sW^{(n)}_S(s)(-\infty,0]$.
For $u\in[t^{(n)}_k,t]$, we have from (\ref{2}) and (\ref{4}) that
\[
\sW^{(n)}_S\bigl(t^{(n)}_k\bigr)(-\infty,0]
=\max_{s \in[\sigma^{(n)}_{k+1}, u]}\sW^{(n)}_S(s)(-\infty,0]
\leq W^{(n)}_S(u),
\]
and hence
$\sW^{(n)}_S(t_k^{(n)})(-\infty,0]
\leq\inf_{u \in[t^{(n)}_k, t]}W^{(n)}_S(u)$.
It follows that
\begin{eqnarray*}
&&\max_{s \in[\sigma^{(n)}_{k+1}, t]}
\Bigl\{\sW^{(n)}_S(s)(-\infty,0]\wedge\inf_{u \in[s,
t]}W_S^{(n)}(u)\Bigr\}\\
&&\qquad\leq
\sW^{(n)}_S\bigl(t^{(n)}_k\bigr)(-\infty,0]
=
\sW^{(n)}_S\bigl(t^{(n)}_k\bigr)(-\infty,0]\wedge\inf_{u \in[t^{(n)}_k,
t]}W^{(n)}_S(u) \\
&&\qquad\leq
\max_{s \in[\sigma^{(n)}_{k+1}, t]}
\Bigl\{\sW^{(n)}_S(s)(-\infty,0]\wedge\inf_{u \in[s,t]}W^{(n)}_S(u)\Bigr\},
\end{eqnarray*}
which establishes (\ref{12}).
\end{pf}
\begin{proposition}\label{p4}
For each $k\geq1$, we have
%
\begin{equation}\label{14}
K^{(n)}(t)=\bigl(W^{(n)}_S\bigl(\tau^{(n)}_{k-1}\bigr)-\bigl(t-\tau^{(n)}_{k-1}\bigr)\bigr)^+,\qquad
\tau^{(n)}_{k-1}\leq t<\sigma^{(n)}_k.
\end{equation}
In particular, $K^{(n)}$ is nonincreasing on
$[\tau^{(n)}_{k-1},\sigma^{(n)}_k)$.
\end{proposition}
\begin{pf}
For all $t \geq0$, we have $K^{(n)}(t)\leq W^{(n)}_S(t)$, and for
$\tau^{(n)}_{k-1}\leq t< \sigma^{(n)}_k$, we further have
from (\ref{5}) that
%
\begin{equation}\label{15}
K^{(n)}(t)\leq W^{(n)}_S(t)=\bigl(W^{(n)}_S\bigl(\tau^{(n)}_{k-1}\bigr)-\bigl(t-\tau
^{(n)}_{k-1}\bigr)\bigr)^+.
\end{equation}
On the other hand, Proposition \ref{P.3}
and (\ref{4}) with $k$ replaced by
$k-1$ imply
\begin{eqnarray*}
&&\max_{s \in[0, \tau^{(n)}_{k-1}]}\Bigl\{\sW^{(n)}_S(s)(-\infty,0]
\wedge\inf_{u \in[s, \tau^{(n)}_{k-1}]}W^{(n)}_S(u)\Bigr\} \\
&&\qquad=
K^{(n)}\bigl(\tau^{(n)}_{k-1}\bigr)
=W^{(n)}_S\bigl(\sigma^{(n)}_{k-1}-\bigr)\vee\max_{s \in[\sigma^{(n)}_{k-1},
\tau^{(n)}_{k-1}]}
\sW^{(n)}_S(s)(-\infty,0]\\
&&\qquad=
W^{(n)}_S\bigl(\tau^{(n)}_{k-1}\bigr).
\end{eqnarray*}
For $t \in[\tau^{(n)}_{k-1}, \sigma^{(n)}_k)$,
it follows from (\ref{5}) and the above equality that
%
\begin{eqnarray}\label{16}
K^{(n)}(t)
&=&
\max_{s \in[0, t]}\Bigl\{\sW^{(n)}_S(s)(-\infty,0]
\wedge\inf_{u \in[s, t]}W^{(n)}_S(u)\Bigr\}
\nonumber\\
&\geq&
\max_{s \in[0, \tau^{(n)}_{k-1}]}
\Bigl\{\sW^{(n)}_S(s)(-\infty,0]\nonumber\\
&&\hspace*{41.3pt}{}
\wedge\inf_{u \in[s, \tau^{(n)}_{k-1}]}
W^{(n)}_S(u)
\wedge\inf_{u \in[\tau^{(n)}_{k-1}, t]}
W^{(n)}_S(u)\Bigr\}
\nonumber\\
&=&
\max_{s \in[0, \tau^{(n)}_{k-1}]}
\Bigl\{\sW^{(n)}_S(s)(-\infty,0]\wedge
\inf_{u \in[s, \tau^{(n)}_{k-1}]}
W^{(n)}_S(u)\Bigr\}
\\
&&{}
\wedge\bigl(W^{(n)}_S\bigl(\tau^{(n)}_{k-1}\bigr)
-\bigl(t-\tau^{(n)}_{k-1}\bigr)\bigr)^+\nonumber\\
&=&
W^{(n)}_S\bigl(\tau^{(n)}_{k-1}\bigr)\wedge
\bigl(W^{(n)}_S\bigl(\tau^{(n)}_{k-1}\bigr)
-\bigl(t-\tau^{(n)}_{k-1}\bigr)\bigr)^+
\nonumber\\
&=&
\bigl(W^{(n)}_S\bigl(\tau^{(n)}_{k-1}\bigr)
-\bigl(t-\tau^{(n)}_{k-1}\bigr)\bigr)^+.
\nonumber
\end{eqnarray}
Equation (\ref{14})
follows from (\ref{15}) and (\ref{16}).
\end{pf}
\begin{remark}\label{r5}
In light of (\ref{U}) and Proposition \ref{P.3},
we have the characterization of
$\tau^{(n)}_{k}$ as
%
\begin{equation}\label{17}\qquad\quad
\tau^{(n)}_{k}=\min\bigl\{t\geq\sigma^{(n)}_k|K^{(n)}(t)
\geq W^{(n)}_S(t)\bigr\}
=\min\bigl\{t\geq\sigma^{(n)}_k|U^{(n)}(t)=0\bigr\}.
\end{equation}
Because $\sigma_{k+1}^{(n)}$ is the time of first arrival
after $\tau_k^{(n)}$, we in fact have
%
\begin{equation}\label{4.22}
U^{(n)}(t)=0,\qquad \tau_k^{(n)}\leq t<\sigma_{k+1}^{(n)}.
\end{equation}
Evaluating (\ref{7}) at
$\sigma^{(n)}_k$ and using
$W^{(n)}_S(\sigma^{(n)}_k-)
\geq\sW^{(n)}_S(\sigma^{(n)}_k)(-\infty,0]$, we
obtain
%
\begin{equation}\label{18}
K^{(n)}\bigl(\sigma^{(n)}_k\bigr)=W^{(n)}_S\bigl(\sigma^{(n)}_k-\bigr).
\end{equation}
But (\ref{5}) and Proposition \ref{p4} show that
%
\begin{equation}\label{C=W}
K^{(n)}\bigl(\sigma_k^{(n)}-\bigr)=W_S^{(n)}\bigl(\sigma_k^{(n)}-\bigr),
\end{equation}
and so
%
\begin{equation}\label{Cconsts}
\bigtriangleup K^{(n)}\bigl(\sigma_k^{(n)}\bigr)=0.
\end{equation}
By contrast $\bigtriangleup K^{(n)}(\tau_k^{(n)})$ can be positive.
Evaluating (\ref{7}) at $\tau^{(n)}_k$ and
using (\ref{4}), we obtain
%
\begin{equation}\label{19}
K^{(n)}\bigl(\tau^{(n)}_k\bigr)=W^{(n)}_S\bigl(\tau^{(n)}_k\bigr).
\end{equation}
In conclusion,
%
\begin{eqnarray}
\label{20}\hspace*{32pt}
K^{(n)}(t)
&=&
K^{(n)}\bigl(\sigma^{(n)}_{k}\bigr)\vee
\max_{s \in[\sigma^{(n)}_{k}, t]}\sW^{(n)}_S(s)(-\infty,0],\qquad
\sigma^{(n)}_{k}\leq t\leq\tau^{(n)}_k,\\
\label{21}
K^{(n)}(t)
&=&
\bigl(K^{(n)}\bigl(\tau^{(n)}_{k-1}\bigr)-\bigl(t-\tau^{(n)}_{k-1}\bigr)\bigr)^+,\qquad
\tau^{(n)}_{k-1}\leq t<\sigma^{(n)}_k.
\end{eqnarray}
\end{remark}
%

\subsection{Dynamics of the reference workload process}\label{subs-refdyn}

The evolutions of $\mathcal{U}^{(n)}$ and
$\mathcal{W}^{(n)}$ are similar; the difference
between them is asymptotically
negligible.
Before proving
the properties of $\mathcal{U}^{(n)}$,
we provide a summary of these
properties. The reader may
work out the evolution of $\mathcal{W}_S^{(n)}$,
$\mathcal{U}^{(n)}$ and $\mathcal{W}^{(n)}$
in Example \ref{Ex4.3} to follow along.
This example appears in detail in
\cite{Festschrift}.
\begin{example}\label{Ex4.3}
Consider a system realization in which
\begin{eqnarray*}
u_1^{(n)} &=& 1,\qquad v_1^{(n)}=4,\qquad L_1^{(n)}=3,\qquad S_1^{(n)}=1,\\
u_2^{(n)} &=& 1,\qquad v_2^{(n)}=4,\qquad L_2^{(n)}=5,\qquad S_2^{(n)}=2,\\
u_3^{(n)} &=& 3,\qquad v_3^{(n)}=2,\qquad L_3^{(n)}=1,\qquad S_3^{(n)}=5,\\
u_4^{(n)} &=& 2,\qquad v_4^{(n)}=1,\qquad L_4^{(n)}=4,\qquad S_4^{(n)}=7,\\
u_5^{(n)} &=& 2,\qquad v_5^{(n)}=1,\qquad L_5^{(n)}=1,\qquad S_5^{(n)}=9.
\end{eqnarray*}
Recall that $\delta_a$ is a unit point mass at $a$.
It is straightforward to compute
\begin{eqnarray*}
\sW^{(n)}_S(t)
&=&
\cases{
0, &\quad $0\leq t<1$,\vspace*{2pt}\cr
(5-t)\delta_{4-t}, &\quad $1\leq t<2$,\vspace*{2pt}\cr
(5-t)\delta_{4-t}+4\delta_{7-t}, &\quad $2\leq t<5$,\vspace*{2pt}\cr
(7-t)\delta_{6-t}+4\delta_{7-t}, &\quad $5\leq t<7$,\vspace*{2pt}\cr
(11-t)\delta_{7-t}+\delta_{11-t}, &\quad $7\leq t <9$,\vspace*{2pt}\cr
2\delta_{-2}+\delta_1+\delta_2, &\quad $t=9$,
}
\\
\sU^{(n)}(t)
&=&
\cases{
0, &\quad $0\leq t<1$,\vspace*{2pt}\cr
(5-t)\delta_{4-t}, &\quad $1\leq t <2$,\vspace*{2pt}\cr
(5-t)\delta_{4-t}+4\delta_{7-t}, &\quad $2\leq t<4$,\vspace*{2pt}\cr
(8-t)\delta_{7-t}, &\quad $4\leq t<5$,\vspace*{2pt}\cr
(6-t)\delta_{6-t}+4\delta_{7-t}, &\quad $5\leq t<6$,\vspace*{2pt}\cr
(10-t)\delta_{7-t}, &\quad $6\leq t<7$,\vspace*{2pt}\cr
(8-t)\delta_{11-t}, &\quad $7\leq t<8$,\vspace*{2pt}\cr
0, &\quad $8\leq t<9$,\vspace*{2pt}\cr
\delta_2, &\quad $t=9$,}
\\
\sW^{(n)}(t)
&=&
\cases{
0, &\quad $0\leq t<1$,\vspace*{2pt}\cr
(5-t)\delta_{4-t}, &\quad $1\leq t<2$,\vspace*{2pt}\cr
(5-t)\delta_{4-t}+4\delta_{7-t}, &\quad $2\leq t<4$,\vspace*{2pt}\cr
(8-t)\delta_{7-t}, &\quad $4\leq t<5$,\vspace*{2pt}\cr
(7-t)\delta_{6-t}+3\delta_{7-t}, &\quad $5\leq t<6$,\vspace*{2pt}\cr
(9-t)\delta_{7-t}, &\quad $6\leq t<7$,\vspace*{2pt}\cr
(8-t)\delta_{11-t}, &\quad $7\leq t<8$,\vspace*{2pt}\cr
0, &\quad $8\leq t<9$,\vspace*{2pt}\cr
\delta_1, &\quad $t=9$.}
\end{eqnarray*}
\end{example}

Recall\vspace*{1pt} that $K^{(n)}$ is the amount of mass
removed from the standard
workload $W_S^{(n)}$ to obtain the reference workload $U^{(n)}$.
To understand the process $K^{(n)}$, we consider
the dynamics of $\sU^{(n)}$.
In the absence of new arrivals,
all atoms of $\sU^{(n)}$ move left with unit speed.
Moreover, the mass
of the leftmost atom of $\mathcal{U}^{(n)}$ decreases with
unit speed until
it vanishes, corresponding to the work being
done on the most urgent job
in queue until it is served to completion [Proposition \ref{P.Un}(i)].
However, if the leftmost atom of $\mathcal{U}^{(n)}$
hits zero, this atom is immediately removed from
$\mathcal{U}^{(n)}$ [Proposition \ref{P.Un}(ii), (v)].
This may be
interpreted as reneging of a customer or deletion of a late
customer from the system. When there is a new arrival at time
$t$ with lead time not smaller than the leftmost point of support
of $\mathcal{U}^{(n)}(t-)$, and this point of support
is strictly positive,
then a mass of the size $v^{(n)}_{A^{(n)}(t)}$ located
at $L^{(n)}_{A^{(n)}(t)}$ is added to $\mathcal{U}^{(n)}(t-)$
[Proposition \ref{P.Un}(iii)].
Similarly, if there is a
new arrival and the leftmost point of the support of $\mathcal{U}^{(n)}$
hits zero at the same time, then both of the above actions take
place [(\ref{20.5}) of Proposition \ref{P.Un}(v)].
This is the case of a simultaneous new arrival and ejection
of a late customer. The EDF system with reneging $\mathcal{W}^{(n)}$
shows the same behavior in all these cases.
However, if a customer arrives to start a new busy period
for $\sU^{(n)}$
or, if at time $t$, there is a new arrival with lead time
more urgent than the leftmost point of the support of
$\mathcal{U}^{(n)}(t-)$ (i.e., we have a ``preemption''),
then the mass $v^{(n)}_{A^{(n)}(t)}$ associated
with the new arrival is distributed in
$[L^{(n)}_{A^{(n)}(t)},\infty)$, or more precisely,
on some atoms of $\mathcal{W}_S^{(n)}(t)$ located on this
half-line, but it is not necessarily located at the
single atom $L^{(n)}_{A^{(n)}(t)}$.
This possibility is described in Lemma \ref{L.Un} and
Proposition \ref{P.Un}(iv).
In this respect, the evolution of $\sU^{(n)}$
differs from that of $\sW^{(n)}$,
for which all the new mass is always
placed at the lead time of the arriving customer.
Example \ref{Ex4.3} illustrates this.

We now begin the rigorous study of $\mathcal{U}^{(n)}$.
As shown in Section \ref{subs-refdef}, the time interval
$[0,\infty)$
can be decomposed into a union of the
disjoint intervals $(\tau^{(n)}_k, \sigma_{k+1}^{(n)}]$
and $(\sigma_k^{(n)}, \tau_k^{(n)}]$, $k \geq0$,
such that $K^{(n)} = W_S^{(n)}-U^{(n)}$ is
nonincreasing on $(\tau^{(n)}_k, \sigma_{k+1}^{(n)}]$ and
nondecreasing on $(\sigma_k^{(n)}, \tau_k^{(n)}]$.
In Lemma \ref{L.Un} below, we analyze the behavior
of $\mathcal{U}^{(n)}$ on
the time intervals $[\tau^{(n)}_{k-1},\sigma^{(n)}_{k}]$,
$k\geq1$, while Proposition \ref{P.Un} describes the dynamics of
$\mathcal{U}^{(n)}$ on the intervals $(\sigma^{(n)}_{k},\tau^{(n)}_{k})$,
$k\geq1$.
The section ends with Corollary \ref{C.eqn4.9}, which describes
the time evolution of the
reference workload process $U^{(n)}$.

We make use of the following elementary
facts about
the standard workload. Since the interarrival
times are
strictly positive, $\bigtriangleup A^{(n)}(t) \in\{0,1\}$, and
%
\begin{equation}\label{e.11.5}
\mathcal{W}_S^{(n)}(t) = \mathcal{W}_S^{(n)}(t-)+ \bigtriangleup A^{(n)}(t)
v^{(n)}_{A^{(n)}(t)}
\delta_{L^{(n)}_{A^{(n)}(t)}},\qquad
t\geq0,
\end{equation}
which implies
%
\begin{equation}\label{e.11.5a}
\bigtriangleup W_S^{(n)}(t) = \bigtriangleup A^{(n)}(t)
v^{(n)}_{A^{(n)}(t)},\qquad
t\geq0.
\end{equation}
For any functions $f$ and $g$ on
$[0,\infty)$ (taking finite or infinite values)
such that whenever $s<t$ and
$t-s$ is small enough,
$f(s)=f(t-)+t-s$ and $g(s)=g(t-)+t-s$, we have
%
\begin{equation}\label{leftlimits}
\lim_{s\uparrow t}\mathcal{W}_S^{(n)}(s)[f(s),g(s)]
= \mathcal{W}_S^{(n)}(t-)[f(t-),g(t-)].
\end{equation}
This is true because the lead times of the customers present
in the standard system decrease with unit rate.
Equation (\ref{leftlimits}) remains valid if the
closed intervals $[f(\cdot),g(\cdot)]$ are replaced by either
$[f(\cdot),g(\cdot))$, $(f(\cdot),g(\cdot)]$ or $(f(\cdot),g(\cdot))$.
These facts will be used repeatedly in the following arguments,
sometimes without mention.
\begin{lemma}\label{L.Un}
Let $k\geq1$. We have
%
\begin{eqnarray}
\label{U=0}
U^{(n)}(t)&=&0,\qquad \tau^{(n)}_{k-1}\leq t<\sigma^{(n)}_{k},
\\
\label{massj}
\bigtriangleup U^{(n)}\bigl(\sigma^{(n)}_{k}\bigr)
&=&
v^{(n)}_{A^{(n)}(\sigma^{(n)}_{k})},\\
\label{shiftr}
\mathcal{U}^{(n)}\bigl(\sigma^{(n)}_{k}\bigr)
\bigl(-\infty, L^{(n)}_{A^{(n)}(\sigma^{(n)}_{k})}\bigr)&=&0.
\end{eqnarray}
\end{lemma}
\begin{pf}
Equation
(\ref{U=0}) follows immediately from (\ref{U}), (\ref{5}) and
Proposition \ref{p4}.
By (\ref{U}), (\ref{Cconsts}), (\ref{e.11.5a})
and the fact that $\bigtriangleup A^{(n)}(\sigma_k^{(n)})=1$,
we have
\[
\bigtriangleup U^{(n)}\bigl(\sigma^{(n)}_k\bigr)
=\bigtriangleup W_S^{(n)}\bigl(\sigma^{(n)}_k\bigr)
- \bigtriangleup K^{(n)}\bigl(\sigma^{(n)}_k\bigr)=
v^{(n)}_{A^{(n)}(\sigma^{(n)}_{k})},
\]
and (\ref{massj}) follows.
For $y<L^{(n)}_{A^{(n)}(\sigma^{(n)}_{k})}$,
(\ref{calU}), (\ref{e.11.5}),
(\ref{Cconsts}) and (\ref{C=W}) imply
\begin{eqnarray*}
\mathcal{U}^{(n)}\bigl(\sigma^{(n)}_{k}\bigr)(-\infty,y]&=&
\bigl[ \mathcal{W}_S^{(n)}\bigl(\sigma^{(n)}_{k}\bigr)(-\infty,y]
- K^{(n)}\bigl(\sigma^{(n)}_k\bigr) \bigr]^+\\
&=&
\bigl[ \mathcal{W}_S^{(n)}\bigl(\sigma^{(n)}_{k}-\bigr)(-\infty,y]
- K^{(n)}\bigl(\sigma^{(n)}_k-\bigr) \bigr]^+\\
&\leq&
\bigl[ W_S^{(n)}\bigl(\sigma^{(n)}_{k}-\bigr)- K^{(n)}\bigl(\sigma^{(n)}_k-\bigr)
\bigr]^+
=
0,
\end{eqnarray*}
and so (\ref{shiftr}) also follows.
\end{pf}

Lemma \ref{L.Un} shows that $\sigma^{(n)}_k$
begins
a busy period for the reference system.
Equation (\ref{17}) implies
that $U^{(n)}(t)>0$ for
$\sigma^{(n)}_k<t<\tau_k^{(n)}$, and thus the intervals
$[\sigma^{(n)}_k,\tau^{(n)}_k)$, $k \geq1$, are precisely the
busy periods
for the reference system.
We analyze the behavior of $\mathcal{U}^{(n)}$
during these
busy periods.
We start with the observation that,
by (\ref{calU}) and Proposition \ref{P.3}, for
$t \in(\sigma^{(n)}_k, \tau_k^{(n)})$,
%
\begin{eqnarray}
\label{e.2}\qquad\quad
\mathcal{U}^{(n)}(t)(-\infty,y]&=&
\Bigl[ \mathcal{W}_S^{(n)}(t)(-\infty,y]\nonumber\\[-8pt]\\[-8pt]
&&\hspace*{2.8pt}{} - \Bigl(
W_S^{(n)}\bigl(\sigma^{(n)}_k-\bigr) \vee
\max_{s \in[\sigma^{(n)}_k, t]}
\mathcal{W}_S^{(n)}(s)(-\infty,0]\Bigr) \Bigr]^+.\nonumber
\end{eqnarray}
In what follows, given $\nu\in\sMM$ and any interval
$\mathcal{I} \subset\R$, we will use
$\nu_{|\mathcal{I}}$ to denote the measure in $\sMM$
that is zero on $\mathcal{I}^c$ and coincides with
$\nu$ on $\mathcal{I}\dvtx
\nu_{|\mathcal{I}} (B) = \nu( B \cap\mathcal{I} )$ for all
$B\in\sB(\R)$.
\begin{proposition}\label{P.Un}
For $k\geq1$ and $\sigma^{(n)}_{k}< t < \tau^{(n)}_k$,
the following five properties hold:
\begin{longlist}
\item
If $\bigtriangleup A^{(n)}(t)=0$
and $E^{(n)}(t-)>0$, then
%
\begin{eqnarray}
\label{Cconst}
\bigtriangleup K^{(n)}(t)&=&0,\\
\label{Uconst}
\bigtriangleup U^{(n)}(t)&=&0.
\end{eqnarray}
In this case, if $\mathcal{U}^{(n)}(t-)\{E^{(n)}(t-)\} >0$,
then both $\mathcal{U}^{(n)}(\cdot)\{E^{(n)}(\cdot)\}$
and $U^{(n)}(t)$ decrease with unit rate in a neighborhood
of $t$. On the other hand,
if
\[
\mathcal{U}^{(n)}(t-)\bigl\{E^{(n)}(t-)\bigr\}=0,
\]
then
$\mathcal{U}^{(n)}(t)=\mathcal{W}_S^{(n)}(t)|_{[E^{(n)}(t),\infty)}$.
\item
If $\bigtriangleup A^{(n)}(t)=0$
and $E^{(n)}(t-)=0$,
then
%
\begin{equation}\label{Cj1}
\mathcal{U}^{(n)}(t-)\{0\} =\bigtriangleup K^{(n)}(t)
= -\bigtriangleup U^{(n)}(t)
\end{equation}
and $\mathcal{U}^{(n)}(t)=\mathcal{W}_S^{(n)}(t)|_{(0,\infty)}$.
\item
If $\bigtriangleup A^{(n)}(t) = 1$ and
$L^{(n)}_{A^{(n)}(t)}\geq E^{(n)}(t-)>0$,
then (\ref{Cconst}) holds,\break
$\bigtriangleup E^{(n)}(t)\geq0$ and
%
\begin{equation}\label{e.14}
\mathcal{U}^{(n)}(t)
=\mathcal{U}^{(n)}(t-)+v^{(n)}_{A^{(n)}(t)}
\delta_{L^{(n)}_{A^{(n)}(t)}}.
\end{equation}
\item
If $\bigtriangleup A^{(n)}(t) = 1$,
$L^{(n)}_{A^{(n)}(t)}< E^{(n)}(t-)$, then (\ref{Cconst}) holds and
%
\begin{eqnarray}
\label{e.15}
L^{(n)}_{A^{(n)}(t)} &\leq& E^{(n)}(t)
\leq E^{(n)}(t-),\\
\label{e.15.1}
\bigtriangleup U^{(n)}(t)&=&v^{(n)}_{A^{(n)}(t)},\\
\label{e.18.5}
\mathcal{U}^{(n)}(t)_{|(E^{(n)}(t-),\infty)}
&=&
\mathcal{U}^{(n)}(t-)_{|(E^{(n)}(t-),\infty)},\\
\label{e.18}
\mathcal{U}^{(n)}(t)\bigl\{E^{(n)}(t-)\bigr\}
&\geq& \mathcal{U}^{(n)}(t-)\bigl\{E^{(n)}(t-)\bigr\}.
\end{eqnarray}
\item
If $\bigtriangleup A^{(n)}(t) = 1$ and
$L^{(n)}_{A^{(n)}(t)}> E^{(n)}(t-)=0$, then
%
\begin{equation}\label{20.5}
\mathcal{U}^{(n)}(t)= \mathcal{U}^{(n)}(t-)+v^{(n)}_{A^{(n)}(t)}
\delta_{L^{(n)}_{A^{(n)}(t)}}
-\mathcal{U}^{(n)}(t-)\{0\} \delta_0.
\end{equation}
\end{longlist}
\end{proposition}
\begin{pf}
Fix $k \geq1$ and
$t\in(\sigma^{(n)}_k, \tau_k^{(n)})$.
We start with the general observation that,
by (\ref{E}) and (\ref{e.2}),
%
\begin{eqnarray}\label{e.E}
&&E^{(n)}(t)=\min\Bigl\{ y\big|\mathcal{W}_S^{(n)}(t)(-\infty,y]\nonumber\\[-8pt]\\[-8pt]
&&\hspace*{76.9pt}>
W_S^{(n)}\bigl(\sigma^{(n)}_k-\bigr) \vee
\max_{s \in[\sigma^{(n)}_k, t]} \mathcal{W}_S^{(n)}(s)(-\infty,0]
\Bigr\},\nonumber
\end{eqnarray}
and because $\mathcal{W}_S^{(n)}(t)$ is purely atomic,
the minimum on the
right-hand side of (\ref{e.E}) is obtained at
the atom of
$\mathcal{W}_S^{(n)}(t)$ located at $y_0=E^{(n)}(t)$.
In particular,
%
\begin{equation}
\label{wsatom}
\mathcal{W}_S^{(n)}(t) \bigl\{E^{(n)}(t)\bigr\} > 0.
\end{equation}
We now consider each of the five different cases of the proposition.

(i) Let $a=E^{(n)}(t-)$.
By (\ref{E}) and (\ref{calU}), for all $s<t$
sufficiently near $t$,
%
\begin{equation}\label{e.3}
\mathcal{W}_S^{(n)}(s)(-\infty,a/2]\leq K^{(n)}(s).
\end{equation}
Also, for $s \in[t-a/2, t)$ sufficiently near
$t$ so that $A^{(n)} (s) = A^{(n)} (t)$ holds
[such $s$ exist due to the assumption
that $\bigtriangleup A^{(n)} (t) = 0$], we have
%
\begin{equation}\label{e.4}
\mathcal{W}_S^{(n)}(t)(-\infty,0]
\leq\mathcal{W}_S^{(n)}(s)(-\infty,a/2].
\end{equation}
The last two relations show that
$\mathcal{W}_S^{(n)}(t)(-\infty,0]$ $\leq K^{(n)}(t-)$,
and so, by Proposition \ref{P.3}, (\ref{Cconst}) holds.
Equation (\ref{Uconst}) follows from (\ref{U}), (\ref{Cconst}),
(\ref{e.11.5a}) and the
assumption $\bigtriangleup A^{(n)} (t) = 0$.
Because $W_S^{(n)}(t)>0$ [see (\ref{wsatom})], $W_S^{(n)}$ decreases at
unit rate in a neighborhood of $t$
[see (\ref{2.15})--(\ref{2.18})].
In addition, (\ref{U}), (\ref{Cconst}) and the fact that,
again by Proposition
\ref{P.3}, $K^{(n)}$ cannot increase on
$[\sigma^{(n)}_{k}, \tau^{(n)}_k]$ except by a
jump and hence is constant in a neighborhood
of $t$, together imply that $U^{(n)}$ also decreases at
unit rate in a neighborhood of $t$.
Furthermore, the nature of the EDF discipline
and (\ref{wsatom}) show that
at $t$,
the standard system is serving a customer with lead time no
greater than $E^{(n)} (t)$.
Combining the above properties with
the fact that $\mathcal{U}^{(n)}(t)_{|(E^{(n)}(t),\infty)} =
\mathcal{W}_S^{(n)}(t)_{|(E^{(n)}(t),\infty)}$ by (\ref{e.19}),
we conclude that if $\mathcal{U}^{(n)}(t-)\{E^{(n)}(t-)\}>0$,
then $\mathcal{U}^{(n)}(\cdot)\{E^{(n)}(\cdot)\}$
decreases with unit rate in a neighborhood of $t$.
On the other hand, if $\mathcal{U}^{(n)}(t-)\{E^{(n)}(t-)\}=0$,
then since $\Delta A^{(n)}(t)=0$, $E^{(n)}$ jumps up at $t$.
Indeed, in this case,
\[
\sW_S^{(n)}(t)\bigl(-\infty,E^{(n)}(t-)\bigr]=\sW_S^{(n)}(t-)
\bigl(-\infty,E^{(n)}(t-)\bigr]=K^{(n)}(t-)= K^{(n)}(t).
\]
This means that
\begin{eqnarray*}
E^{(n)}(t)
&=&
\min\bigl\{y\in\R|\sW_S^{(n)}(t)(-\infty,y]>K^{(n)}(t)\bigr\}\\
&=&
\min\bigl\{y>E^{(n)}(t-)|\sW_S^{(n)}(t)\{y\}>0\bigr\}.
\end{eqnarray*}
It follows that
%
\begin{equation}\label{4.64a}
\sW_S^{(n)}(t)\bigl(E^{(n)}(t-),E^{(n)}(t)\bigr)=0.
\end{equation}
Using the definition of $E^{(n)}(t)$, (\ref{Uconst}),
(\ref{e.19}),
the assumption $\mathcal{U}^{(n)}(t-)\{E^{(n)}(t-)\}$
$=0$,
the assumption $\bigtriangleup A^{(n)}(t)=0$,
and (\ref{4.64a}), we obtain
\begin{eqnarray*}
\mathcal{U}^{(n)}(t)\bigl[E^{(n)}(t),\infty\bigr)
& = &
U^{(n)}(t)
=
U^{(n)}(t-)\\
& = &
\mathcal{U}^{(n)}(t-)\bigl[E^{(n)}(t-),\infty\bigr)\\
& = &\mathcal{U}^{(n)}(t-)\bigl(E^{(n)}(t-),\infty\bigr)\\
& = &
\mathcal{W}_S^{(n)}(t-)\bigl(E^{(n)}(t-),\infty\bigr)\\
& = & \mathcal{W}_S^{(n)}(t)\bigl(E^{(n)}(t-),\infty\bigr)\\
& = &
\sW_S^{(n)}(t)\bigl[E^{(n)}(t),\infty\bigr).
\end{eqnarray*}
From (\ref{e.19}) we see now that
$\mathcal{U}^{(n)}(t)=\sW_S^{(n)}(t)_{|[E^{(n)}(t),\infty)}$.

(ii) By (\ref{17}), (\ref{E}) and (\ref{calU}),
for $s \in(\sigma^{(n)}_{k}, t)$ we have
%
\begin{equation}\label{e.5}
\mathcal{W}_S^{(n)}(s)\bigl(-\infty,E^{(n)}(s)\bigr]> K^{(n)}(s).
\end{equation}
As $s\uparrow t$ in (\ref{e.5}), by
(\ref{leftlimits}), (\ref{e.11.5}),
and the case (ii)
assumptions $\bigtriangleup A^{(n)}(t)=0$
and $E^{(n)}(t-)=0$, we get
%
\begin{equation}\label{e.5.5}
\mathcal{W}_S^{(n)}(t)(-\infty,0]
=\mathcal{W}_S^{(n)}(t-)(-\infty,0]\geq K^{(n)}(t-).
\end{equation}
When combined with Proposition \ref{P.3}, this implies
%
\begin{equation}\label{e.6}
K^{(n)} (t) = K^{(n)} (t-) \vee
\mathcal{W}_S^{(n)}(t)(-\infty,0]=\mathcal{W}_S^{(n)}(t)(-\infty,0].
\end{equation}
By (\ref{E}) and (\ref{calU}), for $s \in(\sigma^{(n)}_{k}, t)$,
\[
\mathcal{U}^{(n)}(s)\bigl\{E^{(n)}(s)\bigr\}
=\mathcal{W}_S^{(n)}(s)\bigl(-\infty,E^{(n)}(s)\bigr] - K^{(n)}(s).
\]
Letting $s\uparrow t$, invoking (\ref{leftlimits}), (\ref{e.5.5})
and (\ref{e.6}),
and recalling the assumption $E^{(n)} (t-) = 0$, we obtain
%
\begin{equation}\label{e.7}\qquad\quad
\mathcal{U}^{(n)}(t-)\{0\}
=\mathcal{W}_S^{(n)}(t-)(-\infty,0] - K^{(n)}(t-) =K^{(n)}(t)-K^{(n)}(t-),
\end{equation}
and the first equality in (\ref{Cj1}) follows. The second equality
in (\ref{Cj1}) follows from
(\ref{U}), (\ref{e.11.5a}) and the assumption $\bigtriangleup A^{(n)}(t)=0$.
Moreover, by (\ref{calU}) and (\ref{e.6}), for every $y \in\R$,
%
\begin{eqnarray}\label{e.28}
\mathcal{U}^{(n)}(t)(-\infty,y]
&=&
\bigl[ \mathcal{W}_S^{(n)}(t)(-\infty,y]
- \mathcal{W}_S^{(n)}(t)(-\infty,0] \bigr]^+\nonumber\\[-8pt]\\[-8pt]
&=&
\mathcal{W}_S^{(n)}(t)_{|(0,\infty)} (-\infty,y].\nonumber
\end{eqnarray}

(iii) Let $a=E^{(n)}(t-)$. We can deduce (\ref{Cconst})
from (\ref{e.3}) and (\ref{e.4}) as in (i),
with the only difference that now (\ref{e.4}),
for $s<t$ sufficiently close to $t$ such that $A^{(n)}(t-)=A^{(n)}(s)$,
follows from the fact that $L^{(n)}_{A^{(n)}(t)}>0$, since
this implies that the work
for the system associated with the customer arriving to
the system at time $t$ does not contribute to
$\mathcal{W}_S^{(n)}(t)(-\infty,0]$. Next,
let $y<a$, let $\ve=(a-y)/2$ and note that by assumption,
$L^{(n)}_{A^{(n)}(t)}\geq a>y+\ve$.
Thus, for $s<t$, $s$ sufficiently close to $t$
[so as to ensure that $A^{(n)}(t-)=A^{(n)}(s)$], we have
$\mathcal{W}_S^{(n)}(t)(-\infty,y]
\leq\mathcal{W}_S^{(n)}(s)(-\infty,y+\ve] \leq K^{(n)}(s)$,
where the last inequality uses (\ref{calU}) and the fact that
$y + \ve< E^{(n)} (t-)$.
Letting $s\uparrow t$, we obtain
$\mathcal{W}_S^{(n)}(t)(-\infty,y]\leq K^{(n)}(t-)$, which,
together with (\ref{Cconst}), shows that $y< E^{(n)}(t)$.
Thus, $E^{(n)}(t-)\leq E^{(n)}(t)$ or, equivalently,
$\bigtriangleup E^{(n)} (t) \geq0$.

We now turn to the proof of (\ref{e.14}).
Equation (\ref{calU}) implies
%
\begin{equation}\label{e.26}
\mathcal{U}^{(n)}(t-)(-\infty,y]
= \bigl[ \mathcal{W}_S^{(n)}(t-)(-\infty,y]- K^{(n)}(t-)\bigr]^+.
\end{equation}
Indeed, for any $y$ such that $\mathcal{W}_S^{(n)}(t-)\{y\}=0$,
(\ref{e.26}) follows from (\ref{calU}),
in which $t$ is replaced by $s<t$,
by taking $s\uparrow t$. However, the family of sets $(-\infty,y]$
with $\mathcal{W}_S^{(n)}(t-)\{y\}=0$ forms a separating
class in $\sB(\R)$, and
so (\ref{e.26}) holds for all $y$.
Moreover, using (\ref{e.11.5}), (\ref{Cconst})
and (\ref{calU}), we see that
%
\begin{eqnarray}
\label{4.66a}\hspace*{32pt}
&&\mathcal{U}^{(n)} (t) (-\infty, y]\nonumber\\[-8pt]\\[-8pt]
&&\qquad=
\bigl[ \mathcal{W}_S^{(n)} (t-) (-\infty, y] - K^{(n)} (t-)
+ v^{(n)}_{A^{(n)} (t)}
\delta_{L_{A^{(n)}(t)}^{(n)} (t)} (-\infty, y] \bigr]^+.
\nonumber
\end{eqnarray}
When combined with (\ref{e.26}), this shows that
%
\begin{equation}\label{4.66b}
\mathcal{U}^{(n)} (t) (-\infty, y]
= \mathcal{U}^{(n)} (t-) (-\infty, y],\qquad
y < L_{A^{(n)}(t)}^{(n)} (t).
\end{equation}
On the other hand, if $y \geq L_{A^{(n)}(t)}^{(n)} (t)$,
then $y\geq E^{(n)}(t-)$ and (\ref{e.26}) becomes
\[
\mathcal{U}^{(n)}(t-)(-\infty,y]
=\sW_S^{(n)}(t-)(-\infty,y]-K^{(n)}(t-).
\]
From (\ref{4.66a}), we now have
%
\begin{equation}\label{4.66c}\quad
\mathcal{U}^{(n)} (t) (-\infty, y]
= \mathcal{U}^{(n)} (t-) (-\infty, y] + v^{(n)}_{A^{(n)} (t)},\qquad
y\geq L_{A^{(n)}(t)}^{(n)}.
\end{equation}
When combined, (\ref{4.66b}) and (\ref{4.66c}) prove (\ref{e.14}).


(iv) We have $L^{(n)}_{A^{(n)}(t)}>0$, and so (\ref{Cconst}) holds
by the same argument as in case (iii),
but now with $a=L^{(n)}_{A^{(n)}(t)}$. The
assumptions $L^{(n)}_{A^{(n)}(t)} < E^{(n)} (t-)$ and
$\bigtriangleup A^{(n)}(t) = 1$, along with the relations
(\ref{e.11.5}), 
(\ref{calU}), (\ref{Cconst}) and the
definition of $E^{(n)}$, imply that
\begin{eqnarray*}
\mathcal{W}_S^{(n)}(t)\bigl(-\infty,E^{(n)}(t-)\bigr]
&=& \mathcal{W}_S^{(n)}(t-)\bigl(-\infty,E^{(n)}(t-)\bigr]
+ v^{(n)}_{A^{(n)}(t)} \\
& > & \mathcal{W}_S^{(n)}(t-)\bigl(-\infty,E^{(n)}(t-)\bigr] \\
& \geq& K^{(n)} (t-) \\
& = & K^{(n)} (t).
\end{eqnarray*}
Invoking (\ref{calU}) again, this shows that
$\mathcal{U}^{(n)} (t) (-\infty, E^{(n)} (t-)] > 0$, which
implies $E^{(n)}(t)\leq E^{(n)}(t-)$.
Now, let $y<a=L^{(n)}_{A^{(n)}(t)}$ and let $\ve= (a-y)/2$.
Then, combining (\ref{e.11.5}),
the inequalities $y+\ve<a<E^{(n)}(t-)$ and
(\ref{Cconst}), we obtain
\[
\mathcal{W}_S^{(n)}(t)(-\infty,y]
\leq\mathcal{W}_S^{(n)}(t-)(-\infty,y+\ve] \leq K^{(n)}(t-)
=K^{(n)}(t).
\]
This shows that $y<E^{(n)}(t)$, which proves (\ref{e.15}).
In addition, by (\ref{U}) and (\ref{Cconst}), we have
\begin{eqnarray*}
U^{(n)}(t)
&=&
W_S^{(n)}(t)-K^{(n)}(t)\\
&=&
W_S^{(n)}(t- )+v^{(n)}_{A^{(n)}(t)}-K^{(n)}(t-)\\
&=&U^{(n)}(t-)+v^{(n)}_{A^{(n)}(t)},
\end{eqnarray*}
and (\ref{e.15.1}) follows.
Furthermore, since $E^{(n)} (t) \leq E^{(n)}(t-)$ by
(\ref{e.15}), the relations
(\ref{e.19}), (\ref{e.11.5}) and the assumption
$L^{(n)}_{A^{(n)}(t)}<E^{(n)}(t-)$ imply
\begin{eqnarray*}
\mathcal{U}^{(n)}(t)_{|(E^{(n)}(t-),\infty)}
& = &
\mathcal{W}_S^{(n)}(t)_{|(E^{(n)}(t-),\infty)} \\
& = &
\mathcal{W}_S^{(n)}(t-)_{|(E^{(n)}(t-),\infty)}\\
&=&
\mathcal{U}^{(n)}(t-)_{|(E^{(n)}(t-),\infty)}.
\end{eqnarray*}
This establishes (\ref{e.18.5}).

Finally, to prove (\ref{e.18}), we will consider two cases.

\textit{Case} I. $E^{(n)}(t)<E^{(n)}(t-)$.

By (\ref{e.19}), we know that
\[
\mathcal{U}^{(n)}(t)\bigl\{E^{(n)}(t-)\bigr\}
=\mathcal{W}_S^{(n)}(t)\bigl\{E^{(n)}(t-)\bigr\}.
\]
In turn, when combined with (\ref{e.26})
and the definition of $E^{(n)}$, this shows that
\begin{eqnarray*}
&&\mathcal{U}^{(n)}(t-)\bigl\{E^{(n)}(t-)\bigr\}\\
&&\qquad=
\mathcal{U}^{(n)}(t-)\bigl(-\infty, E^{(n)}(t-)\bigr]
- \mathcal{U}^{(n)}(t-)\bigl(-\infty, E^{(n)}(t-)\bigr) \\
&&\qquad =
\mathcal{W}_S^{(n)}(t-)\bigl(-\infty,E^{(n)}(t-)\bigr] -K^{(n)}(t-)\\
&&\qquad\quad{}
- \bigl[\mathcal{W}_S^{(n)}(t-)\bigl(-\infty,E^{(n)}(t-)\bigr)
-K^{(n)}(t-)\bigr]^+\\
&&\qquad\leq
\mathcal{W}_S^{(n)}(t-)\bigl(-\infty,E^{(n)}(t-)\bigr]
-\mathcal{W}_S^{(n)}(t-)\bigl(-\infty,E^{(n)}(t-)\bigr) \\
&&\qquad=
\mathcal{W}_S^{(n)}(t)\bigl\{E^{(n)}(t-)\bigr\} \\
&&\qquad=
\mathcal{U}^{(n)}(t)\bigl\{E^{(n)}(t-)\bigr\},
\end{eqnarray*}
and so (\ref{e.18}) holds.

\textit{Case} II. $E^{(n)}(t)=E^{(n)}(t-)$.

By (\ref{calU}), (\ref{e.11.5}), (\ref{Cconst}),
(\ref{e.26}) and the definition of $E^{(n)}$,
\begin{eqnarray*}
\mathcal{U}^{(n)}(t)\bigl\{E^{(n)}(t)\bigr\}
&=&
\mathcal{U}^{(n)}(t)\bigl(-\infty,E^{(n)}(t)\bigr]\\
&=&
\mathcal{W}_S^{(n)}(t)\bigl(-\infty,E^{(n)}(t)\bigr] -K^{(n)}(t)\\
&=&
\mathcal{W}_S^{(n)}(t-)\bigl(-\infty,E^{(n)}(t-)\bigr]
+v^{(n)}_{A^{(n)}(t)}-K^{(n)}(t-)\\
&=&
\mathcal{U}^{(n)}(t-)\bigl(-\infty,E^{(n)}(t-)\bigr] +v^{(n)}_{A^{(n)}(t)}\\
&=&
\mathcal{U}^{(n)}(t-)\bigl\{E^{(n)}(t-)\bigr\}+v^{(n)}_{A^{(n)}(t)},
\end{eqnarray*}
which establishes (\ref{e.18}) in this case as well.
Since $E^{(n)} (t) \leq E^{(n)} (t-)$,
the two cases above are
exhaustive, and so (\ref{e.18}) is proved.

(v) Equation (\ref{e.28}) holds by the same
argument as in (ii), but where now the
equality in (\ref{e.5.5}) follows from the
fact that $L^{(n)}_{A^{(n)}(t)}>0$.
Let $\mathcal{U}^{(n)}_1(t)
\define\mathcal{U}^{(n)}(t-)+v^{(n)}_{A^{(n)}(t)}
\delta_{L^{(n)}_{A^{(n)}(t)}}
-\mathcal{U}^{(n)}(t-)\{0\} \delta_0$.
We want to show that $\mathcal{U}^{(n)}(t)=\mathcal{U}^{(n)}_1(t)$.
By (\ref{E>0}), $\mathcal{U}^{(n)}(t)$
and $\mathcal{U}^{(n)}(t-)$ are supported on
$(0,\infty)$ and $[0,\infty)$, respectively. Thus,
%
\begin{equation}\label{butlast}
\mathcal{U}^{(n)}(t)(-\infty,y]
=\mathcal{U}^{(n)}_1(t)(-\infty,y]=0,\qquad y\leq0.
\end{equation}
By (\ref{e.19}) and the fact that $E^{(n)}(t-) = 0$,
we have
$\mathcal{U}^{(n)}(t-)_{|(0,\infty)}
= \mathcal{W}^{(n)}_S(t-)_{|(0,\infty)}$.
The last two statements, along with (\ref{e.11.5}),
(\ref{e.28}) and another application of
(\ref{e.19}), show that
\begin{eqnarray*}
\mathcal{U}_1^{(n)}(t)_{|(0,\infty)}
&=&
\mathcal{U}^{(n)}(t-)_{|(0,\infty)}
+v^{(n)}_{A^{(n)}(t)} \delta_{L^{(n)}_{A^{(n)}(t)}} \\
&=&
\mathcal{W}_S^{(n)}(t-)_{|(0,\infty)}
+v^{(n)}_{A^{(n)}(t)} \delta_{L^{(n)}_{A^{(n)}(t)}} \\
&=&
\mathcal{W}_S^{(n)}(t)_{|(0,\infty)} \\
&=&
\mathcal{U}^{(n)}(t)_{|(0,\infty)}.
\end{eqnarray*}
This, together with (\ref{butlast}),
shows that $\mathcal{U}^{(n)}(t)=\mathcal{U}^{(n)}_1(t)$.
\end{pf}

The last result of this section concerns the
evolution of $U^{(n)}$.
Despite the different ways in which arriving mass is distributed
in the system with reneging and the reference system,
in both systems one can keep track of the total mass
in system by beginning with the arrived mass (which
is the same in both systems), subtracting the reduction
in mass due to service (which occurs continuously at unit rate
per unit time whenever mass is present),
and subtracting the mass that
has become late and been deleted. In particular,
a simple mass balance shows that
%
\begin{equation}
\label{WAccount}
W^{(n)}(t)
= V^{(n)}\bigl(A^{(n)}(t)\bigr)
-\int_0^t\ind_{\{W^{(n)}(s)>0\}} \,ds
- \renw(t),
\end{equation}
where we recall that $R_W^{(n)}$ is the total amount of
reneged work in the reneging system,
which admits the representation (\ref{def-renw}),
$\renw(t) = \sum_{0<s\leq t}\sW^{(n)}(s-)\{0\}$,
for all $t \in[0,\infty)$.
We now show that the following analogous relation
holds for the reference workload:
%
\begin{equation}
\label{UAccount}
U^{(n)}(t) =
V^{(n)}\bigl(A^{(n)}(t)\bigr)
-\int_0^t\ind_{\{U^{(n)}(s)>0\}} \,ds
- \renu(t),
\end{equation}
where
%
\begin{equation}
\label{def-renu}
\renu(t) \define\sum_{0<s\leq t}\sU^{(n)}(s-)\{0\}.
\end{equation}
Also, for notational convenience, we
set $\renw(0-) = \renu(0-) = 0$.
\begin{corollary}\label{C.eqn4.9}
For every $t\geq0$, equation (\ref{UAccount}) holds.
Moreover, $\renu= \kpn$ and hence
%
\begin{equation}
\label{cor-UAccount}
U^{(n)} = N^{(n)} + I_U^{(n)} - \kpn,
\end{equation}
where, for $ t \geq0$,
%
\begin{equation}
\label{def-idleu}
I_U^{(n)} (t) \define\int_{0}^t \ind_{\{U^{(n)}(s)=0\}} \,ds.
\end{equation}
\end{corollary}
\begin{pf}
For $t \geq0$, let $\widetilde{U}^{(n)}(t)$
be equal to the right-hand side of
(\ref{UAccount}). By (\ref{U=0}) of Lemma \ref{L.Un},
we have $U^{(n)}(0)=0 = \widetilde{U}^{(n)} (0)$.
Moreover, for every $k \geq1$,
by Lemma \ref{L.Un} and the definition of
$\sigma^{(n)}_k$, it follows that
$U^{(n)} (t-) = U^{(n)} (t) = 0$ and
$\bigtriangleup V^{(n)} (A^{(n)} (t)) = 0$ for
$t \in(\tau_{k-1}^{(n)}, \sigma^{(n)}_k)$,
$U^{(n)} (\sigma^{(n)}_k-) = 0$ and
$\bigtriangleup U^{(n)} (\sigma^{(n)}_k) =
\bigtriangleup V^{(n)} (A^{(n)} (\sigma^{n}_k))$.
When compared with the right-hand side of (\ref{UAccount}),
this shows that $U^{(n)}$ and $\widetilde{U}^{(n)}$ are both
flat on $(\tau_{k-1}^{(n)}, \sigma_{k}^{(n)})$,
with an upward jump at $\sigma_k^{(n)}$ of size
$\bigtriangleup V^{(n)} (A^{(n)} (\sigma_k^{(n)}))$.
Thus, to prove the corollary, it suffices to show that
the increments of $\widetilde{U}^{(n)}$
and $U^{(n)}$ on the intervals
$(\sigma_{k}^{(n)}, \tau_k^{(n)}]$, $k \geq1$,
coincide.

Fix $k\geq1$.
We first show that
%
\begin{equation}\label{deltaU}
\Delta\widetilde{U}^{(n)}\bigl(\tau_k^{(n)}\bigr)=\Delta U^{(n)}\bigl(\tau_k^{(n)}\bigr).
\end{equation}
Equality
(\ref{17}) shows that there cannot be an arrival
at time $\tau_k^{(n)}$, for such an arrival
would have a positive lead time
and hence increase $W_S^{(n)}$ without increasing
$K^{(n)}$ (see Proposition \ref{P.3}). In other words,
$\Delta A^{(n)}(\tau_k^{(n)})=0$.
Because there is no arrival at $\tau_k^{(n)}$,
the measure-valued process $\sW_S^{(n)}$ is continuous
at $\tau_k^{(n)}$. Taking the limit in (\ref{calU})
as $t\uparrow\tau_k^{(n)}$,
we obtain
\[
\sU^{(n)}\bigl(\tau_k^{(n)}-\bigr)(-\infty,0]
=\bigl[\sW_S^{(n)}\bigl(\tau_k^{(n)}\bigr)(-\infty,0]
-K^{(n)}\bigl(\tau_k^{(n)}-\bigr)\bigr]^+
=\Delta K^{(n)}\bigl(\tau_k^{(n)}\bigr),
\]
where the last equality is a consequence of (\ref{20}).
However, (\ref{E>0}) implies that
$\sU^{(n)}(\tau_k^{(n)}-)(-\infty,0)\leq\lim_{t\uparrow\tau_k^{(n)}}
\sU^{(n)}(t)(-\infty,0)=0$,
so $\Delta\widetilde{U}^{(n)}(\tau_k^{(n)})=
-\sU^{(n)}(\tau_k^{(n)}-)\{0\}=-\Delta K^{(n)}(\tau_k^{(n)})$.
From (\ref{U}) and the continuity of $W_S^{(n)}$
at $\tau_k^{(n)}$, we see
that $-\Delta K^{(n)}(\tau_k^{(n)})$ is also equal to
$\Delta U^{(n)}(\tau_k^{(n)})$, and (\ref{deltaU}) is proved.

We next show that $\Delta\widetilde{U}^{(n)}(t)=\Delta U^{(n)}(t)$
for $t \in(\sigma_{k}^{(n)}, \tau_k^{(n)})$.
If $E^{(n)} (t-) > 0$, then the definitions of $E^{(n)}$ and
$\widetilde{U}^{(n)}$, and
statements (i), (iii) and (iv) of Proposition~\ref{P.Un} show that
\[
\bigtriangleup U^{(n)} (t) = \bigtriangleup\widetilde{U}^{(n)} (t)
= \bigtriangleup V^{(n)} \bigl(A^{(n)} (t)\bigr).
\]
On the other hand, if $E^{(n)} (t-) = 0$,
then properties (ii) and (v) of Proposition
\ref{P.Un} and the definition of $\widetilde{U}^{(n)}$ show that
\[
\bigtriangleup U^{(n)} (t)
= \bigtriangleup\widetilde{U}^{(n)} (t)
= \bigtriangleup A^{(n)} (t) v_{ A^{(n)} (t)}^{(n)}
- \mathcal{U}^{(n)} (t-) \{0\}.
\]

Now, let $S^{(n)}$ be the (random) set of
times $s\geq0$
for which $U^{(n)}(s)>0$ and at least one of the
following three properties holds:
\begin{eqnarray*}
\bigtriangleup A^{(n)}(s)&>&0,\\
E^{(n)}(s-)&=&0
\end{eqnarray*}
or
\[
\mathcal{U}^{(n)}(s-)\bigl\{E^{(n)}(s-)\bigr\}=0.
\]
Suppose $U^{(n)} (s) > 0$.
If $E^{(n)} (s-) = 0$, then the fact that $E^{(n)}(s)>0$ by
(\ref{E>0}) implies $\bigtriangleup E^{(n)} (s) > 0$, while if
$\mathcal{U}^{(n)} (s-) \{E^{(n)} (s-)\} = 0$, the definition of
$E^{(n)}(s)$ implies that
$\bigtriangleup(\mathcal{U}^{(n)} (s) \{E^{(n)} (s)\}) > 0$.
Thus, the set $S^{(n)}$ is countable, and on the
set $\{ s \in(\sigma_k^{(n)},t)\dvtx U^{(n)} (s) > 0 \} \setminus S^{(n)}$,
the process $U^{(n)}$ decreases with unit rate by
Proposition \ref{P.Un}(i). Therefore, the total amount of
this decrease on any time interval of the
form
$(\sigma_k^{(n)}, t)$ equals
$\int_{\sigma_k^{(n)}}^t \ind_{\{U^{(n)} (s) > 0\}} \,ds$,
which\vspace*{-2pt}
coincides with the absolutely continuous part of
$\widetilde{U}^{(n)}(t)-\widetilde{U}^{(n)}(\sigma_k^{(n)}-)$
on the same interval. This concludes the proof of
(\ref{UAccount}).

Adding and subtracting $t$ to (\ref{UAccount}),
by the definition
(\ref{2.15}) of the netput process $N^{(n)}$ and the nonnegativity
of $U^{(n)}$, we obtain
%
\begin{equation}\label{4.72a}
U^{(n)}(t) = N^{(n)} (t) + \int_0^t\ind_{\{U^{(n)}(s) = 0\}} \,ds
- R_U^{(n)} (t),
\end{equation}
while substituting (\ref{eq-Kdcomp})
and (\ref{2.18}) into (\ref{U}), we have
\[
U^{(n)} (t) = N^{(n)} (t) + I_S^{(n)} (t) + \knn(t) - \kpn(t)
\]
for $t \geq0$. On the other hand, we know that
\[
\int_{[0,\infty)} \ind_{\{U^{(n)}(s) > 0\}} \,dI_S^{(n)} (s)
= 0 \quad\mbox{and}\quad
\int_{[0,\infty)} \ind_{\{U^{(n)}(s) > 0\}} \,d\knn(s) = 0,
\]
where the former equality holds because
$W_S^{(n)} \geq U^{(n)}$ by (\ref{eq-uws}), and
$I_S^{(n)}$ increases only at times when $W_S^{(n)}$ is zero,
while the latter holds by (\ref{complementarity}).
From the last three displays, we conclude that
%
\begin{equation}\label{4.72b}
\int_{[0,\infty)} \ind_{\{U^{(n)}(s) > 0\}} \,d \renu(s)
= \int_{[0,\infty)} \ind_{\{U^{(n)}(s) > 0\}} \,d \kpn(s).
\end{equation}
On the other hand,
since $U^{(n)} (s) = 0$ implies $\Delta A^{(n)} (s) =0$, from
properties (i) and (ii) of Proposition \ref{P.Un} and the fact that
$\renu$ is a pure jump process with
$\Delta\renu(t) = \mathcal{U}^{(n)} (t-) \{0\}$,
it follows that
\[
\int_{[0,\infty)} \ind_{\{U^{(n)}(s) = 0\}} \,d \renu(s) = \int
_{[0,\infty)} \ind_{\{U^{(n)}(s) = 0\}} \,d \kpn(s).
\]
Together, the last two equalities imply
$\renu=\kpn$, which, when substituted into (\ref{UAccount}),
yields (\ref{cor-UAccount}).
\end{pf}

\section{The reneging system} \label{S5a}

In this section
we bound the difference in workload between the
pre-limit reference and reneging systems---Lemma \ref{L.step2}
provides a lower bound, while Lemma \ref{L5.3} provides
an upper bound. The proof of the upper bound uses
an optimality property of EDF that may be of independent
interest.
\begin{theorem}\label{t.EDFopt}
Let $\pi$ be a service policy for a single-station,
single-customer-class
queueing system with reneging
such that the customer arrival times to this system
do not have a finite accumulation point.
Let $R_\pi(t)$ be the amount of work removed from
this system up to time $t$ due to lateness.
Let $R_W (t)$ be the amount of work removed
due to lateness up to time $t$ from the EDF system
with reneging and the same
interarrival times, service times and lead times
as in the former system.
Then for every $t\geq0$, we have
%
\begin{equation}\label{e.EDFopt}
R_W (t)\leq R_\pi(t).
\end{equation}
\end{theorem}

The proof of Theorem \ref{t.EDFopt} is deferred to the \hyperref[app]{Appendix}.
The related fact
that the EDF protocol minimizes the number of late customers in the
$G/M/c$ queue was proved in \cite{pat}, and the main idea of
our proof is similar to that of \cite{pat}. However,
our argument is pathwise and the only
assumption on the distribution of the system stochastic primitives
that we impose is that customer arrivals do not have a finite
accumulation point.
This assumption is clearly satisfied almost
surely by a $GI/G/1$ queue.

\subsection{Comparison results}
\label{subs-comparison}

In this section, we establish bounds on the
difference between the processes $U^{(n)}$ and $W^{(n)}$.
In Section \ref{subs-htwork},
this difference will be shown to be
negligible in the heavy traffic limit.
We start with Lemma \ref{L.step2}
showing that $W^{(n)}\leq U^{(n)}$, which
implies that
$\renu\leq\renw$ (see Corollary \ref{c.latecomp}).

In the proofs of these results, we will make
frequent use of the observation
that, by (\ref{WAccount}) and (\ref{UAccount}),
%
\begin{eqnarray}
\label{wudiff}
W^{(n)} (t) -U^{(n)} (t)
&=& \int_0^t \ind_{\{U^{(n)}(s) > 0\}} \,ds
-\int_0^t \ind_{\{W^{(n)} (s) > 0\}} \,ds\nonumber\\[-8pt]\\[-8pt]
&&{} + \renu(t) - \renw(t)\nonumber
\end{eqnarray}
for $t \in[0,\infty)$.
\begin{lemma}\label{L.step2}
For every $t\geq0$, we have
%
\begin{equation}\label{e.step2}
W^{(n)}(t)\leq U^{(n)}(t).
\end{equation}
\end{lemma}
\begin{pf}
Let
%
\begin{equation}\label{tau61}
\tau\define\min\bigl\{t\geq0\dvtx W^{(n)}(t)> U^{(n)}(t)\bigr\}.
\end{equation}
If $\tau=+\infty$, then (\ref{e.step2}) holds.
Assume $\tau<+\infty$. In this case,
we claim that the minimum on the right-hand
side of (\ref{tau61}) is attained.
Indeed, (\ref{wudiff}) and the fact that
$\renu$ and $\renw$ are pure jump
processes show that the only way that $W^{(n)} - U^{(n)}$ can
become strictly positive is via a jump.
Thus $W^{(n)} (\tau) > U^{(n)} (\tau)$.
Since $\mathcal{W}(\tau) (-\infty, 0]
= \sU(\tau) (-\infty, 0] = 0$ (in fact, this equality
holds for any time $t$),
this means there must exist a $y>0$ such that
%
\begin{equation}\label{y61}
\mathcal{W}^{(n)}(\tau)(y,\infty)> \mathcal{U}^{(n)}(\tau
)(y,\infty).
\end{equation}
Let
%
\begin{equation}\label{tau061}\quad
\tau_0 \define\inf\bigl\{t \in[0,\tau]\dvtx
\mathcal{W}^{(n)}(t)(y+\tau-t,\infty)
> \mathcal{U}^{(n)}(t)(y+\tau-t,\infty)\bigr\}.
\end{equation}
By (\ref{y61}), the above infimum is over a nonempty set.
Lemma \ref{L.Un} and Proposition~\ref{P.Un} imply that the
only difference in the dynamics of $\mathcal{W}^{(n)}$ and
$\mathcal{U}^{(n)}$ is that the arriving mass $v^{(n)}_k$
is concentrated at $L^{(n)}_k$ in the case of the EDF
system with reneging and
distributed in $[L^{(n)}_k,\infty)$ in the reference
system. On the other hand, in both systems at time
$t\in[0,\tau]$, no mass leaves the
interval $(y+\tau-t,\infty)$ due to lateness.
This implies that
$\mathcal{W}^{(n)}(t)(y+\tau-t,\infty)
- \mathcal{U}^{(n)}(t)(y+\tau-t,\infty)$, $t\in[0,\tau]$,
has no positive jumps and therefore
%
\begin{equation}\label{attau061}
\mathcal{W}^{(n)}(t)(y+\tau-\tau_0,\infty)
= \mathcal{U}^{(n)}(t)(y+\tau-\tau_0,\infty).
\end{equation}
By (\ref{y61}) and (\ref{attau061}), $\tau_0<\tau$. Thus,
there exists $t\in(\tau_0,\tau)$, where $t-\tau_0$ is
arbitrarily small and
%
\begin{equation}\label{261}
\mathcal{W}^{(n)}(t)(y+\tau-t,\infty)
> \mathcal{U}^{(n)}(t)(y+\tau-t,\infty).
\end{equation}
However, we claim that (\ref{attau061}) and (\ref{261})
imply that for all $t\in(\tau_0,\tau)$, where $t-\tau_0$
is small enough, it must be that
%
\begin{eqnarray}
\label{361}
\mathcal{W}^{(n)}(t)(0,y+\tau-t]&>&0,\\
\label{461}
\mathcal{U}^{(n)}(t)(0,y+\tau-t]&=&0.
\end{eqnarray}
Indeed, if (\ref{361}) is false, then the left-hand side of
(\ref{261}) is equal to $W^{(n)} (t)$,
and consequently decreases with unit speed as long as
it is nonzero in
some time interval beginning with $\tau_0$.
Similarly, if (\ref{461}) is false,
the right-hand side of (\ref{261})
is constant on some interval beginning with $\tau_0$.
In both cases, due to (\ref{attau061}), (\ref{261}) cannot hold for
$t\in(\tau_0,\tau)$ with $t-\tau_0$ arbitrarily small.
But (\ref{261})--(\ref{461}) yield $W^{(n)}(t)> U^{(n)}(t)$ for some
$t<\tau$, which contradicts (\ref{tau61}).
\end{pf}
\begin{corollary}\label{c.latecomp}
For every $t\geq0$,
%
\begin{equation}\label{c621}
\renu(t)\leq\renw(t).
\end{equation}
Moreover, for $k\geq1$ and $t\geq\sigma^{(n)}_k$,
%
\begin{equation}\label{c622}
\renu(t)- \renu\bigl(\sigma^{(n)}_k-\bigr)\leq\renw(t)- \renw\bigl(\sigma^{(n)}_k-\bigr).
\end{equation}
\end{corollary}
\begin{pf}
Lemma \ref{L.step2} and (\ref{wudiff}) imply that for $0 \leq s \leq t$,
%
\begin{equation}
\label{uwdiff2}
\bigl( \renu(t) - \renu(s ) \bigr) -
\bigl( \renw(t) -\renw(s ) \bigr)
\leq U^{(n)} (s) - W^{(n)} (s).
\end{equation}
Substituting $s = 0$ into (\ref{uwdiff2}) and using the fact that
$\renu(0) = \renw(0) = U^{(n)} (0) = W^{(n)} (0) =0$,
we obtain (\ref{c621}).
Likewise, for $0 \leq s \leq\sigma_k^{(n)} \leq t$,
taking limits as $s$ tends to $\sigma_k^{(n)}$ in (\ref{uwdiff2}),
and using the fact that $U^{(n)} ( \sigma_k^{(n)} -) = W^{(n)} (\sigma
_k^{(n)} -) = 0$,
which follows from (\ref{U=0}), Lemma \ref{L.step2} and the
nonnegativity of
$W^{(n)}$, we obtain (\ref{c622}).
\end{pf}

The proofs of Lemma \ref{L.step2} and
Corollary \ref{c.latecomp}
show the following more general (and intuitively obvious) fact:
if all customers in the EDF system with reneging get larger deadlines,
this results in a
larger workload at every time $t$ and a smaller total amount of mass
removed from the system
due to lateness in the time interval $[0,t]$.

We now establish an inequality between the frontiers in both systems.
\begin{lemma}\label{L.fors3}
For every $t\geq0$ such that $U^{(n)}(t)>0$, we have
%
\begin{equation}\label{e.fors3}
E^{(n)}(t)\leq F^{(n)}(t).
\end{equation}
\end{lemma}
\begin{pf}
Subtracting (\ref{calU}) from (\ref{U}), we see that for any $y \in
\R$,
%
\begin{eqnarray}\label{use1}\qquad
\mathcal{U}^{(n)} (t) (y, \infty) & = & W_S^{(n)} (t) - K^{(n)} (t) -
\bigl[ \mathcal{W}_S^{(n)} (t) (-\infty, y] - K^{(n)} (t) \bigr]^+
\nonumber\\[-8pt]\\[-8pt]
& \leq&
\mathcal{W}_S^{(n)} (t) (y, \infty).\nonumber
\end{eqnarray}
Now, assume that for some $t$ we have $F^{(n)}(t)<E^{(n)}(t)$. In this case,
%
\begin{eqnarray} \label{W>=U}
W^{(n)}(t)&\geq& \mathcal{W}^{(n)}(t)\bigl\{C^{(n)}(t)\bigr\} + \mathcal
{W}^{(n)}(t)\bigl(F^{(n)}(t),\infty\bigr)\nonumber\\
&=& \mathcal{W}^{(n)}(t)\bigl\{C^{(n)}(t)\bigr\} + \mathcal
{V}^{(n)}(t)\bigl(F^{(n)}(t),\infty\bigr)\nonumber\\
&\geq& \mathcal{W}^{(n)}(t)\bigl\{C^{(n)}(t)\bigr\} + \mathcal
{V}^{(n)}(t)\bigl[E^{(n)}(t),\infty\bigr)\nonumber\\[-8pt]\\[-8pt]
&\geq& \mathcal{W}^{(n)}(t)\bigl\{C^{(n)}(t)\bigr\} + \mathcal
{W}_S^{(n)}(t)\bigl[E^{(n)}(t),\infty\bigr)\nonumber\\
&\geq& \mathcal{W}^{(n)}(t)\bigl\{C^{(n)}(t)\bigr\} + \mathcal
{U}^{(n)}(t)\bigl[E^{(n)}(t),\infty\bigr)\nonumber\\
&\geq& U^{(n)}(t), \nonumber
\end{eqnarray}
where the second line follows from the fact that none of the customers
in the EDF
system with reneging that have lead times at time $t$ greater than
$F^{(n)}(t)$ has received any service up to time $t$, the second-last
inequality follows from (\ref{use1}) and the last
line holds due to the equality $U^{(n)}(t)=\mathcal
{U}^{(n)}(t)[E^{(n)}(t),\infty)$.
When combined with the assumption that $U^{(n)} (t) > 0$, this implies that
$W^{(n)}(t)>0$.
This, in turn, implies that $\mathcal{W}^{(n)}(t)\{C^{(n)}(t)\} > 0$ because
the residual service time of the currently served customer is strictly positive.
Thus, the last inequality in (\ref{W>=U}) is strict, which contradicts
(\ref{e.step2}).
\end{pf}

Let $D^{(n)}(t)$ be the amount of work deleted by the EDF system
with reneging
in the time interval $[0,t]$ that is associated with customers
whose lead times upon arrival
were smaller than the
value of the frontier at the time of their arrival.
In the proof of the next lemma,
we will make use of the elementary fact that
by the definition of
$F^{(n)}$ we have
%
\begin{equation}\label{frontdyn}
F^{(n)}(t_1)-(t_2-t_1)\leq F^{(n)}(t_2),\qquad
S_1^{(n)}\leq t_1\leq t_2.
\end{equation}
\begin{lemma}\label{stepI}
For every $t\geq0$,
%
\begin{equation}\label{e.De0}
U^{(n)}(t) - W^{(n)}(t)\leq D^{(n)}(t).
\end{equation}
\end{lemma}
\begin{pf}
If $t\in[\tau^{(n)}_{k-1},\sigma^{(n)}_k)$
for some $k\geq1$, then $U^{(n)}(t)=0$ by (\ref{U=0}).
Thus, by (\ref{6}), it suffices to prove
(\ref{e.De0}) on $[\sigma^{(n)}_k,\tau^{(n)}_k)$ for every $k\geq1$.
Let $k\geq1$.
Suppose that (\ref{e.De0}) is false for some
$t\in[\sigma^{(n)}_k,\tau^{(n)}_k)$. Let
%
\begin{equation}\label{e.Dtau}
\tau\define\min\bigl\{t\in\bigl[\sigma^{(n)}_k,\tau^{(n)}_k\bigr)
|U^{(n)}(t) - W^{(n)}(t)> D^{(n)}(t)\bigr\}.
\end{equation}
We first argue that the minimum on the right-hand side
of (\ref{e.Dtau}) is attained. Indeed, by (\ref{wudiff}) and
Lemma \ref{L.step2}, it is clear that $U^{(n)} - W^{(n)}$
cannot increase except by a jump that is
due to lateness in the EDF system with reneging.
Thus, we have $\mathcal{W}^{(n)}(\tau-)\{0\}>0$ and
%
\begin{equation}\label{e.attau}
U^{(n)}(\tau) - W^{(n)}(\tau)>D^{(n)}(\tau).
\end{equation}
Also, (\ref{U=0}), (\ref{massj}) and Lemma \ref{L.step2} imply that
$U^{(n)}(\sigma^{(n)}_k) =
\bigtriangleup U^{(n)}(\sigma^{(n)}_k)
= \bigtriangleup W^{(n)}(\sigma^{(n)}_k)
=W^{(n)}(\sigma^{(n)}_k)$, 
so $\sigma^{(n)}_k<\tau$.
In particular, (\ref{e.Dtau}) implies
%
\begin{equation}\label{e.beforetau}
U^{(n)}(\tau-) - W^{(n)}(\tau-)\leq D^{(n)}(\tau-).
\end{equation}

Let $k_0$ be the index of the customer arriving
at time $\sigma_k^{(n)}$, that is,
$S^{(n)}_{k_0}=\sigma^{(n)}_k$.
Let $k_1\geq k_0$ be the index of a customer who reneges
in the reneging system at time $\tau$. There must be such
a customer, and there may in fact be more than one such customer.
The amount of work associated with all such customers at
time $\tau$ is $\sW^{(n)}(\tau-)\{0\}$, and we seek to
show that this work is bounded above by
$\bigtriangleup D^{(n)}(\tau)$.
We have $S_{k_1}^{(n)}\in[\sigma_k^{(n)},\tau)$
and $L_{k_1}^{(n)}-(\tau-S_{k_1}^{(n)})=0$.
The subsequent analysis is divided into two cases.

\textit{Case} I.
For every customer $k_1$ chosen as just described, assume
there is a customer $\ell$ arriving in the time interval
$[\sigma_k^{(n)},S_{k_1}^{(n)}]$ who is at least
as urgent as customer $k_1$ when customer $k_1$ arrives
but whose associated mass in the reference system
is at least partly assigned so that upon the arrival
of customer $k_1$, this mass is to the right of $L_{k_1}^{(n)}$.
In other words, $\ell\in[k_0,k_1]$,
$L_{\ell}^{(n)}-(S_{k_1}^{(n)}-S_{\ell}^{(n)})
\leq L_{k_1}^{(n)}$ and
$\bigtriangleup\mathcal{W}^{(n)}(S^{(n)}_\ell)\{L^{(n)}_\ell\}
> \bigtriangleup\mathcal{U}^{(n)}(S^{(n)}_\ell)
[L^{(n)}_l,L^{(n)}_{k_1}+S^{(n)}_{k_1}-S^{(n)}_{\ell}]$.
In this case,
$\bigtriangleup\mathcal{U}^{(n)}(S^{(n)}_\ell)
(L^{(n)}_{k_1}+S^{(n)}_{k_1}-S^{(n)}_\ell,\infty)>0$.
Indeed, by Lemma \ref{L.Un} and Proposition
\ref{P.Un}(iv) (describing the only case in which part of the
mass of a new customer is
distributed by the reference workload to a point
other than its lead time)
$\bigtriangleup U(S^{(n)}_\ell)=v^{(n)}_\ell$ and
$\bigtriangleup\mathcal{U}^{(n)}(S^{(n)}_\ell)
(-\infty,L^{(n)}_\ell)=0$
[see (\ref{massj}), (\ref{shiftr}), (\ref{e.15}), (\ref{e.15.1})
and (\ref{E})]. Let
$s>L^{(n)}_{k_1}+S^{(n)}_{k_1}-S^{(n)}_\ell$
satisfy $\bigtriangleup\mathcal{U}^{(n)}(S^{(n)}_\ell)\{s\}>0$.
Such a point $s$ exists since the measure
$\mathcal{U}^{(n)}(S^{(n)}_\ell)$ is discrete.

If $\ell>k_0$ (e.g., $\ell=k_1$), then, by
(\ref{e.18.5}) in Proposition \ref{P.Un}(iv)
and Lemma \ref{L.fors3}, we have
$s\leq E^{(n)}(S^{(n)}_\ell-)
\leq F^{(n)}(S^{(n)}_\ell-)\leq F^{(n)}(S^{(n)}_\ell)$.
Thus, by (\ref{frontdyn}),
$L^{(n)}_{k_1}<s-(S^{(n)}_{k_1}-S^{(n)}_\ell)
\leq F^{(n)}(S^{(n)}_\ell)-(S^{(n)}_{k_1}-S^{(n)}_\ell)
\leq F^{(n)}(S^{(n)}_{k_1})$.

If $\ell=k_0$, then, because
$\mathcal{U}^{(n)}(S^{(n)}_{k_0})\{s\}>0$,
we have $\mathcal{W}_S^{(n)}(S^{(n)}_{k_0})\{s\}>0$ by the
definition of $\mathcal{U}^{(n)}$. In this case
$\mathcal{W}^{(n)}(S^{(n)}_{k_0})\{s\}=0$, because
$W^{(n)}\equiv0$ on $[\tau^{(n)}_{k-1},\sigma^{(n)}_k)$
by (\ref{U=0}) and Lemma \ref{L.step2}, so
$\mathcal{W}^{(n)}(S^{(n)}_{k_0})
=\mathcal{W}^{(n)}(\sigma^{(n)}_k)
=v^{(n)}_{k_0}\delta_{L^{(n)}_{k_0}}$ and
$s>L^{(n)}_{k_1}+S^{(n)}_{k_1}-S^{(n)}_{k_0}\geq L^{(n)}_{k_0}$
by the definitions of $\ell$ and $s$. Thus, a customer with
lead time equal to $s$
at time $S^{(n)}_{k_0}$ has already been in service in the
EDF system with reneging, so
$L^{(n)}_{k_1}+S^{(n)}_{k_1}-S^{(n)}_{k_0}<s
\leq F^{(n)}(S^{(n)}_{k_0})$
and consequently, by (\ref{frontdyn}),
$L^{(n)}_{k_1}
< F^{(n)}(S^{(n)}_{k_0})-(S^{(n)}_{k_1}-S^{(n)}_{k_0})
\leq F^{(n)}(S^{(n)}_{k_1})$.

Thus, regardless of the value of
$\ell$, $L^{(n)}_{k_1}<F^{(n)}(S^{(n)}_{k_1})$.
In other words, under the Case I assumption,
every customer $k_1$ who becomes
late at time $\tau$ in the EDF system with reneging
arrived with initial
lead time smaller than the value of $F^{(n)}$ at
the time of its arrival. The work associated with
these customers deleted at time $\tau$ is
$\bigtriangleup D^{(n)}(\tau)$.
We conclude that
$\sW^{(n)}(\tau-)\{0\}=\bigtriangleup D^{(n)}(\tau)$.
However, by (\ref{wudiff}), we have
$\bigtriangleup(U^{(n)}-W^{(n)})(\tau)
\leq\sW^{(n)}(\tau-)\{0\}$,
and so $\bigtriangleup(U^{(n)}-W^{(n)})(\tau)
\leq\bigtriangleup D^{(n)}(\tau)$.
This, together with (\ref{e.beforetau}), contradicts
(\ref{e.attau}).

\textit{Case} II.
For a customer $k_1$ chosen as described above,
assume that every customer $\ell$ arriving in
the time interval $[\sigma_k^{(n)},S_{k_1}^{(n)}]$ who is
as least as urgent as customer $k_1$ when customer $k_1$
arrives has all its associated mass initially assigned
in the reference system
to the interval $(0,L_{k_1}^{(n)}+S_{k_1}^{(n)}-S_{\ell}^{(n)}]$
upon arrival.
Customers $\ell$ who are less urgent then $k_1$ must have
lead times satisfying
$L_{\ell}^{(n)}>L_{k_1}^{(n)}+S_{k_1}^{(n)}-S_{\ell}^{(n)}$,
and hence the mass brought by such customers must be initially
assigned to the half-line
$(L_{k_1}^{(n)}+S_{k_1}^{(n)}-S_{\ell}^{(n)},\infty)$
in both systems. Then for every
$t\in[\sigma^{(n)}_{k},S^{(n)}_{k_1}]$, we have
%
\begin{equation}\label{e.!}
\mathcal{W}^{(n)}(t)\bigl(0,L^{(n)}_{k_1}-\bigl(t-S^{(n)}_{k_1}\bigr)\bigr]
\leq\mathcal{U}^{(n)}(t)\bigl(0,L^{(n)}_{k_1}-\bigl(t-S^{(n)}_{k_1}\bigr)\bigr],
\end{equation}
as we now explain.
Under the Case II assumption the arrival of new
mass is the same on both sides of (\ref{e.!}).
Furthermore, disregarding lateness and new arrivals,
both sides of (\ref{e.!}) decrease at unit rate
so long as they are nonzero. Finally, by (\ref{c622})
the amount of late work removed from the EDF system with reneging
in the time interval $[\sigma_k^{(n)},t]$ is greater
than or equal to the amount of late work removed from
$\sU^{(n)}$ in this time interval. Therefore,
(\ref{e.!}) holds for
every $t\in[\sigma_k^{(n)}, S_{k_1}^{(n)}]$.

We claim that (\ref{e.!}) in fact
holds for all $t\in[\sigma_k^{(n)},\tau)$.
Suppose this is not the case. Let
%
\begin{eqnarray}\label{e.Deta}
&&\eta\define\inf\bigl\{t\in\bigl[S_{k_1}^{(n)},\tau\bigr)|
\mathcal{W}^{(n)}(t)\bigl(0,L^{(n)}_{k_1}-\bigl(t-S^{(n)}_{k_1}\bigr)\bigr]\nonumber\\[-8pt]\\[-8pt]
&&\hspace*{92.3pt}
>\mathcal{U}^{(n)}(t)\bigl(0,L^{(n)}_{k_1}-\bigl(t-S^{(n)}_{k_1}\bigr)\bigr]\bigr\}.\nonumber
\end{eqnarray}
The strict inequality in (\ref{e.Deta}) can occur only
because of an arrival at time $t$ which brings mass
to the interval $(0,L_{k_1}^{(n)}-(t-S_{k_1}^{(n)})]$
under the $\sW^{(n)}$ measure but not under the $\sU^{(n)}$
measure. The arrival at time $k_1$ does not have this property
because the Case II assumption applies to $\ell=k_1$.
Therefore, $\eta>S_{k_1}^{(n)}$.

Also, for $t\in[S^{(n)}_{k_1},\tau)$,
%
\begin{equation}\label{e.lhs>0}
\mathcal{W}^{(n)}(t)\bigl\{L^{(n)}_{k_1}-\bigl(t-S^{(n)}_{k_1}\bigr)\bigr\} > 0,
\end{equation}
because the customer $k_1$ is present in the EDF system
with reneging at time $t$. By (\ref{E}),
(\ref{e.lhs>0}) and the definition of $\eta$, we have
$E^{(n)}(t)\leq L^{(n)}_{k_1}-(t-S^{(n)}_{k_1})$
for $t\in[S^{(n)}_{k_1},\eta)$.
Thus,
$E^{(n)}(t-)\leq L^{(n)}_{k_1}-(t-S^{(n)}_{k_1})$
for $t\in(S^{(n)}_{k_1},\eta]$.\vspace*{1pt}
We argue that this implies that the
amounts of mass arriving in both the
EDF system with reneging and the
reference workload\vspace*{-2pt} at any time $t\in(S^{(n)}_{k_1},\eta]$
with lead times upon arrival less than or equal to
$L^{(n)}_{k_1}-(t-S^{(n)}_{k_1})$ are the same.
Indeed, Proposition \ref{P.Un}, especially (\ref{e.18.5}),
implies that no mass arriving at
time $t$ with lead time smaller than $E^{(n)}(t-)$ in
the EDF system with reneging
is distributed to lead times greater than $E^{(n)}(t-)$
by the reference workload. Also, Proposition
\ref{P.Un}(iii) and (v) imply that the mass arriving at
time $t$ with lead time greater than
or equal to $E^{(n)}(t-)$ is distributed in the same
way by the EDF system with reneging and the reference system.
By the same argument as in the case of
$t\in[\sigma^{(n)}_{k},S^{(n)}_{k_1}]$, we conclude that
(\ref{e.!}) holds for $t\in[S^{(n)}_{k_1},\eta]$,
which contradicts the
definition of $\eta$. We have shown that (\ref{e.!})
holds for $t\in[\sigma^{(n)}_{k},\tau)$.

Letting\vspace*{1pt} $t\uparrow\tau$ in
(\ref{e.!}) and using the fact that
$L^{(n)}_{k_1}-(\tau-S^{(n)}_{k_1})=0$, we get
$\mathcal{W}^{(n)}(\tau-)\{0\} \leq\mathcal{U}^{(n)}(\tau-)\{0\}$.
Thus, by (\ref{wudiff}),
$\bigtriangleup(U^{(n)} - W^{(n)})(\tau)
= \mathcal{W}^{(n)}(\tau-)\{0\} - \mathcal{U}^{(n)}(\tau-)\{0\}
\leq0$ which,
together with (\ref{e.beforetau}) and the fact
that $D^{(n)}$ is nondecreasing, contradicts (\ref{e.attau}).
\end{pf}

For the sake of the next proof, we define a sequence of auxiliary \textit{hybrid 
systems} (with the same stochastic primitives as in the
case of the EDF systems described in Section \ref{S2.1}) as follows.
The hybrid system gives priority to the jobs whose lead times upon
arrival are smaller
than the current frontier $F^{(n)}$ in the corresponding EDF system
with reneging. In other words, for each $k$, the
$k$th customer arriving at the hybrid system joins the high-priority
class if and only if
%
\begin{equation}\label{highp}
L^{(n)}_k< F^{(n)} \bigl(S^{(n)}_k\bigr).
\end{equation}
The system processes high-priority customers according to the FIFO
service discipline. When the priority class empties, the system goes
idle until
either another high-priority customer arrives and the system resumes
service in the manner described above, or the corresponding
EDF system with reneging finishes serving the~customers who have
received priority in the hybrid system.
Here, we are using the fact that the high-priority customers leave the
hybrid system
before they leave the EDF system with reneging, which is a consequence
of the optimality
of the EDF discipline established in Theorem \ref{t.EDFopt}.
(We have slightly abused the terminology here, identifying the $k$th
customer in the hybrid system with the corresponding customer from the
EDF system with reneging, while,
formally, only the random variables $u^{(n)}_k$, $v^{(n)}_k$ and
$L^{(n)}_k$ associated with these customers are the same.)
Whenever the EDF system with reneging finishes serving
a batch of customers who have received high priority in the
hybrid system, both systems then serve
the low-priority class using the
EDF discipline until the next high-priority customer arrives.
In both systems, if a customer is present when his
deadline passes, he leaves the queue immediately,
regardless of his class.
The measure-valued workload process associated with the hybrid
system will be denoted by $\mathcal{W}^{(n)}_H$.
\begin{lemma}\label{L5.3}
For every $t\geq0$, we have
%
\begin{eqnarray}
\label{e.step3c}
&&U^{(n)}(t) - W^{(n)}(t)\nonumber\\
&&\qquad\leq
\sum_{k=1}^{A^{(n)}(t)}
v^{(n)}_k \wedge
\bigl( \mathcal{W}^{(n)}\bigl(S^{(n)}_k -\bigr)\bigl(0,F^{(n)}\bigl(S_k^{(n)}\bigr)\bigr)
+ v^{(n)}_k -L^{(n)}_k \bigr)^+\\
&&\qquad\quad\hspace*{53.2pt}{}\times
\ind_{\{L^{(n)}_k< F^{(n)}(S_k^{(n)})\}} .\nonumber
\end{eqnarray}
\end{lemma}
\begin{pf}
By Lemma \ref{stepI}, it suffices to show that $D^{(n)}(t)$ is not
greater than the right-hand side of (\ref{e.step3c}).
By Theorem \ref{t.EDFopt}, $D^{(n)}(t)$, the amount of unfinished work
associated with customers who arrived with lead times smaller than
$F^{(n)}$ and
were deleted in the time interval $[0,t]$ by the EDF system with
reneging, is not greater than the unfinished work associated with these
customers and
deleted by the corresponding hybrid system.
Note that the customers with lead times satisfying (\ref{highp})
form a priority class in both the
EDF system with reneging and the hybrid system, and so their
service is not affected by the presence of other customers.
Furthermore, unfinished work associated with deleted
customers who arrived with lead times greater than
or equal to $F^{(n)}$ is the same in both systems.

For each $k$, if (\ref{highp}) holds,
then the $k$th customer of the hybrid system belongs to the
high-priority class.
Moreover, if, for some $l<k$, $L^{(n)}_l< F^{(n)} (S^{(n)}_l)$, then,
by (\ref{frontdyn}), $L^{(n)}_l- (S^{(n)}_k-S^{(n)}_l)$,
the lead time of the $l$th customer at time $S^{(n)}_k$, is smaller
than $F^{(n)} (S^{(n)}_k)$.
Thus, if (\ref{highp}) holds, the $k$th customer waits at most
$\mathcal{W}^{(n)}_H(S^{(n)}_k-)(0,F^{(n)} (S^{(n)}_k) )$ time units
before he starts receiving service. (His waiting time may actually be
smaller because
some of the high-priority customers in queue who have arrived before
him may renege before they are served to completion.)
We have
%
\begin{equation}\label{orred}
\mathcal{W}^{(n)}_H\bigl(S^{(n)}_k-\bigr)\bigl(0,F^{(n)} \bigl(S^{(n)}_k\bigr) \bigr) \leq\mathcal
{W}^{(n)}\bigl(S^{(n)}_k-\bigr)\bigl(0,F^{(n)} \bigl(S^{(n)}_k\bigr)
\bigr),
\end{equation}
because, in both systems under consideration, the arrivals with
lead times smaller than $F^{(n)}$ and the
corresponding work associated
with them are the same, the server serves these customers
with rate $1$ as long as they are present in the system, but, by
Theorem \ref{t.EDFopt}, the amount of
unfinished work associated with these customers and
deleted by the
EDF system with reneging is not greater
than the work deleted by the hybrid system.
Thus, if (\ref{highp}) holds, the time required for the hybrid system
to fully serve the $k$th customer is at most
$\mathcal{W}^{(n)}(S^{(n)}_k-)(0,F^{(n)} (S^{(n)}_k) )+v^{(n)}_k$. Therefore,
under assumption (\ref{highp}), the
unfinished work deleted by the hybrid system due to
lateness of the $k$th customer is at most
$v^{(n)}_k \wedge
( \mathcal{W}^{(n)}(S^{(n)}_k-)(0,F^{(n)}(S_k^{(n)})) + v^{(n)}_k
-L^{(n)}_k )^+$.
Thus, the amount of work associated with high-priority customers
deleted by the hybrid system up to time $t$ is bounded above by
the right-hand side of (\ref{e.step3c}).
\end{pf}


\section{Heavy traffic analysis}\label{S5}

In Sections \ref{subs-htwork} and \ref{subs-htren}, respectively,
we identify the heavy traffic limit of the
scaled workload and the scaled reneged work in the reneging system.
In both cases, this is done by first considering the reference
system, which is easier to analyze, and
then using the bounds derived
in Section \ref{subs-comparison} to show that the limits
in both systems coincide. For the heavy traffic analysis of the
reference system, we will find it useful
to introduce the following scaled quantities:
%
\begin{eqnarray}\label{6.0a}
\widehat{U}^{(n)} (t) &\define& \frac{1}{\sqrt{n}}U^{(n)}(nt),\qquad
\widehat{R}_U^{(n)} (t) \define\frac{1}{\sqrt{n}}
\renu(nt),\nonumber\\[-8pt]\\[-8pt]
\widehat{K}_+^{(n)} (t) &\define& \frac{1}{\sqrt{n}} K_+^{(n)}
(nt),\nonumber
\end{eqnarray}
and, for every Borel set $B \subset\R$,
%
\begin{equation}\label{6.0b}
\widehat{\sU}^{(n)} (t) (B)
\define\frac{1}{\sqrt{n}}\sU^{(n)} (nt) \bigl(\sqrt{n} B\bigr).
\end{equation}
Also, define
%
\begin{equation}
\label{def-ustar}
\mathcal{U}^* \define\Phi(\sW^*_S)
\quad\mbox{and}\quad
U^* (\cdot) \define\mathcal{U^*}(\cdot) (\R) = \Phi(\sW^*_S)(\R).
\end{equation}

\subsection{Proofs of main results concerning the workload}
\label{subs-htwork}

\subsubsection{\texorpdfstring{Proof of Theorem \protect\ref{C2.1}}{Proof of Theorem 3.4}}
\label{subsub-pre}

In Lemma \ref{L.drbm},
we use the continuity property of the mapping\vadjust{\goodbreak} $\Phi$ established in
Lemma \ref{lem-Phi}, along with the
characterization of the heavy traffic limit of
the workload measure-valued process in the
standard system,
to identify the heavy traffic limit of the
workload in the
reference system. Let
$\Lambda_{H(0)}\dvtx
D[0,\infty)\rightarrow D[0,\infty)$ be
the mapping
defined, for every $\phi\in D[0,\infty)$ and $t\geq0$, by
%
\begin{equation}
\label{def-lambda}
\Lambda_{H(0)}(\phi)(t)\define\phi(t)- \sup_{s \in[0,t]}
\Bigl[ \bigl( \phi(s)-H(0)\bigr)^
+\wedge\inf_{u \in[s,t]} \phi(u) \Bigr].
\end{equation}
If $\phi$
is nonnegative, then by Theorem 1.4
from \cite{doublereflection},
$\Lambda_{H(0)}(\phi)$ is the function in $D[0,\infty)$
obtained by double reflection of $\phi$ at
$0$ and $H(0)$. In other words,
$\Lambda_{H(0)}(\phi)$ takes values
in $[0,H(0)]$ and has the unique decomposition
%
\begin{equation}\label{decomp}
\Lambda_{H(0)}(\phi)=\phi-\kappa_++\kappa_-,
\end{equation}
where $\kappa_{\pm}$ are nondecreasing RCLL functions
satisfying $\kappa_{\pm}(0-)=0$ and
%
\begin{eqnarray}\label{lc}
\int_{[0,\infty)}\ind_{\{\Lambda_{H(0)}(\phi)(s)<H(0)\}}
\,d\kappa_+(s)&=&0,\nonumber\\[-8pt]\\[-8pt]
\int_{[0,\infty)}\ind_{\{\Lambda_{H(0)}(\phi)(s)>0\}}
\,d\kappa_-(s)&=&0.\nonumber
\end{eqnarray}
\begin{lemma}\label{L.drbm}
The process $U^*$ satisfies
%
\begin{equation}\label{ustardr}
U^* = \Lambda_{H(0)} (W_S^*)
\end{equation}
and has the same distribution as $W^*$.
Moreover,
${\widehat{\sU}^{(n)}} \Rightarrow{\sU}^* = \Phi(\sW^*_S)$ and
$\widehat{U}^{(n)}\Rightarrow W^*$ as $n\rightarrow\infty$.
\end{lemma}
\begin{pf}
By the definition of $U^*$ and
$\Phi$ given in (\ref{def-ustar})
and (\ref{Phi}),
respectively,
\[
U^* (t) = \Phi(\sW_S^*) (\R) (t) =
W_S^* (t) - \sup_{s \in[0,t]} \Bigl[
\mathcal{W}_S^* (-\infty, 0]
\wedge\inf_{u \in[s,t]} W_S^* (u) \Bigr],\qquad
t\geq0.
\]
Since (\ref{5.1})--(\ref{5.3S}) imply
$\sW_S^*(t) (-\infty, 0] = (W_S^* (t)- H(0))^+$
for every $t \geq0$,
this shows that
$U^* = \Lambda_{H(0)} (W_S^*)$.
By the characterization of
$W_S^*$ given at the end of Section \ref{S2.2},
$\Lambda_{H(0)}(W^*_S)$ is a Brownian motion
with variance $(\alpha^2 +\beta^2)\lambda$ per unit time
and drift $-\gamma$, reflected at $0$ and $H(0)$.
This proves the first claim.

Next, using the definition
$\mathcal{U}^{(n)} = \Phi(\sW_S^{(n)})$ and the scaling
properties of $\Phi$, it is easy to see that
$\widehat{\sU}^{(n)} = \Phi(\widehat{\sW}_S^{(n)}).$
Since, by Theorem \ref{T2.2S}, we know
that $\widehat{\sW}_S^{(n)} \Rightarrow\sW_S^*$,
where $\sW_S^*$ is continuous and
$\sW_S^*(t)$ has a continuous distribution for every~$t$,
an application of the continuous mapping 
theorem, along with the continuity property
of $\Phi$ stated in Lemma \ref{lem-Phi}, shows
that
$\widehat{\sU}^{(n)} \Rightarrow\Phi(\sW_S^*)$.
This, in particular, implies that
$\widehat{U}^{(n)} = \widehat{\sU}^{(n)} (\R) \Rightarrow U^*$.
Since $U^*$ has the same distribution as $W^*$,
this proves the lemma.
\end{pf}

We identify the heavy traffic limit of the
workload in the reneging system. We start with
Proposition \ref{P.crush},
which states that the number of customers in the
EDF system with reneging having lead times not greater
than the current frontier and the work associated with
these customers are negligible under heavy traffic scaling.
Then, in Corollary \ref{C.U=W}, we
use the comparison results established in
Section \ref{subs-comparison} to
show that the workloads in the reference and reneging systems are
equal with high probability and so
their heavy traffic limits coincide.
\begin{proposition}\label{P.crush}
The processes $\widehat{\mathcal{W}}^{(n)} (0,\Fhat^{(n)}]$
and $\widehat{\mathcal{Q}}^{(n)}(0,\Fhat^{(n)}]$ converge
in distribution to zero as $n\rightarrow\infty$.
\end{proposition}

This result holds for the same
reason that state-space collapse occurs for
priority queues, an idea that can be traced
back to \cite{whitt71}. Specifically, in our model,
due to the nature of the EDF service discipline,
the entire capacity of
the server is always devoted to work that
lies to the left or at the frontier,
as long as the system is nonempty.
Thus the process $\mathcal{W}^{(n)} (0,F^{(n)}]$
is equal to the
workload in a single-server
$GI/G/1$ queue that has\vspace*{1pt} netput process
$\mathcal{V}^{(n)}(t)(-\infty,F^{(n)}(t)] - t$,
$t \geq0$.
By showing that $F^{(n)} (t) < \sqrt{n}y^*$, one shows that
this (high-priority) queue is in
light traffic as $n \ra\infty$, and so its
diffusion scaling\vspace*{1pt} vanishes in the limit.
Since a rigorous proof that
$\widehat{\mathcal{W}}^{(n)} [\Chat^{(n)},\Fhat^{(n)}]\Rightarrow0$
and
$\widehat{\mathcal{Q}}^{(n)} [\Chat^{(n)},\Fhat^{(n)}]\Rightarrow0$
would be very
similar to the proofs of Proposition 3.6 and\vspace*{1pt}
Corollary 3.8 in \cite{bogdan},
we omit the details.
We note that $\widehat{\mathcal{W}}^{(n)} (0,\Chat^{(n)})
=\widehat{\mathcal{Q}}^{(n)} (0,\Chat^{(n)})= 0$ by
definition.
\begin{corollary}\label{C.U=W}
Let $T>0$. As $n\rightarrow\infty$,
%
\begin{equation}\label{ee.U=W}
\PP\bigl[U^{(n)}(t)=W^{(n)}(t),
0\leq t\leq n T\bigr]\rightarrow1.
\end{equation}
\end{corollary}
\begin{pf}
Because customers with strictly positive
lead times do not renege, we have
$\mathcal{W}^{(n)}(S^{(n)}_k-)(0,F^{(n)}(S_k^{(n)}))
\leq\mathcal{W}^{(n)}(S^{(n)}_k)(0,F^{(n)}(S_k^{(n)}))$
for $k\geq1$.
Thus, by Lemmas \ref{L.step2} and \ref{L5.3},
to prove (\ref{ee.U=W}),
it suffices to show that as $n\rightarrow\infty$,
\[
\PP\bigl[\mathcal{W}^{(n)}\bigl(S^{(n)}_k\bigr)\bigl(0,F^{(n)}\bigl(S_k^{(n)}\bigr)\bigr)
+ v^{(n)}_k \leq L^{(n)}_k,
1\leq k \leq A^{(n)}(nT) \bigr]\rightarrow1.
\]
However, this follows from the fact that, by (\ref{w*}),
\[
\max_{1\leq k \leq A^{(n)}(nT)} v^{(n)}_k
= \sqrt{n} \max_{0\leq t\leq T}
\bigtriangleup\widehat{N}^{(n)}_S(t) = o\bigl(\sqrt{n}\bigr),
\]
the inequalities $L^{(n)}_k\geq\sqrt{n}y_*$, $y_*>0$,
and Proposition \ref{P.crush}.
\end{pf}

Theorem \ref{C2.1} now follows immediately
from Lemma \ref{L.drbm} and Corollary \ref{C.U=W}.

\subsubsection{\texorpdfstring{Proofs of Proposition \protect\ref{P2.1}
and Theorem \protect\ref{T2.2}}{Proofs of Proposition 3.5 and Theorem 3.6}}
\label{S.mainp}

We present the proofs of the remaining
two limit theorems concerning the measure-valued
workload processes.
For this, we need two preliminary results.
The first, Lemma
\ref{L.F>0}, is that the frontier in the reneging
system is strictly positive with high probability.
The second result,
Proposition \ref{P.input}, is a recap of a
result established in \cite{bogdan}.
\begin{lemma}\label{L.F>0}
Let $T>0$. As $n\rightarrow\infty$,
%
\begin{equation}\label{e.U=W}
\PP\bigl[F^{(n)}(t)>0, 0\leq t\leq n T\bigr]\rightarrow1.
\end{equation}
\end{lemma}
\begin{pf}
Let $0\leq t\leq n T$. If $W^{(n)}(t)>0$,
then $F^{(n)}(t)$ is not smaller than the
lead time of the currently served customer,
so $F^{(n)}(t)>0$. If $W^{(n)}(t)=0$, then
the customer indexed by $A^{(n)}(t)$ has already been in service, so
%
\begin{eqnarray}\label{Flbound}
F^{(n)}(t)
&\geq&
L^{(n)}_{A^{(n)}(t)}-\bigl(t-S^{(n)}_{A^{(n)}(t)}\bigr)
\nonumber\\
&\geq&
\sqrt{n}y_*-u^{(n)}_{A^{(n)}(t)+1}\\
&\geq&
\sqrt{n}y_*-\max_{1\leq k\leq A^{(n)}(nT)+1} u^{(n)}_k.
\nonumber
\end{eqnarray}
However, $\max_{1\leq k\leq A^{(n)}(nT)+1} u^{(n)}_k=o(\sqrt{n})$
by (\ref{An}) (in particular, by the fact that $S^*$ has
continuous sample paths), so (\ref{Flbound}) implies (\ref{e.U=W}).
\end{pf}
\begin{proposition}[(Proposition 3.4 \cite{bogdan})]\label{P.input}
Let $-\infty<y_0<y^*$ and $T>0$ be given. As $n\rightarrow\infty$,
\begin{eqnarray*}
\sup_{y_0\leq y\leq y^*} \sup_{0\leq t\leq T}\bigl|
\widehat{\mathcal{V}}^{(n)}(t)(y,\infty)+H\bigl(y+\sqrt{n}t\bigr)-H(y)\bigr|
&\inprob&0,\\
\sup_{y_0\leq y\leq y^*} \sup_{0\leq t\leq T}
\bigl| \widehat{\mathcal{A}}^{(n)}(t)(y,\infty)
+\lambda H\bigl(y+\sqrt{n}t\bigr)- \lambda H(y)\bigr| &\inprob&0.
\end{eqnarray*}
\end{proposition}
\begin{pf*}{Proof of Proposition \ref{P2.1}} Let $T>0$.
We will show that $ \widehat{F}^{(n)} \Rightarrow F^*$
in $D_{\R}[0,T]$.
By definition, $y^*-\sqrt{n}t\leq\Fhat^{(n)}(t)\leq y^*$.
Thus, by Proposition \ref{P.input} and the fact that $H(y)=0$
for $y\geq y^*$,
%
\begin{equation}\label{e.p3.9}
\sup_{0\leq y\leq y^*} \sup_{0\leq t\leq T}\bigl|
\widehat{\mathcal{V}}^{(n)}(t)\bigl(\Fhat^{(n)}(t)\vee y,\infty\bigr)
-H\bigl(\Fhat^{(n)}(t)\vee y\bigr)\bigr|
\inprob0.
\end{equation}
Putting $y=0$ in (\ref{e.p3.9}) and using Lemma \ref{L.F>0},
we obtain
%
\begin{equation}\label{e.m2}
\sup_{0\leq t\leq T}\bigl|
\widehat{\mathcal{V}}^{(n)}(t)\bigl(\Fhat^{(n)}(t),\infty\bigr)
-H\bigl(\Fhat^{(n)}(t)\bigr)\bigr|
\inprob0.
\end{equation}
For any $t\geq0$,
%
\begin{eqnarray}\label{e.bycrush}
\What^{(n)}(t)
&=&
\widehat{\mathcal{W}}^{(n)}(t)\bigl(0,\Fhat^{(n)}(t)\bigr]
+
\widehat{\mathcal{W}}^{(n)}(t)\bigl(\Fhat^{(n)}(t),\infty\bigr)\nonumber\\[-8pt]\\[-8pt]
&=&
\widehat{\mathcal{W}}^{(n)}(t)\bigl(0,\Fhat^{(n)}(t)\bigr]
+ \widehat{\mathcal{V}}^{(n)}(t)\bigl(\Fhat^{(n)}(t),\infty\bigr),
\nonumber
\end{eqnarray}
where the second line follows from the fact that none of the customers
in the EDF system with reneging with lead times at time $t$ greater
than $F^{(n)}(t)$ has received any service up\vspace*{1pt} to time $t$. This,
together with Proposition \ref{P.crush} and Theorem~\ref{C2.1}, yields
$\widehat{\mathcal{V}}^{(n)}(\Fhat^{(n)},\infty)\Rightarrow W^*$. Thus,
by (\ref{e.m2}), we have $H(\Fhat^{(n)})\Rightarrow W^*$ in
$D_{\R}[0,T]$. Applying the continuous function $H^{-1}$ to both sides
of this relation and using (\ref{5.2}), we obtain $ \widehat{F}^{(n)}
\Rightarrow F^*$ in $D_{\R}[0,T]$.
\end{pf*}
\begin{pf*}{Proof of Theorem \ref{T2.2}}
Define a mapping $\psi\dvtx\R\rightarrow\sMM$ by the formula
$\psi(x)(B)\define \int_{B\cap[x,\infty)}(1-G(\eta)) \,d\eta$ for
$x\in\R$ and $ B\in\sB(\R)$. It is easy\vspace*{1pt} to see that $\psi$ is
continuous. Hence, by Proposition \ref{P2.1},
%
\begin{equation}\label{e.fin1}
\bigl(\psi\bigl(\Fhat^{(n)}\bigr),\lambda\psi\bigl(\Fhat^{(n)}\bigr)\bigr)
\Rightarrow(\psi(F^*),\lambda\psi(F^*))
=(\mathcal{W}^*,\mathcal{Q}^*).
\end{equation}
Let $T>0$. We claim that
%
\begin{eqnarray}
\label{e.fin2}
\sup_{0\leq y \leq y^*} \sup_{0\leq t\leq T}\bigl|
\widehat{\mathcal{W}}^{(n)}(t)(y,\infty)
-\psi\bigl(\Fhat^{(n)}(t)\bigr)(y,\infty)\bigr|
&\inprob&0,\\
\label{e.fin2a}
\sup_{0\leq y \leq y^*} \sup_{0\leq t\leq T}\bigl|
\widehat{\mathcal{Q}}^{(n)}(t)(y,\infty)
-\lambda\psi\bigl(\Fhat^{(n)}(t)\bigr)(y,\infty)\bigr|
&\inprob&0.
\end{eqnarray}
Indeed, reasoning as in (\ref{e.bycrush}), we see that,
for $0\leq y \leq y^*$ and $0\leq t\leq T$,
\begin{eqnarray*}
&&\bigl| \widehat{\mathcal{W}}^{(n)}(t)(y,\infty)
-H\bigl(\Fhat^{(n)}(t)\vee y\bigr)\bigr|\\
&&\qquad\leq
\bigl| \widehat{\mathcal{W}}^{(n)}(t)
\bigl(\Fhat^{(n)}(t)\vee y,\infty\bigr)-H\bigl(\Fhat^{(n)}(t)\vee y\bigr)\bigr|
+ \widehat{\mathcal{W}}^{(n)}(t)\bigl(0,\Fhat^{(n)}(t)\bigr]\\
&&\qquad=
\bigl| \widehat{\mathcal{V}}^{(n)}(t)\bigl(\Fhat^{(n)}(t)
\vee y,\infty\bigr)-\psi\bigl(F^{(n)}(t)\bigr)(y,\infty)\bigr|
+ \widehat{\mathcal{W}}^{(n)}(t)\bigl(0,\Fhat^{(n)}(t)\bigr].
\end{eqnarray*}
Therefore, (\ref{e.fin2}) follows from (\ref{e.p3.9})
and Proposition \ref{P.crush}. A similar argument gives
(\ref{e.fin2a}).
We have\vspace*{1pt}
$\widehat{\mathcal{W}}^{(n)}(t)(-\infty,0]
= \widehat{\mathcal{W}}^{(n)}(t)(y^*,\infty)
= \widehat{\mathcal{Q}}^{(n)}(t)(-\infty,0]
= \widehat{\mathcal{Q}}^{(n)}(t)(y^*,\infty)
= 0$
and, by Lemma \ref{L.F>0},
$\PP[\psi(\Fhat^{(n)}(t))(-\infty,0]
=0, 0\leq t\leq T]\rightarrow1$
as $n\rightarrow\infty$. Also, $\psi(x)(y^*,\infty)=0$
for every $x\in\R$.
Thus, (\ref{e.fin1})--(\ref{e.fin2a}) imply that
$(\widehat{\mathcal{W}}^{(n)},\widehat{\mathcal{Q}}^{(n)})
\Rightarrow(\mathcal{W}^*,\mathcal{Q}^*)$
in $D_{\sMM}[0,T]$.
\end{pf*}

\subsection{The heavy traffic limit of the reneged work process}
\label{subs-htren}

In this section, we identify the limit of the
sequence $\{\hrenw, n \in\N\}$, thereby
proving Theorem \ref{T3.1}.
To do this, it is convenient to show
that many of the processes
under consideration can be put
on a common probability space so that
certain weak limits established earlier can be replaced by almost sure
limits.
\begin{lemma}
\label{lem-probspace}
The processes $\sWhat_S^{(n)}$, $\sUhat^{(n)}$, $\What^{(n)}$, $n
\in\N$,
$\sW_S^*$, $\sU^*$ and $W^*$
can be defined on a common probability space $(\Omega, \mathcal{F},
\PP)$ such
that $\PP$ almost surely, as $n \ra\infty$,
%
\begin{eqnarray}
\label{swhatas}
\sWhat^{(n)}_S &\rightarrow& \sW^*_S,
\\
%
\label{whatas}
\What_S^{(n)} &\ra& W_S^*,
\\
%
\label{swhatlate}
\sWhat^{(n)}_S(\cdot)(-\infty,0] & \rightarrow &
\sW_S^*(\cdot)(-\infty,0]
=\bigl(W^*_S(\cdot)-H(0)\bigr)^+,
\\
%
\label{uhat}
\Uhat^{(n)} &\rightarrow& U^*
\end{eqnarray}
and
%
\begin{equation}\label{what}
\What^{(n)}\rightarrow W^* \define U^*,
\end{equation}
where $\What^{(n)}_S=\sWhat^{(n)}_S(\R)$,
$W^*_S=\mathcal{W}^*_S(\R)$, $\widehat{U}^{(n)}=\widehat{\mathcal
{U}}^{(n)}(\R)$
and $U^*=\mathcal{U}^*(\R)$. Furthermore, $W^*_S$ is a Brownian
motion with variance
$(\alpha^2+\beta^2)\lambda$ per unit time and drift $-\gamma$,
reflected at $0$, while $U^*$ is a doubly reflected Brownian
motion on $[0,H(0)]$, also with variance $(\alpha^2+\beta^2)\lambda$
per unit time and drift $-\gamma$. In particular,
%
\begin{equation}\label{ustardcomp}
U^*= \Lambda_{H(0)} (W_S^*) = W^*_S-K_+^*+K_-^*,
\end{equation}
where $K_{\pm}^*$ are the unique RCLL nondecreasing functions satisfying
$K_{\pm}^*(0)=0$ and
%
\begin{equation}\label{ustarlc}
\int_{[0,\infty)}\ind_{\{U^*(s)<H(0)\}} \,dK^*_+(s)=0,
\int_{[0,\infty)}\ind_{\{U^*(s)>0\}} \,dK^*_-(s)=0.
\end{equation}
The almost sure limits in (\ref{swhatas})--(\ref{what})
hold uniformly on compact intervals.
\end{lemma}
\begin{pf}
Recall from Theorem \ref{T2.2S} that
$\sWhat^{(n)}_S\Rightarrow\sW^*_S$. Using the Skorokhod representation
theorem, we construct the model primitives\vspace*{2pt} $u_j^{(n)}$,
$v_j^{(n)}$ and $L_j^{(n)}$ for $j\in\N$ and $n\in\N$ on a common
probability space $(\Omega,\mathcal{F},\PP)$ such that the sequence of
processes $\sWhat^{(n)}_S$, $n\in\N$, and the limiting process
$\sW^*_S$ are defined on this space and (\ref{swhatas}) holds. Here and
below the almost sure convergences are in the $J_1$ topology on
$D_\mathcal{M}[0,\infty)$ or $D_{\R}[0,\infty)$, and since the limits
are continuous in every case, this is equivalent to uniform convergence
on compact intervals. Since the mapping $f\dvtx D_\mathcal{M}[0,\infty)
\mapsto D_{\R}[0,\infty)$ given by $f(\mu)(\cdot)= \mu(\cdot)(\R)$ is
continuous, we have~(\ref{whatas}). Under $\PP$ the measure-valued
process $\sW^*_S$ constructed on $\Omega$ has the same distribution as
the process $\sW^*_S$ appearing in Theorem \ref{T2.2S}, and thus
$\sW^*_S$ takes values in the set of measure-valued process of the form
$\int_{B\cap[F^{o}_S(t),\infty)}(1-G(y))\,dy$ for some RCLL process
$F^{o}_S(t)$. However,
$W^*_S(t)=\int_{\R\cap[F^o_S(t),\infty)}(1-G(y))\,du =H(F^o_S(t))$; hence
$F^o_S(t)=F^*_S(t)$ is given by (\ref{5.2S}). In other words, with
$F^*_S$ defined by (\ref{5.2S}), the first equation in (\ref{5.3S})
holds. Due to Proposition \ref{P2.1S}, the\vspace*{1pt} above argument also shows
that under $\PP$, $W_S^*$ is a Brownian motion with variance $(\alpha^2 +\beta^2)
\lambda$ per unit time and drift $-\gamma$. In addition, since for each
$t$, the measure $\sW^*_S(t)$ is nonatomic, we have (\ref{swhatlate}).

Now, following (\ref{refworm}) and (\ref{def-ustar}), we set $\sU^{(n)}
= \Phi(\sW_S^{(n)})$ and $\sU^*=\Phi(\sW_S^*)$. Also, as defined in
(\ref{6.0b}), let $\sUhat^{(n)}$ be the scaled version of $\sU^{(n)}$,
and let $\Uhat^{(n)}$ and $\Uhat^*$ be as defined in the statement of
the lemma. Then $\sUhat^{(n)}, \Uhat^{(n)}, n \in\N$, $\sU^*$ and $U^*$
are also defined on $(\Omega, \mathcal{F}, \PP)$ and (\ref{uhat})
follows from Lemma \ref{lem-Phi}. This implies (\ref{what}). Since $U^*
= \Phi(\sW^*_S)(\R)=\Lambda_{H(0)}(W^*_S)$, the characterization of
$U^*$ as a doubly reflected Brownian motion that satisfies relations
(\ref{ustardcomp}) and (\ref{ustarlc}) is a consequence of the
statements following (\ref{def-lambda}), in particular, (\ref{decomp})
and (\ref{lc}).

Since the model primitives $u_j^{(n)}$, $v_j^{(n)}$ and $L_j^{(n)}$ for
$j\in\N$ and $n\in\N$ are all defined on $(\Omega,\mathcal{F},\PP)$, so
are the workload process $W^{(n)}$ and its scaled version
$\What^{(n)}$. Corollary \ref{C.U=W} implies that $\Uhat^{(n)}$ and
$\What^{(n)}$ have the same limit, and hence (\ref{what}), the almost
sure counterpart to Theorem \ref{C2.1}, holds.
\end{pf}

The assertion of Theorem \ref{T3.1} is that
%
\begin{equation}\label{T3.1claim}
\widehat{R}_W^{(n)}\Rightarrow K_+^*,
\end{equation}
where $K_+^*$ is the local time for $U^*$ at $H(0)$ from
(\ref{ustardcomp}). For $T < \infty$, define
%
\begin{equation}
\label{def-zn}
\mathcal{Z}_n (T) \define\bigl\{ \hrenu(t) = \hrenw(t),
0 \leq t \leq T \bigr\}.
\end{equation}
From the workload evolution equations (\ref{WAccount}) and
(\ref{UAccount}), it follows that if $\widehat{U}^{(n)} (t) =
\widehat{W}^{(n)} (t)$ for $t \in[0,T]$, then $\hrenu(t) = \hrenw(t)$
for $t\in[0,T]$. Hence, by Corollary \ref{C.U=W}, we know that for
every $T < \infty$, $\PP(\mathcal{Z}_n(T)) \ra1$ as $n \ra\infty$,
which shows that the limits in distribution of $\hrenu$ and $\hrenw$,
$n\in\N$, must coincide (if they exist). Further, since $\hkpn=\hrenu$
by Corollary \ref{C.eqn4.9}, these must be equal to the limit in
distribution of $\hkpn, n \in\N$. Hence, to complete the proof of
Theorem \ref{T3.1}, it suffices to show that
%
\begin{equation}\label{T3.6p}
\hkpn\Rightarrow K_+^*.
\end{equation}

For $n \in\N$ and $k \geq1$, recall the definitions of $\tau
_{k-1}^{(n)}$ and
$\sigma_k^{(n)}$ given in (\ref{2}) and (\ref{3}), respectively, and
define
$\htaun_{k-1} \define\frac{1}{n}
\tau_{k-1}^{(n)}$ and
$\hsigman_k \define
\frac{1}{n} \sigma_k^{(n)}$.
Applying the heavy traffic scaling to (\ref{def-kpn}),
it is easy to see that for $t \geq0$,
%
\begin{eqnarray}
\label{def-hkpn}
\hkpn(t) &=& \sum_{k \in\N} \Bigl[
\widehat{W}_S^{(n)} \bigl(\hsigman_k -\bigr) \vee
\max_{s \in[\hsigman_k,t\wedge\htaun_k]}
\widehat{\sW}_S^{(n)}(s) (-\infty,0]\nonumber\\[-8pt]\\[-8pt]
&&\hspace*{151.3pt}{}
- \widehat{W}_S^{(n)} \bigl(\hsigman_k -\bigr) \Bigr] .\nonumber
\end{eqnarray}
Keeping in mind the limits in (\ref{swhatas})
and (\ref{swhatlate}), we introduce the related process
%
\begin{equation}
\label{def-hyn}
\hyn(t) \define\sum_{k \in\N} \Bigl[
W_S^* \bigl(\hsigman_k\bigr) \vee
\max_{s \in[\hsigman_k,t\wedge\htaun_k]}
\bigl( W^*_S (s) - H(0)\bigr)^+
- W_S^* \bigl(\hsigman_k \bigr) \Bigr]\hspace*{-32pt}
\end{equation}
for $t \geq0$, and denote the difference by
%
\begin{equation}
\label{def-ven}
\ven(t) \define\hyn(t) - \hkpn(t)\qquad
\forall t \geq0.
\end{equation}
Then
$\hyn$ is nondecreasing and continuous, and
$\ven$ is an RCLL process.

In the next two lemmas, we show that $\hyn$ increases
only when $U^*$ is at $H(0)$
and that the difference $\ven$ between
$\hyn$ and $\hkpn$ is negligible in heavy traffic.
The main reason for introducing the sequence $\hyn$,
$n \in\N$, is that
it facilitates the proof of the former property.
\begin{lemma}
\label{lem-hyn} For every $n \in\N$,
$\hyn$ and $\hkpn$ are constant on each interval
$\idle$, $k \geq1$. Moreover,
%
\begin{equation}
\label{prop-hyn}
\int_{[0,T]} \ind_{\{U^{*} (t) < H(0)\}} \,d \hyn(t) = 0.
\end{equation}
\end{lemma}
\begin{pf}
Fix $n \in\N$.
The first statement follows immediately from
(\ref{def-hkpn}), (\ref{def-hyn}), and
the fact that the intervals $[\htaun_{k-1},\hsigman_k)$
and $[\hsigman_k,\htaun_k)$,
$k \geq1$, form a disjoint covering of $[0,\infty)$.
Now, fix $k \geq1$ and let $J_k^{(n)}$ be the set of
points $t \in[\hsigman_k,\htaun_k)$ such that
%
\begin{equation}
\label{cond}
W_S^* \bigl(\hsigman_k \bigr)
\leq\max_{s \in[\hsigman_k,t]} \bigl(W^*_S (s) - H(0)\bigr)^+
= W^*_S (t) - H(0).
\end{equation}
Since $W_S^*$ is continuous, $J_k^{(n)}$ is closed,
and so its complement in $[\sigma_k^{(n)},\tau_k^{(n)})$ is
the union of a countable number of open intervals,
with possibly one half-open interval of the form
$[\hsigman_k, a)$ for some
$a > \hsigman_k$. From the explicit formula for
$\hyn$ given in (\ref{def-hyn}), it is easy to deduce that
$\hyn$ is also constant on each such interval.
Thus, to establish (\ref{prop-hyn}), it only remains
to show that for each $k \geq1$,
%
\begin{equation}
\label{jkn}
\int_{J_k^{(n)}} \ind_{\{U^*(t) < H(0)\}} \,d\hyn(t) = 0.
\end{equation}

Fix $t \in J_k^{(n)}$ and note that by the equality in
(\ref{ustardcomp}) and the definition (\ref{def-lambda})
of $\Lambda_{H(0)}$, we have $U^* (t) = W_S^*(t) - K^*(t)$, where
%
\begin{equation}
\label{eq-kstar}
K^* (t) \define\sup_{s \in[0,t]}
\Bigl[ \bigl(W^*_S (s)-H(0)\bigr)^+
\wedge\inf_{u \in[s,t]} W_S^*(u)\Bigr].
\end{equation}
Also, note that
\begin{eqnarray*}
\sup_{s \in[0,\hsigman_k)}
\Bigl[ \bigl(W^*_S (s)-H(0)\bigr)^+
\wedge\inf_{u \in[s,t]} W_S^*(u)\Bigr]
&\leq&
\sup_{s \in[0,\hsigman_k)}
\inf_{u \in[s,t]} W_S^*(u)\\
&\leq&
W_S^* \bigl(\hsigman_k\bigr),
\end{eqnarray*}
and that the equality in (\ref{cond}) implies
\[
\sup_{s \in[\hsigman_k,t]}\Bigl[ \bigl(W^*_S (s)-H(0)\bigr)^+
\wedge\inf_{u \in[s,t]} W_S^*(u)\Bigr]
= W^*_S(t)- H(0).
\]
Since $K^*(t)$ is equal to the maximum of the
quantities on the left-hand side of the last two displays,
we conclude that
\[
K^* (t) \leq W_S^* \bigl(\hsigman_k\bigr)
\vee\bigl( W^*_S(t)- H(0)\bigr) = W^*_S (t) - H(0),
\]
where the equality follows from the inequality in (\ref{cond}).
This, when combined with the fact that $U^*(t) \in[0,H(0)]$,
shows that $U^*(t) = W_S^* (t) - K^*(t) = H (0)$ for all
$t \in J_k^{(n)}$, which proves (\ref{jkn}).
\end{pf}

We recall some standard definitions that will be
used in the next lemma.
Given $f \in\mathcal{D}[0,\infty)$ and
$0 \leq t_1 \leq t_2 < \infty$, the
\textit{oscillation of
$f$ over $[t_1,t_2]$}
is
\[
\osc(f; [t_1,t_2]) \define
\sup\{ |f(t)-f(s)|\dvtx t_1 \leq s \leq t \leq t_2 \},
\]
and the \textit{modulus of continuity of $f$ over $[0,T]$} is
\[
w_f(\delta;[0,T]) \define
\sup\{ |f(t) - f(s)|\dvtx 0 \leq s \leq t \leq T, |t-s|\leq\delta\}.
\]
%
%
%
\begin{lemma}
\label{lem-ven}
As $n \ra\infty$, $\ve^{(n)} \inprob0$.
\end{lemma}
\begin{pf}
Fix $T>0$ and
let $\eta>0$ be arbitrarily small. By the
Kolmogorov--\v{C}entsov theorem (see, e.g., Theorem 2.8, page 53 of
\cite{ks}), we can construct a positive,
increasing deterministic function
$\theta(\cdot)$ satisfying
$\lim_{\delta\downarrow0}\theta(\delta)=0$ and
majorizing the modulus of continuity $w_{W^*}(\cdot;[0,T])$
of the reflected Brownian motion $W^*$
over $[0,T]$ on a set
$\widetilde{\Omega}$ with
$\PP( \widetilde{\Omega}) \geq1-\eta$.

For each subsequence in $\N$, there is a sub-subsequence
$\sS$ along which the limits (\ref{swhatas})--(\ref{what})
hold $\PP$-almost surely.
We choose $\widetilde{\Omega}$ so that these
limits hold along $\sS$ for all $\omega\in\widetilde{\Omega}$.

In what follows, for $n\in\sS$, we denote
$\mathcal{Z}_n(T)$ simply by
$\mathcal{Z}_n$, and evaluate all processes below
at a fixed
$\omega\in\mathcal{Z}_n \cap\widetilde{\Omega}$.
Choose $\sma< y_*/3$,
and let $n_0\in\sS$ be such that for all $n\in\sS$,
$n\geq n_0$,
%
\begin{eqnarray}
\label{wlim1}
\sup_{t \in[0,T]} \bigl| \widehat{\sW}_S^{(n)}(t) (-\infty,0]
- \bigl(W_S^*(t) - H(0)\bigr)^+ \bigr| &\leq& \sma,
\\
%
\label{wlim2}
\sup_{t \in[0,T]} \bigl| \widehat{W}^{(n)}_S (t-)
- W^*_S (t) \bigr| &\leq& \sma,\nonumber\\[-8pt]\\[-8pt]
\sup_{t \in[0,T]} \bigl| \widehat{W}^{(n)} (t-)
- W^* (t) \bigr| &\leq& \sma,\nonumber
\\
%
\label{wlim3}
\sup_{t\in[0,T]}\bigl|\widehat{U}^{(n)}(t-)-U^*(t)\bigr|
&\leq& \sma.
\end{eqnarray}
From the definitions (\ref{def-hkpn}) and
(\ref{def-hyn}), respectively, of $\hkpn$ and $\hyn$ it is clear that,
for every $k \in\N$ such that $\tau_k^{(n)}\leq T$,
\[
\sup_{t\in[\hsigman_k,\htaun_k]}
\bigl| \hyn(t) - \hyn\bigl(\hsigman_k -\bigr)
- \bigl( \hkpn(t) - \hkpn\bigl(\hsigman_k-\bigr) \bigr)\bigr|
\leq2\sma.
\]
Define
\begin{eqnarray*}
J_n
& \define&
\bigl\{ k \in\N\dvtx \hkpn\bigl(\htaun_k\bigr) - \hkpn\bigl(\hsigman_k-\bigr) > 0,
\htaun_k \leq T \bigr\}, \\
\widetilde{J}_n
& \define&
\bigl\{ k \in\N\dvtx \hyn\bigl(\htaun_k\bigr) - \hyn\bigl(\hsigman_k-\bigr) > 0,
\htaun_k \leq T \bigr\},
\end{eqnarray*}
and let $c^{(n)}$ be the cardinality of
$J^{(n)} \cup\widetilde{J}^{(n)}$.
Since $\hkpn$ and $\hyn$ are both constant on
intervals of the form $\idle$, $k \geq1$
(see Lemma \ref{lem-hyn}), we have
%
\begin{equation}
\label{ineq-kpnyn}
\overline{\ve}{}^{(n)} (T)
\define
\sup_{s \in[0,T]} \bigl|\hyn(s)-\hkpn(s) \bigr|
\leq2c^{(n)} \sma.
\end{equation}
We now claim that
%
\begin{equation}
\label{ineq-osc}
k \in[J_n \cup\widetilde{J}_n]
\quad\Rightarrow\quad\osc\bigl(W^*, \bigl[\hsigman_k,\htaun_k\bigr]\bigr)
\geq\frac{y_*}{3}.
\end{equation}
We defer the proof of the claim and instead
first show that the lemma follows from this
claim.
Let $\theta^{-1}(\cdot)$ denote the inverse of
$\theta$ and define $M \define T/\theta^{-1} (y_*/3) < \infty$.
From the claim, we conclude that
if $k\in[J_n \cup\widetilde{J}_n]$ then
$\htaun_k - \hsigman_k \geq\theta^{-1} (y_*/3) > 0$, which
in turn implies that $c^{(n)} \leq M$. Substituting this
into (\ref{ineq-kpnyn}), we conclude that for
every $\sma> 0$, there exists
$n_0(\sma)\in\sS$ such that
for all $n\in\sS$, $n \geq n_0(\sma)$,
\[
\PP\bigl(\overline{\ve}^{(n)} (T) > 2M \sma\bigr)
\leq\PP( \mathcal{Z}_n^c \cup\widetilde{\Omega}^c)
\leq\PP(\mathcal{Z}_n^c)+\eta.
\]
Taking limits as $n\ra\infty$ through
$\sS$ and using the fact
that $\PP(\mathcal{Z}_n) \ra1$, we conclude that
$\overline{\ve}^{(n)}(T) \inprob0$.
We have shown that for each subsequence in $\N$,
there is a sub-subsequence along which
$\overline{\ve}^{(n)}(T) \inprob0$.
It follow that $\overline{\ve}^{(n)}(T) \inprob0$,
where the limit is taken over all $n\in\N$,
and this proves the lemma.

We now turn to the proof of the claim (\ref{ineq-osc}).
Note first that by the definition of $H(0)$ and
$y_*$, we have $H(0)\geq y_*$.
If $k \in\widetilde{J}_n$, then Lemma \ref{lem-hyn} shows that
$U^* (t) =H(0)$ for some $t \in\busy$.
By the equality $\widehat{U}^{(n)} (\hsigman_k-) = 0$ proved in
Lem\-ma~\ref{L.Un} and (\ref{wlim3}),
this implies that the oscillation of
$U^*$ on $\busy$ is no less than $H(0) - \sma\geq y_*/3$.
Since $W^*=U^*$, the conclusion in (\ref{ineq-osc}) holds.

Finally, suppose $k\in J_n$. Since $\hkpn= \hrenu=\hrenw$,
we have
\[
\hrenw\bigl(\htaun_k\bigr) - \hrenw\bigl(\hsigman_k-\bigr)>0,
\]
that is, the deadline of a customer in the reneging system
expires during the unscaled time interval
$[\sigma_k^{(n)},\tau_k^{(n)}]$.
Since
%
\begin{equation}\label{6.43}
W^{(n)}\bigl(\sigma_k^{(n)}-\bigr)=0
\end{equation}
[because
$U^{(n)} (\sigma^{(n)}_k-) =0$ and, by Lemma \ref{L.step2},
$W^{(n)} \leq U^{(n)}$], this customer must arrive
during the interval $[\sigma_k^{(n)},\tau_k^{(n)})$.
Since his initial lead time is greater than or equal to
$\sqrt{n} y_*$, there is a time
$nt_0\in[\sigma_k^{(n)},\tau_k^{(n)})$
when this customer has lead time
exactly $\sqrt{n} y_*$.
After time $nt_0$, this customer cannot be preempted
by new arrivals, all of which have initial lead times
greater than or equal to $\sqrt{n}y_*$.
At time $nt_0$, the work that must be completed before
this customer is served to completion is
at most
$\sW^{(n)}(nt_0)(0,\sqrt{n}y_*]$.
Since this customer
becomes late, we must have
$W(nt_0)\geq\sW^{(n)}(nt_0)(0,\sqrt{n}y_*]> \sqrt{n}y_*$,
or equivalently,
$\What^{(n)}(t_0)\geq\sWhat^{(n)}(t_0)(0,y_*]>y_*$.
By right continuity,
$\What^{(n)}((t_0+\nu)-)>y_*$
for some $\nu>0$ so small that
$t_0+\nu\leq\htaun_k$.
From the second inequality in (\ref{wlim2}) and
the fact that $\What^{(n)}(\hsigman_k-)=0$
[the scaled version of (\ref{6.43})],
we conclude that
\[
W^*(t_0+\nu)-W^*\bigl(\hsigman_k\bigr)\geq\frac{y_*}{3},
\]
and this gives us the conclusion in (\ref{ineq-osc}).
\end{pf}
\begin{pf*}{Proof of Theorem \ref{T3.1}}
Fix $T < \infty$.
Let $\newsmall^{(n)} \define U^* - \widehat{U}^{(n)}$, and let
$\overline{\newsmall}{}^{(n)} \define
{\sup_{s\in[0,T]}} |U^*(s) - \widehat{U}^{(n)}(s)|$.
According to (\ref{U}) and (\ref{eq-Kdcomp}),
\[
U^{(n)}=W_S^{(n)}-K_+^{(n)}+K_-^{(n)}.
\]
We scale this equation to obtain
%
\begin{equation}\label{ustar}
U^*=\What_S^{(n)}+\delta^{(n)}-\Khat_+^{(n)}+\Khat_-^{(n)}
=\What_S^{(n)}+\delta^{(n)}+\varepsilon^{(n)}
-\hyn+\Khat_-^{(n)},\hspace*{-28pt}
\end{equation}
where [cf. (\ref{def-knn})]
\[
\Khat_-^{(n)}(t)
\define
-\sum_{k\in\N}\bigl[\bigl(\What_S^{(n)}\bigl(\htaun_{k-1}\bigr)
-\bigl(\hsigman_k\wedge t-\htaun_{k-1}\bigr)\bigr)^+
-\What_S^{(n)}\bigl(\htaun_{k-1}\bigr)\bigr],
\]
$\hkpn$ is defined by (\ref{def-hkpn})
and $\ve^{(n)}$ is defined by (\ref{def-ven}).
According to (\ref{complementarity}),
$\int_0^T\ind_{\{\Uhat{}^{(n)}(t)>0\}} \,d\Khat_-^{(n)}(t)=0$,
which implies
%
\begin{equation}\label{-comp}
\int_0^T\ind_{\{U^*(t)>\overline{\newsmall}^{(n)}\}}
d\Khat_-^{(n)}(t)=0.
\end{equation}

Since $\What_S^{(n)}+\delta^{(n)}+\epsilon^{(n)}\Rightarrow W_S^*$
due to (\ref{whatas}), (\ref{uhat}) and Lemma \ref{lem-ven},
and, by (\ref{ustardcomp}), $U^*$ is obtained by applying the
Skorokhod map on $[0,H(0)]$ to $W^*_S$, the convergence
(\ref{T3.6p}) is an immediate consequence of
(\ref{ustar}), (\ref{-comp}), Lemmas \ref{lem-hyn},
\ref{lem-ven} and the invariance
principle for reflected Brownian motions. However, since we are in
a particularly simple setting here, we will provide a direct
proof without invoking the general invariance principle.

We choose $n_0$ so that $\overline{\newsmall}{}^{(n_0)}<H(0)/3$
and recursively define stopping times $\rho_0=0$,
and for $k\geq1$,
\[
\nu_k=
\min\biggl\{t\geq\rho_{k-1}\Big|
U^*(t)=\frac{2H(0)}{3}\biggr\},\qquad
\rho_k=
\min\biggl\{t\geq\nu_{k}\Big|U^*(t)
=\frac{H(0)}{3}\biggr\}.
\]
Then $0=\rho_0<\nu_1<\rho_1<\nu_2<\cdots$ and
$\lim_{k\rightarrow\infty}\rho_k=\lim_{k\rightarrow\infty}\nu_k
=\infty$.
For $n\geq n_0$, $\Khat_-^{(n)}$ is constant on each
of the intervals $[\nu_k,\rho_k]$. Moreover,
Lemma \ref{lem-hyn} implies that for each $k$,
$\hyn$ is constant on each of the intervals $[\rho_{k-1},\nu_k]$.
For $t\in[\nu_k,\rho_k]$, we have from (\ref{ustar}),
(\ref{whatas}), (\ref{uhat}) and Lemma \ref{lem-ven} that
\begin{eqnarray*}
\hyn(t)-\hyn(\nu_k)
&=&
\What_S^{(n)}(t)-U^*(t)+\delta^{(n)}(t)+\varepsilon^{(n)}(t)
\\
&&{}
-\What_S^{(n)}(\nu_k)+U^*(\nu_k)-\delta^{(n)}(\nu_k)
-\varepsilon^{(n)}(\nu_k)\\
&\inprob&
W^*_S(t)-U^*(t)-\bigl(W^*_S(\nu_k)-U^*(\nu_k)\bigr).
\end{eqnarray*}
It follows that, uniformly for $t\in[0,T]$,
$\hyn(t)$ converges in probability to
%
\begin{equation}\label{YtoK}
\sum_{k\in\N}
\bigl[W_S^*\bigl((t\vee\nu_k)\wedge\rho_k\bigr)
-U^*\bigl((t\vee\nu_k)\wedge\rho_k\bigr)-
\bigl(W_S^*(\nu_k)-U^*(\nu_k)\bigr)\bigr].\hspace*{-28pt}
\end{equation}
However, (\ref{ustarlc}) implies that for each $k$,
$K_-^*$ is constant on $[\nu_k,\rho_k]$, and $K_+^*$
is constant on $[\rho_{k-1},\nu_k]$.
Therefore, (\ref{ustardcomp}) implies
that for $t\in[\nu_k,\rho_k]$,
\[
K_+^*(t)-K_+^*(\nu_k)
=W_S^*(t)-U^*(t)-\bigl(W_S^*(\nu_k)-U^*(\nu_k)\bigr).
\]
This implies that the expression in (\ref{YtoK})
is $K_+^*(t)$. But $\hyn$ and $\Khat_+^{(n)}$
have the same limit in probability because of Lemma
\ref{lem-ven}, and we conclude that
%
\begin{equation}\label{endp}
\max_{t\in[0,T]}\bigl|\Khat_+^{(n)}(t)- K_+^*(t)\bigr|
\inprob0.
\end{equation}
Convergence in probability
implies weak convergence, and we have
(\ref{T3.6p}).
\end{pf*}

\section{Performance evaluation and simulations}\label{S.sim}

We use the heavy traffic approximations of this paper to evaluate the
performance of the system with reneging and compare this to the system
in which all customers are served to completion. The predictions of the
theory, derived in Section \ref{Subsection7.3} and compared to
simulations in Section \ref{Subsection7.2}, are predicated on the
assumption that one can interchange the limit as $n\rightarrow\infty$
and the limit as time goes to infinity of the fraction of reneged work.
A formal proof would require a coupling argument such as that found in
\cite{ZZ}. The simulation results attest to the accuracy of the
approximations derived in Section~\ref{Subsection7.3} and also show the
great difference in performance between the reneging and nonreneging
systems.

\subsection{Derivation of theory
predictions}\label{Subsection7.3}

We derive formulas (\ref{frw})--(\ref{flcr}). We begin with one of the
main results of this paper, Theorem \ref{C2.1}, which states that the
limiting scaled workload in the reneging system is a reflected Brownian
motion in $[0,H(0)]$ with drift. More specifically,
%
\begin{equation}\label{push}
W^*(t)=W_S^*(t)-K_+^*(t)+K_-^*(t),
\end{equation}
where $W_S^*(t)$ is a reflected Brownian motion on $[0,\infty)$ with
variance $\sigma^2=\lambda(\alpha^2+\beta^2)$ per unit time and drift
$-\gamma$, $K_-^*$ is the nondecreasing process starting at
$K_-^*(0)=0$ that grows only when $W^*=0$, and $K_+^*$ is the
nondecreasing process starting at $K_+^*(0)=0$ that grows only when
$W^*=H(0)$. We further saw in Theorem \ref{T3.1} that $K_+^*(t)$ is the
limit of the scaled workload that reneges prior to time $t$ in the
diffusion scaling, that is, $\sqrt{n}K_+^*(t)$ is approximately the
(unscaled) workload that reneges in the $n$th system prior to time
$nt$.
\begin{lemma}[(\cite{harrison}, Proposition 5, page 90)]
\label{L7.1}
We have
%
\begin{equation}\label{mstar}
\lim_{t\ra\infty}\frac{1}{t} K_+^*(t)=
\lim_{t\rightarrow\infty}\frac{1}{t}\E K_+^*(t)
=\cases{
\dfrac{\gamma}{e^{2\gamma H(0)/\sigma^2}-1}, &\quad if $\gamma\neq
0$,\vspace*{2pt}\cr
\dfrac{\sigma^2}{2 H(0)}, &\quad if $\gamma=0$.}
\end{equation}
\end{lemma}
\begin{pf}
The first equality in (\ref{mstar}) is a consequence of the fact that
$W^*$ has a stationary distribution [see (\ref{7.5a}) below]. For the
proof of the second equality, recall that $W_S^*$ has the decomposition
(\ref{2.16c}). Let $f$ be a $C^2$ function. Applying It\^o's formula to
$f(W^*(t))$ and taking expectations, we obtain
%
\begin{eqnarray}\label{feq}
&&f'(0)\E[I_S^*(t)+K_-^*(t)]
-f'(H(0))\E K_+^*(t)
\nonumber\\
&&\qquad=
\E\int_0^t\biggl[\gamma f'(W^*(s))
-\frac12\sigma^2f''(W^*(s))\biggr]\,ds\\
&&\qquad\quad{}
+\E f(W^*(t))-f(0).\nonumber
\end{eqnarray}
Taking $f(x)=x$, we obtain
$\E[I_S^*(t)+K_-^*(t)]-\E K_+^*(t)=\gamma t+\E W^*(t)-f(0)$.
If $\gamma\neq0$, we may take
$f(x)=\frac{\sigma^2}{2\gamma}e^{2\gamma x/\sigma^2}$
in (\ref{feq}), which leads to the equation
$\E[I_S^*(t)+K_-^*(t)]
-e^{2\gamma H(0)/\sigma^2}\E K_+^*(t)=
\frac{\sigma^2}{2\gamma}(\E e^{2\gamma W^*(t)/\sigma^2}-1)$.
Solving these equations for $\E K_+^*(t)$,
we obtain the second equality in (\ref{mstar}) for
$\gamma\neq0$. To obtain this equality for
$\gamma=0$, we take $f(x)=x^2$.
\end{pf}

According to\vspace*{1pt} (\ref{3.10}), the work that arrives to
the $n$th system by time $nt$ is
$V^{(n)}(A^{(n)}(nt))=\sqrt{n}\widehat{N}^{(n)}(t)+nt$.
But,
$\widehat{N}^{(n)}$
is approximately $N^*$, and hence
\[
\lim_{t\rightarrow\infty}
\frac{\sqrt{n}\widehat{N}^{(n)}(t)+nt}{nt}
\approx
\lim_{t\rightarrow\infty}
\frac{\sqrt{n}N^*(t)+nt}{nt}
=\biggl(1-\frac{\gamma}{\sqrt{n}}\biggr).
\]
Therefore, if $\gamma\neq0$,
the long-run fraction of reneged work is approximately
\begin{eqnarray*}
\lim_{t\rightarrow\infty}
\frac{\sqrt{n} K_+^*(t)}{V^{(n)}(A^{(n)}(nt))}
&=&
\frac{1}{\sqrt{n}}\lim_{t\rightarrow\infty}
\frac{1}{t} K_+^*(t)\cdot
\lim_{t\rightarrow\infty}
\biggl(\frac{\sqrt{n}\widehat{N}^{(n)}(t)+nt}{nt}\biggr)^{-1}\\
&\approx&
\frac{\gamma/\sqrt{n}}
{(1-\gamma/\sqrt{n})(e^{2\gamma H(0)/\sigma^2}-1)}.
\end{eqnarray*}
Finally, (\ref{2.4}) implies that
the expected lead time in the $n$th system is
$\E L_j^{(n)}=\int_{0}^{\infty}(1-G(y/\sqrt{n}))\,dy
=\sqrt{n} H(0)$.
Using this formula and (\ref{3.3}),
we conclude that the fraction of work that reneges
in the $n$th system when $\gamma\neq0$ is approximately
%
\begin{equation}\label{LostWork}
\frac{1-\rho^{(n)}}
{\rho^{(n)}(e^{2(1-\rho^{(n)})\E L_j^{(n)}/\sigma^2}-1)}
=\frac{1-\rho^{(n)}}{\rho^{(n)}(e^{\theta\Dbar}-1)},
\end{equation}
where
\[
\theta=\frac{2(1-\rho^{(n)})}{\sigma^2}
\approx\frac{2\gamma}{\sqrt{n}\sigma^2},\qquad
\Dbar=\E L^{(n)}_j=\sqrt{n}H(0).
\]
We have suppressed the dependence of $\theta$
and $\Dbar$ on $n$, which will remain fixed.
If $\gamma=0$, then in place of
(\ref{LostWork}) we have
$\frac{\sigma^2}{2\Dbar}$.
We have established (\ref{frw}) and (\ref{frwo}).
\begin{remark}\label{RQ}
Corollary \ref{C3.5a} also implies that the limiting
scaled queue length process is
$\lambda W^*$, which is a doubly reflected Brownian
motion in $[0,\lambda H(0)]$ with drift $-\gamma\lambda$
and variance per unit time $\lambda^2\sigma^2$.
This incorrectly suggests that $\lambda\sqrt{n} K_+^*(t)$
is approximately the number of
customers who renege in the
$n$th system prior to~$nt$. The simulations indicate
that this naive interpretation of Corollary~\ref{C3.5a}
applied to the queue length process is incorrect, as does
the following heuristic.

According to \cite{harrison}, Proposition 5, page 90,
if $\gamma\neq0$, the stationary density for $W^*$ is
%
\begin{equation}\label{7.5a}
\varphi^*(x)\define
\cases{
\dfrac{2\gamma e^{-2\gamma x/\sigma^2}}
{\sigma^2(1-e^{-2\gamma H(0)/\sigma^2})}, &\quad if $0\leq x\leq
H(0)$,\vspace*{2pt}\cr
0, &\quad otherwise,}
\end{equation}
whereas the stationary density is uniform on $[0,H(0)]$
if $\gamma=0$.
Therefore, for $\gamma\neq0$ and $t$ large,
the density of $W^{(n)}(nt)\approx\sqrt{n}W^*(t)$
is approximately
\[
\varphi(w)=\frac{1}{\sqrt{n}}\varphi^*\bigl(w/\sqrt{n}\bigr)
=\cases{
\dfrac{\theta e^{-\theta w}}{1-e^{-\theta\Dbar}},
&\quad if $0\leq w\leq\Dbar$,\vspace*{2pt}\cr
0, &\quad otherwise.}
\]
We have suppressed the dependence of
$\varphi$ on $n$.

Suppose now that the lead times of arriving
customers are not random. Then in the $n$th system,
all lead times are equal to $\sqrt{n}H(0)=\Dbar$.
In this case, the EDF policy serves customers in
order of arrival (FIFO).
Suppose the workload in queue is $W$
at the time of arrival of a customer whose service requirement
is $V$. Recall that the expected service
time is $1/\mu^{(n)}$, and because $n$ is fixed, we suppress
it and write $\E V=1/\mu$.
The arriving\vspace*{1pt} customer will be served to completion if
and only if $W+V\leq\Dbar$.
Suppose further that the arrival process
$A^{(n)}$ is Poisson, so that according to
the PASTA property (``{\underline P}oisson
{\underline a}rrivals {\underline s}ee
{\underline t}ime {\underline a}verages'';
see \cite{asm}, Theorem 6.7, page 218),
an arriving customer will encounter a workload
$W$ having approximately the distribution $\varphi$.
The probability the arriving customer eventually reneges is thus
\[
\PP\{W>\Dbar-V\}=
\E[\PP\{W>\Dbar-V|V\}]=
\E\biggl[\int_{(\Dbar-V)^+}^{\Dbar}\varphi(w) \,dw\biggr].
\]
Because $\Dbar$ is of order $\sqrt{n}$ and $V$ is of order $1$,
we have $(\Dbar-V)^+=\Dbar-V$ with high probability.
Using this approximation, we complete the calculation
for the case $\gamma\neq0$ to obtain
%
\begin{equation}\label{PRen}
\PP\{\mbox{Customer reneges}\}
\approx\frac{1}{e^{\theta\Dbar}-1}(\E e^{\theta V}-1).
\end{equation}

If the customer reneges, then
work $V+W-\Dbar>0$ is lost.
The expected lost work is
\[
\E[V+W-\Dbar|\mbox{Customer reneges}]\approx
\E\biggl[\frac{\int_{(\Dbar-V)^+}^{\Dbar}(V+w-\Dbar)
\varphi(w) \,dw}
{\PP\{\mbox{Customer reneges}\}}\biggr].
\]
Again using the approximation $(\Dbar-V)^+\approx\Dbar-V$,
we obtain
\begin{eqnarray*}
\E[V+W-\Dbar|\mbox{Customer reneges}]
&\approx&
\frac{1}{\theta}-\frac{\E V}{\E e^{\theta V}-1}\\
&\approx&
\frac{1}{\theta}-\frac{\E V}{\theta\E V+1/2\theta^2\E[V^2]
+O(n^{-3/2})}\\
&\approx&
\frac{\E[V^2]}{2\E V}.
\end{eqnarray*}
The last expression is, perhaps not surprisingly,
the formula for the average residual
lifetime of a renewal cycle (see \cite{Ross}, Example 3.6(b),
pages 80 and 81).
Consequently, when lead times are constant
and the arrival process is Poisson, we should expect the total number of
customers reneging in $[0,t]$ times the expected amount of work
lost per reneging customer to approximately equal the total amount of
work lost by reneging in $[0,t]$. If we divide both by the total number of
customer arrivals in $[0,t]$ and take limits as
$t \rightarrow\infty$, we find
%
\begin{eqnarray}\label{7.6}
&&\mbox{Fraction of lost customers in reneging system}
\nonumber\\
&&\qquad\approx
\frac{\mbox{Fraction of lost work in reneging system}}
{\E[V+W-\Dbar|\mbox{Customer reneges}]}*\E V\\
&&\qquad\approx
\frac{2(\E V)^{2}}{\E[V^2]}
\times
(\mbox{Fraction of lost work in reneging system}).\nonumber
\end{eqnarray}
This is (\ref{flcr}) with $\E V=\frac{1}{\mu}$ and
$E[V^2]=\beta^2+\frac{1}{\mu^2}$.\vspace*{1pt}

If $V$ is exponentially
distributed, hence $\E[V^{2}] = 2(\E V )^2$, then (\ref{7.6}) implies that
the fraction of customers who renege will
be approximately the fraction of work that reneges.
See Figure \ref{fig:4} for simulations that confirm this
assertion.
On the other hand, if $V$ is nonrandom, hence equal to
its mean $1/\mu$, then (\ref{7.6}) predicts that
the fraction of customers who renege will be
twice the fraction of work that reneges.
See Figure \ref{fig:5} for simulations that confirm
this assertion.
Both these conclusions hold irrespective of the value of $\lambda$.

\begin{figure}

\includegraphics{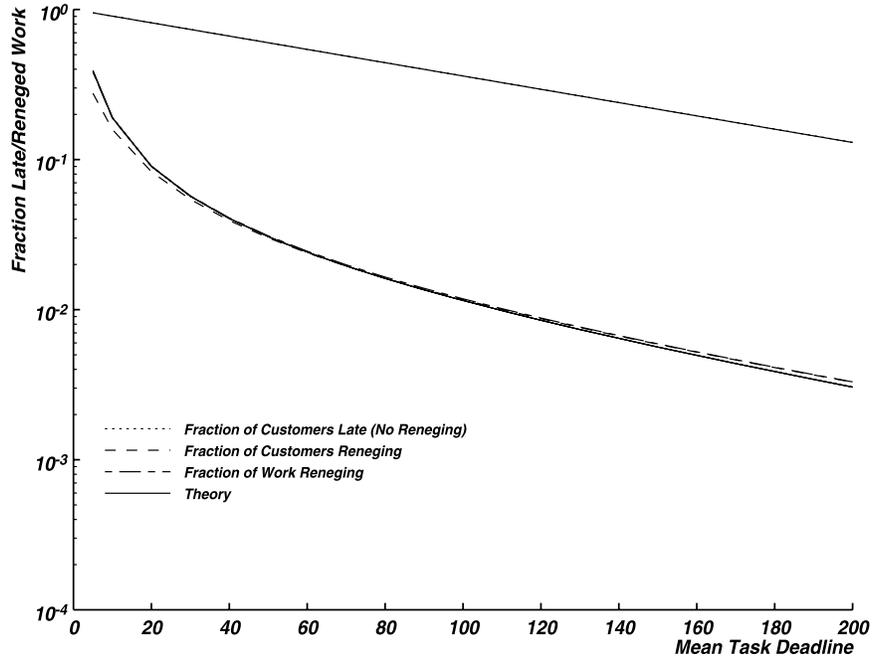}

\caption{$M/M/1$ queue.}
\label{fig:4}
\end{figure}

\begin{figure}

\includegraphics{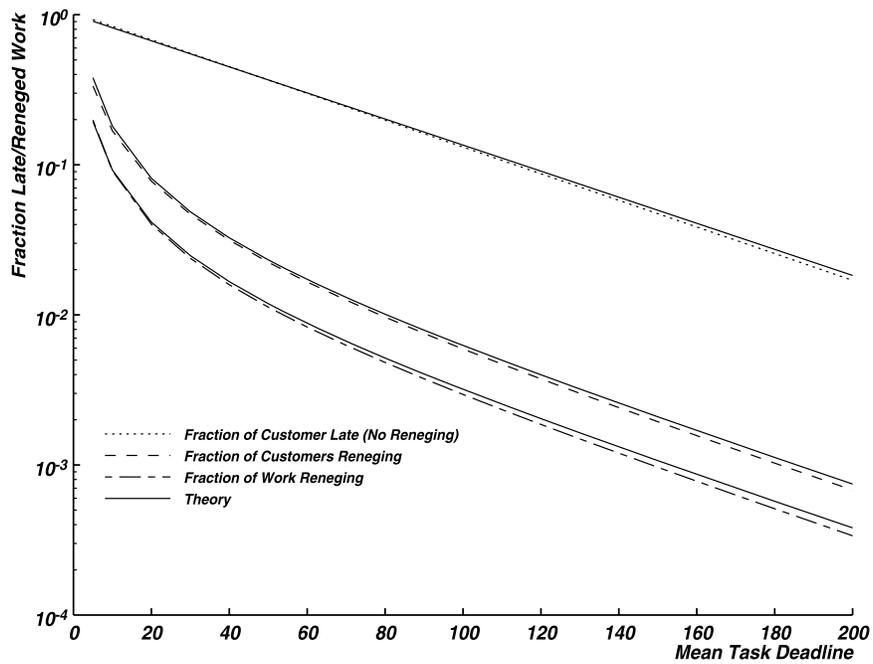}

\caption{$M/D/1$ queue.}
\label{fig:5}
\end{figure}

The last conclusion is inconsistent with a naive
interpretation
of Corollary \ref{C3.5a}, according to which
work reneges at a rate $1/\lambda$ times the
rate of customer reneging. Since work arrives at a rate
$\E V\approx1/\lambda$ times the rate of customer arrivals,
this naive interpretation of Corollary \ref{C3.5a} would
say that the fraction of work reneging would
approximately agree with the fraction of customers reneging
regardless of the distribution of $V$.
\end{remark}

We next turn our attention to the performance of
the standard (nonreneging) system.
Recall from (\ref{w*}) that
the scaled workload process when all customers are served
to completion converges to $W_{S}^{*}$, a reflected Brownian
motion with drift $-\gamma$ (we now assume $\gamma>0$
in order to have a stationary distribution)
and variance~$\sigma^{2}$.
In particular, $W^{(n)}_S(nt)\approx\sqrt{n}W^*_S(t)$.
The stationary density for $W^*_S$ is
\[
\varphi^*_S(x)\define\cases{
\dfrac{2\gamma}{\sigma^2}e^{-2\gamma x/\sigma^2}, &\quad
if $x\geq0$,\vspace*{2pt}\cr
0, &\quad otherwise,}
\]
and so for large $t$, the density of $W^{(n)}(nt)$
is approximately
\[
\varphi_S(w)=\frac{1}{\sqrt{n}}\varphi^*_S\bigl(w/\sqrt{n}\bigr)
=\cases{
\theta e^{-\theta w}, &\quad if $w\geq0$,\cr
0, &\quad otherwise.}
\]
Consequently, the long-run fraction of time $W^{(n)}$
spends above level $\Dbar$ is $e^{-\theta\Dbar}$.
The workload level at which
the limiting frontier reaches $0$ is $H(0)$, and hence
it is approximately the case that the
$n$th system sees lateness if and only if
$W^{(n)}$ exceeds $\Dbar=\sqrt{n}H(0)$.
In other words, the theory predicts that
%
\begin{eqnarray} \label{lwcs2}
&&\mbox{Fraction of late customers in standard system}
\nonumber\\
&&\qquad=
\mbox{Fraction of late work in standard system}\\
&&\qquad=
e^{-\theta\Dbar}.\nonumber
\end{eqnarray}
We are using here the result for $GI/G/1$ queues that
\begin{eqnarray*}
&&\lim_{n \rightarrow\infty}\lim_{T \rightarrow
\infty}\frac{1}{T}\int_{0}^{T}\ind_{\{\widehat{W}_{S}^{(n)}(t) >
H(0)\}}\,dt \\
&&\qquad = \lim_{T \rightarrow\infty}\lim_{n \rightarrow
\infty}\frac{1}{T}\int_{0}^{T}\ind_{\{\widehat{W}_{S}^{(n)}(t) >
H(0)\}}\,dt \\
&&\qquad = \lim_{T \rightarrow\infty}
\frac{1}{T}\int_{0}^{T}\ind_{\{W_{S}^{*}(t) > H(0)\}}\,dt,
\end{eqnarray*}
a result that grows out of the work of Kingman
\cite{Kingman61,Kingman62} (see \cite{GamarnikZeevi}
for a general result that specializes to the
case under consideration).

It is important to compare the fraction of work that
reneges in the reneging system, given by (\ref{LostWork}),
with the fraction of work that is
late in the standard (nonreneging) system.
The ratio of these quantities of lost/late work is
%
\begin{equation}\label{Ratio}
\frac{\mbox{Lost work in reneging system}}
{\mbox{Late work in standard system}}
\approx\frac{e^{\theta\Dbar}}{e^{\theta\Dbar}-1}
\biggl(\frac{1-\rho^{(n)}}{\rho^{(n)}}\biggr).
\end{equation}
The parameter $\theta$ is $O(1/\sqrt{n})$, $\theta\Dbar$ is
$O(1)$, and $1-\rho^{(n)}$ is $O(1/\sqrt{n})$.
Thus the ratio in (\ref{Ratio}) is $O(1/\sqrt{n})$.
\begin{remark}\label{R7.3}
If lead times are a nonrandom
constant $\Dbar$, EDF reduces to first-in-first-out,
and the fraction of lost customers in an M/G/1 queue
with $0<\rho<1$ is
$(1-\rho)\PP\{W>\Dbar\}/(1-\rho\PP\{W>\Dbar\})$,
where $W$ is the steady-state workload in the corresponding
nonreneging M/G/1 queue
(see \cite{BT}). In the heavy traffic limit of our model,
$\PP\{W>\Dbar\}=e^{-\theta\Dbar}$ [see the derivation of
(\ref{lwcs2})].
Recalling that $1-\rho=O(1/\sqrt{n})$
in (\ref{frw}), we observe that this is consistent with (\ref{frw}).
\end{remark}

\subsection{Simulation results}\label{Subsection7.2}

We conducted a simulation study to assess the accuracy of these
approximations and to compare the performance of the systems with and
without reneging. Two systems were considered, an $M/M/1$ system
presented in Figure \ref{fig:4} and an
M/D/1 system presented in Figure \ref{fig:5}.
In both cases, $\lambda= 0.5$
and $\frac{1}{\mu} =
1.96$, and
so the traffic intensity is $\rho= 0.98$. These parameter values
result in $\theta= 0.010202$ for the $M/M/1$ case and $\theta= 0.02$ for
the M/D/1 case. The initial deadline
distribution is uniform
on $[5,B]$ with the mean deadline
$\Dbar= \frac{5+B}{2}$, varying from $B=5$
(constant deadlines)
to $B=200$. The data points
are the simulation results averaged over one billion customer arrivals
per case. The curves that are superimposed on the data are the
theoretical values, $e^{-\theta\Dbar}$ for the case in which
customers are served to completion (the standard system), and equations
(\ref{frw}) and (\ref{flcr}) for the fraction of work lost and the
fraction of lost customers for the reneging system. Equation
(\ref{flcr}) is derived in
Remark \ref{RQ}
under the assumption of constant deadlines.
Nevertheless, we apply it for the variable deadline case in the
simulation study. The fraction of late work or late customers for the
system in which customers are served to completion is also presented
to compare its performance with that of the reneging system.

The $M/M/1$ results are presented in
Figure \ref{fig:4}
with the fraction of customers missing their deadlines,
the fraction of customers reneging, and the fraction of work
reneging plotted
on a log scale on the $y$-axis against the mean deadline on the $x$-axis.
There is nearly perfect
agreement between the theoretical
approximation and the simulation. In fact,
one cannot see the plot of ``Fraction of Customers Late (No Reneging)''
because it coincides with the ``Theory'' plot
at the top of the figure. Similarly, one can see only
parts of the plots of ``Fraction of Customers Reneging'' and
``Fraction of Work Reneging'' because they coincide
with the ``Theory'' plot in the middle of the figure.
One can see the linear form
for the case of service to completion. Furthermore,
the simulation confirms the prediction of (\ref{frw})--(\ref{7.3})
that for sufficiently
large values of $\Dbar$, the performance
of the reneging system is parallel on a log scale to that of the
standard system with the two curves separated by approximately $0.02$.
This corresponds to a
reduction in work that misses its deadline by a factor of $40$ to $50$.

Figure \ref{fig:5} presents the
results for the M/D/1 system. The results are
qualitatively identical to those of Figure
\ref{fig:4}, except the fits of the
theoretical curves are not as exact as the fits for the $M/M/1$ system;
it appears that now the value $\theta= 0.02$ is
slightly too small and hence the theory slightly
overestimates the fraction of work that misses its deadline,
especially when the mean deadline is large.
Also, the lost or late work and the
customer loss or lateness fractions
are significantly smaller than for the $M/M/1$
system owing to the reduction in variability of the customer
service time distribution. The reduction in missed deadlines between
the two systems for large values of $\Dbar$ is again a factor of 40
to 50. In both figures, it is clear that there are one to two orders
of magnitude of improvement in the overall performance of the system
resulting from stopping service on customers when their deadlines expire.

\begin{appendix}\label{app}
\section*{Appendix: Optimality of EDF}

\vspace*{-8pt}

\begin{pf*}{Proof of Theorem \ref{t.EDFopt}}
Let $\pi$ be a service policy and let $t_0$ be the first time
$\pi$ deviates from the EDF policy, either because it idles
when there is work present, or it serves a customer other
than the customer present with the smallest lead time.
Let $j$ be the index of the customer with the smallest
lead time at time $t_0$.

We consider first the case that $\pi$ idles at
time $t_0$. In this case, we
define $\rho(\pi)$ to be the policy that emulates $\pi$ except
as noted below.
From time $t_0$, whenever $\pi$ idles, $\rho(\pi)$
serves customer $j$,
at least until time $t_1$, when customer $j$
leaves the $\rho(\pi)$ system because either
$\rho(\pi)$ serves customer $j$
to completion or else the deadline of customer $j$ elapses.
From time $t_1$, $\rho(\pi)$ idles
if $\pi$ serves customer $j$.
We will show that
for $t\geq0$,
\setcounter{equation}{0}
\begin{equation}\label{eqA.1}
R_{\rho(\pi)}(t)\leq R_\pi(t).
\end{equation}
Let $v_k(t)$ [resp., $v^{\rho}_k(t)$] be the residual service
time of the $k$th customer at time $t$ under $\pi$
[resp., $\rho(\pi)$].
In particular, if $d_k$ is the deadline of the $k$th customer,
then $v_k(d_k-)$ [resp., $v^{\rho}_k(d_k-)$]
is the work corresponding
to this customer
that is deleted by $\pi$ [resp., $\rho(\pi)$]
due to lateness, and
%
\begin{equation}\label{eqA.2}
R_{\rho(\pi)}(t) = \sum_{k\dvtx d_k \leq t} v^{\rho}_k(d_k-),\qquad
R_\pi(t) = \sum_{k\dvtx d_k \leq t} v_k(d_k-).
\end{equation}
By the definition of $\rho(\pi)$, for $t\geq0$
and $k\neq j$, we have
%
\begin{equation}\label{eqA.3}
v^{\rho}_k(t)= v_k(t),
\end{equation}
whereas
%
\begin{equation}\label{eqA.4}
v_j^{\rho}(t)\leq v_j(t).
\end{equation}
Summing (\ref{eqA.3}) over $k\neq j$, invoking (\ref{eqA.4})
and (\ref{eqA.2}),
we obtain (\ref{eqA.1}).

We next consider the case that at time $t_0$, $\pi$ serves
customer $i\neq j$. In this case, we define $\rho(\pi)$
to be the policy
that emulates $\pi$ except as noted below.
From time~$t_0$, whenever
$\pi$ serves customer $i$, $\rho(\pi)$ serves customer $j$, at least
until time~$t_1$, when $\rho(\pi)$ serves customer $j$ to completion
or the deadline of customer $j$ elapses. From time $t_1$,
$\rho(\pi)$ serves customer $i$ if $\pi$ serves customer $j$,
provided customer $i$ is present in the system under $\rho(\pi)$.
If $\pi$ serves customer $j$ and customer $i$ is not present
under $\rho(\pi)$, then $\rho(\pi)$ idles.
We again have (\ref{eqA.2}) and (\ref{eqA.4}),
whereas (\ref{eqA.3}) now holds only for $k\notin\{i,j\}$.
If the $i$th customer is served to completion under
$\rho(\pi)$, then $v^{\rho}_i(d_i-)=0$, and
(\ref{eqA.3}) for $k\notin\{i,j\}$, and (\ref{eqA.4})
imply that (\ref{eqA.1}) holds for all $t$. It remains
to consider the case that the
$i$th customer becomes late under $\rho(\pi)$. In this case
(\ref{eqA.3}) for $k\notin\{i,j\}$ and (\ref{eqA.4})
imply that (\ref{eqA.1}) holds for $t\in[0,d_i)$.
Let $w_1$ denote the work done by $\rho(\pi)$ on the
$j$th customer when $\pi$ works on the $i$th
customer in the interval $[t_0,t_1)$.
Let $w_2$ be the work done by $\rho(\pi)$
on customer $i$ in the time interval $[t_1,\infty)$
while $\pi$ works on
customer $j$ in this time interval. Finally, let $w_3$
be the work done by $\pi$ on customer $j$
in the time interval $[t_1,\infty)$ while $\rho(\pi)$ is idle.
Then $v_j^{\rho}(d_j-)+w_1=v_j(d_j-)+w_2+w_3$ and
$v_i^{\rho}(d_i-)+w_2=v_i(d_i-)+w_1$, which implies
%
\begin{equation}\label{eqA.5}
v_j^{\rho}(d_j-)+v_i^{\rho}(d_i-)=v_j(d_j-)+v_i(d_i-)+w_3.
\end{equation}
We argue by contradiction that $w_3$ cannot be positive.
If $w_3$ were positive, then at some time $t\geq t_1$,
$\pi$ serves customer $j$ and customer $i$ is not in the
$\rho(\pi)$ system. This implies that
$d_j>t$, and since by assumption,
$d_i>d_j$, the absence of customer $i$ in the $\rho(\pi)$ system
means that this system has served customer $i$ to completion.
We conclude that $v_i^{\rho}(d_i-)=0$. On the other hand,
customer $j$ is also not in the $\rho(\pi)$ system at time
$t\geq t_1$,
and so $v_j^{\rho}(d_j-)=0$ as well.
The left-hand side of (\ref{eqA.5})
is zero, and hence $w_3$ must be zero. We conclude that
%
\begin{equation}\label{eqA.6}
v_j^{\rho}(d_j-)+v_i^{\rho}(d_i-)
=v_j(d_j-)+v_i(d_i-).
\end{equation}
Since $d_j<d_i$, if $t\geq d_i$,
then (\ref{eqA.3}) for $k\notin\{i,j\}$ and
(\ref{eqA.6}) imply (\ref{eqA.1}).

Starting from the service policy $\pi$, we
have obtained a service policy $\rho(\pi)$
that either is work conserving until the departure of
customer $j$ or else gives customer $j$ priority
over customer $i$ until the departure of customer $j$.
However, immediately after time $t_0$, the policy $\pi$
may serve some customer $k\notin\{i,j\}$, and hence $\rho(\pi)$
also serves $k$ at this time, although customer
$j$ is more urgent. Therefore, we apply $n$ iterations of
the mapping $\rho$, where $n$ is the number
of customers in the $\pi$ system at time $t_0$, and thereby
obtain a policy that is work-conserving
and serves in EDF order at least until the first time after
$t_0$ that there is a departure or an arrival. We have
$R_{\rho^n(\pi)}(t)\leq R_{\pi}(t)$ for all $t\geq0$.

By assumption, for each $t$ the number of system arrivals $A(t)$
by time $t$ is finite. Hence the maximum
number of customers in the system over the interval $[0,t]$
is bounded by $A(t)$, and the number of arrivals and departures
up to time $t$ is bounded by $2A(t)$, irrespective of
the service policy.
Thus, if we start with any
policy $\pi$, the number of iterations of the mapping $\rho$
required to obtain a policy that is work conserving
and serves in EDF order up to time $t$ is finite.
Under this policy the amount of work removed by lateness up
to time $t$ is the same as for the EDF system
in the theorem, and hence (\ref{e.EDFopt}) holds.
\end{pf*}
\setcounter{theorem}{0}
\begin{remark}
In the above proof we have implicitly assumed that
$\pi$ [and thus $\rho(\pi)$] never serves more
than one customer at the same time.
This assumption simplifies the exposition of
the argument, and the
generality of Theorem \ref{t.EDFopt}
is sufficient for
this paper. However, the proof can be generalized
to policies permitting simultaneous
service of customers (e.g., processor sharing).
In this case, in the construction of
$\rho(\pi)$ we must additionally take the
rates at which customers receive
service into account. For example, the difference in the rates
with which the $j$th customer
receives service under $\rho(\pi)$ and $\pi$ in the
time interval $[t_0,t_1)$ must be equal
to the rate of service of
the $i$th customer under $\pi$ in
this time interval, the rates of service
of all other customers in this time interval
under $\pi$ and $\rho(\pi)$ must be the same, etc.
\end{remark}
\end{appendix}


\printaddresses

\end{document}